# A Novel Meshless Method Based on the Virtual Construction of Node Control Domains for Porous Flow Problems

Xiang Rao[1, 2, 3]*
[1]School of Petroleum Engineering, Yangtze University, Wuhan 430100, P. R. China
[2]Key Laboratory of Drilling and Production Engineering for Oil and Gas, Hubei Province, Wuhan Hubei China, 430100, P. R. China
[3]Cooperative Innovation Center of Unconventional Oil and Gas (Ministry of Education & Hubei Province), Yangtze University, Wuhan, 430100, P. R. China
*Corresponding author: raoxiang0103@163.com, raoxiang@yangtzeu.edu.cn

**Abstract:** In this paper, a novel meshless method that can handle porous flow problems with singular source terms is developed by virtually constructing the node control domains. By defining the connectable node cloud, this novel meshless method uses the integral of the diffusion term and generalized difference operators to derive overdetermined equations of the node control volumes. An empirical method of calculating reliable node control volumes and a triangulation-based method to determine the connectable point cloud are developed. NCDMM only focuses on the volume of the node control domain rather than the specific shape, so the construction of node control domains is called virtual, which will not increase the computational cost. To our knowledge, this is the first time to construct node control volumes in the meshless framework, so this novel method is named a node control domains-based meshless method, abbreviated as NCDMM, which can also be regarded as an extended finite volume method (EFVM). Taking two-phase porous flow problems as an example, the NCDMM discrete schemes meeting local mass conservation are derived by integrating the generalized finite difference schemes of governing equations on each node control domain. Finally, existing commonly used low-order finite volume method (FVM) based nonlinear solvers for various porous flow models can be directly employed in the proposed NCDMM, significantly facilitating the general-purpose applications of the NCDMM. Theoretically, the proposed NCDMM has the advantages of previous meshless methods for discretizing computational domains with complex geometries, as well as the advantages of traditional low-order FVMs for stably handling a variety of porous flow problems with local mass conservation. Four numerical cases are implemented to test the computational accuracy, efficiency, convergence, and good adaptability to the calculation domain with complex geometry and various boundary conditions.

## 1. Introduction

Up to now, numerical methods are mainly divided into mesh-based and meshless types according to the discretization of the computational domain. For reservoir porous flow problems, due to the large size of the computational domain high-resolution methods are rarely applied, while the mesh-based low-order FVM is the most widely used. The significant factor limiting the mesh-based numerical methods is to generate a high-quality mesh to discretize the computational domain with complex geometry. Under this background, various meshless methods are proposed to obtain the high-accuracy calculation by only utilizing a node cloud instead of the mesh to discretize the calculation domain. However, due to the lack of the node control volume in previous meshless methods, less attention is paid to problems with singular source terms which is very common in reservoir porous flow problems, and the important local mass conservation in fluid flow simulation, so previous meshless methods still do not occupy a sufficient share in flow problems, especially in reservoir numerical simulation.

The generalized finite difference method (GFDM) is a popular domain-type meshless method [1, 2, 3]. This method uses a node cloud to discretize the computational domain and divides it into intersecting node influence subdomains. The spatial derivatives of functions in the governing equation can be expressed as the difference scheme of nodal values using the local Taylor expansion of multivariate functions and the weighted least square method in the subdomain, overcoming the mesh dependence of the traditional finite difference method (FDM). Up to now, this method has been widely used to solve seismic wave propagation problem [4], Burgers' equations [5, 6], shallow water equations [7, 8], problems of plates [9], heat conduction analysis [10, 11], elasticity problems [12, 13, 14], water wave interactions problems [15], soil mechanics [16, 17], flow problems [18, 19, 20, 21, 22], liquid-vapor phase transitions [23], and various scientific and engineering problems [24, 25, 26]. Rao et al. [27, 28] introduced the single point upstream (SPU) scheme [29, 30] in reservoir simulators to handle the convection term in the flow governing equations with a stable upwind effect and applied the upwind GFDM to solve the heat and mass transfer and two-phase porous flow problems. The good computational performances indicate that GFDM may realize high-quality meshless solutions to porous flow problems, and has great potential in reservoir numerical simulation.

For the reservoir flow problems, there are often point source terms (e.g. production and injection wells) in the reservoir domain [31], because the radius of the drilling wells is very small compared with the size of the reservoir domain. In the traditional FDM for reservoir numerical simulation [32, 33], the singular source term can be approximated by the nonsingular average source strength in the grid volume which contains the point source. The singular Dirac function representing the point sources is integrated to become nonsingular in the finite element method (FEM) [34, 35] and FVM [36, 37, 38, 39], thus the singular source or sink term is handled naturally. The traditional FVM is most widely used in reservoir numerical simulation, including the commonly used reservoir numerical simulation platforms CMG, ECLIPSE, AD-GPRS [40, 41], Matlab Reservoir Simulator Toolbox (MRST) [42] because it can adapt to wider mesh types than FDM, meet local mass conservation, and have simpler implementation for multiphase/multicomponent flow compared with FEM. However, FVM can directly handle closed boundary conditions, but it struggles with more complex boundary conditions. For the computational domain with complex geometry (complex boundary geometry, fractures, caves, flow barriers, etc.), the high-quality mesh division is often difficult.

GFDM can yield generalized finite difference schemes by using local Taylor expansion and moving least squares when the calculation domain is only discretized by a point cloud, which motivates us to develop a meshless reservoir numerical simulator that can flexibly discretize the calculation domain with complex geometry and directly handle various boundary conditions. However, because the node cloud in previous meshless methods (including GFDM) lacks geometric information like node control volume, making it is difficult to efficiently deal with singular source terms (e.g. production and injection wells) and employ the existing nonlinear solver of various porous flow models.

Therefore, by introducing node control volumes to the meshless point cloud, this study aims to develop a node control domain-based meshless method (NCDMM) that can handle singular source terms and satisfy local mass conservation. Theoretically, the proposed NCDMM has the advantages of previous meshless methods for discretizing complex domains, as well as the advantages of traditional FVM for stably handling a variety of porous flow problems. Overall, the proposed NCDMM is expected to become a new mainstream numerical method and provide a novel meshless tool to handle various scientific and engineering problems, especially reservoir flow problems.

This paper is organized as follows. Section 2.1 briefly reviews the basic GFDM theory. Section 2.2

presents the virtual construction method of node control domains, including the definitions of the node control domain and the connectable point cloud, and derivation of the calculation equations of node control volumes. Section 2.3 develops an empirical method of calculating node control volumes, and Section 2.4 gives a triangulation-based method of determining the connectable point cloud. Section 2.5 derives the NCDMM discrete scheme by integrating the GFDM discrete scheme on the node control domain and shows that NCDMM satisfies the local mass conservation. Section 2.6 introduces the treatment of various boundary conditions. Section 2.7 conducts a convergence analysis of NCDMM. In Section 3, four numerical test cases are implemented to analyze the computational performances of the proposed NCDMM. The main conclusions are summarized in Section 4. Some future work is given in Section 5.

## 2. Methodology

### 2.1 A brief review of GFDM

GFDM is a meshless method, that only uses a point cloud to discretize the computational domain. In GFDM, a node influence domain for each node is defined, and the influence domains of different nodes may intersect, i.e., $\Omega = \bigcup_i I_i$, where, $\Omega$ is the calculation domain and $I_i$ is the influence domain of node *i*.

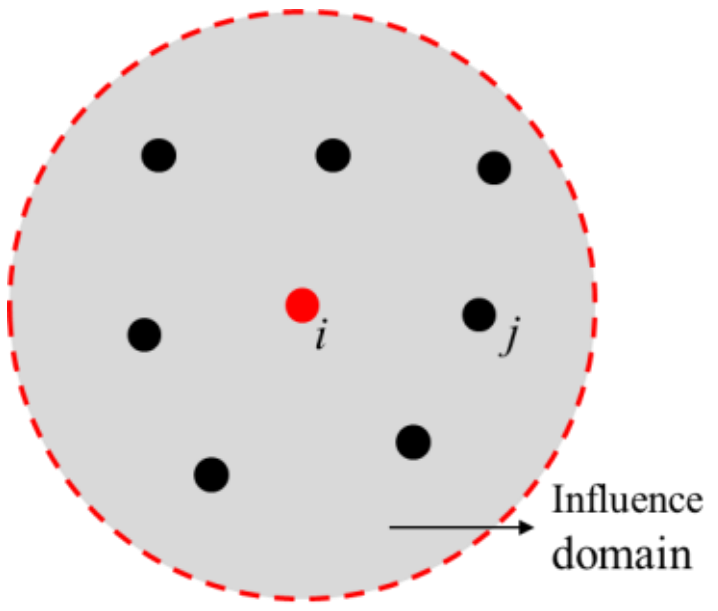

Fig. 1 A sketch of node influence domain and containing nodes

As shown in Fig. 1, let node *i* be a node in $\Omega$, and the coordinate is $(x_i, y_i)$. $I_i$ contains another $n_i$ nodes except node *i*, which are denoted as $j$ ($j = 1, 2, 3 \cdots n_i$), and corresponding coordinates are $(x_j, y_j)$ ($j = 1, 2, 3 \cdots n_i$). Using Taylor expansion of the unknown function at node *i*, it is obtained that:

$$u(x_j, y_j) = u(x_i, y_i) + \Delta x_j \left.\frac{\partial u}{\partial x}\right|_i + \Delta y_j \left.\frac{\partial u}{\partial y}\right|_i + \frac{1}{2}\left( (\Delta x_j)^2 \left.\frac{\partial^2 u}{\partial x^2}\right|_i + 2\Delta x_j \Delta y_j \left.\frac{\partial^2 u}{\partial x \partial y}\right|_i + (\Delta y_j)^2 \left.\frac{\partial^2 u}{\partial y^2}\right|_i \right) + O(r^3) \quad (1)$$

where $\Delta x_j = x_j - x_i$, $\Delta y_j = y_j - y_i$.

Simply denote $u_i = (x_i, y_i)$, $u_{x,i} = \left.\frac{\partial u}{\partial x}\right|_i$, $u_{y,i} = \left.\frac{\partial u}{\partial y}\right|_i$, $u_{xx,i} = \left.\frac{\partial^2 u}{\partial x^2}\right|_i$, $u_{xy,i} = \left.\frac{\partial^2 u}{\partial x \partial y}\right|_i$, $u_{yy,i} = \left.\frac{\partial^2 u}{\partial y^2}\right|_i$.

Define the weighted error function $B(\mathbf{D}_u)$:

$$B(\mathbf{D}_u) = \sum_{j=1}^{n} \left[ \left( u_i - u_j + \Delta x_j u_{x,i} + \Delta y_j u_{y,i} + \frac{1}{2}(\Delta x_j)^2 u_{xx,i} + \frac{1}{2}(\Delta y_j)^2 u_{yy,i} + \Delta x_j \Delta y_j u_{xy,i} \right) \omega_j \right]^2 \quad (2)$$

where $\mathbf{D}_u = (u_{x,i}, u_{y,i}, u_{xx,i}, u_{yy,i}, u_{xy,i})^T$, $\omega_j = \omega(\Delta x_j, \Delta y_j)$ is the weight function. Commonly used weight functions include exponential functions, potential functions, cubic splines, quartic splines, etc. [43]. Benito et al., [1], Gavete et al., [3], and Chen et al. [19] demonstrated that different weight functions have little impact on the GFDM calculation results, and currently quartic spline function in Eq. (3) is often used. However, in this work, it is found that the weight function in Eq. (4) can make the values of node control volumes more accurate in most cases. In fact, the weight function in Eq. (4) was widely used in the early GFDM work [3], and this weight function can yield the well-known nine-point finite difference scheme in GFDM, which makes it stand out in so many weight functions. Section 2.3 will give an Eq. (4)-based empirical method of calculating node control volumes, and numerical examples in Section 3 will illustrate that compared with using the weight function in Eq. (3), using the weight function in Eq. (4) can obtain higher-accuracy node control volumes and simulation results.

$$\omega_j = \begin{cases} 1-6\left(\frac{r_j}{r_m}\right)^2+8\left(\frac{r_j}{r_m}\right)^3-3\left(\frac{r_j}{r_m}\right)^4 & r_j \le r_m \\ 0 & r_j > r_m \end{cases} \tag{3}$$

$$\omega_j = \begin{cases} 1\Big/\left(\frac{r_j}{r_m}\right)^3 & r_j \le r_m \\ 0 & r_j \le r_m \end{cases} \tag{4}$$

where $r_j$ is the Euclidean distance from node *j* to the central node *i*, and $r_m$ is the radius of the influence domain of node *i*.

To minimize the weighted error function, the partial derivatives of the weighted error function are required equal to zero, that is,

$$\frac{\partial B(\mathbf{D}_u)}{\partial u_{x,i}}=0,\frac{\partial B(\mathbf{D}_u)}{\partial u_{y,i}}=0,\frac{\partial B(\mathbf{D}_u)}{\partial u_{xx,i}}=0,\frac{\partial B(\mathbf{D}_u)}{\partial u_{yy,i}}=0,\frac{\partial B(\mathbf{D}_u)}{\partial u_{xy,i}}=0 \tag{5}$$

Eq. (5) are sorted into linear equations as follows:

$$\mathbf{A}\mathbf{D}_u=\mathbf{b} \tag{6}$$

where $\mathbf{A}=\mathbf{L}^T\boldsymbol{\omega}\mathbf{L}$, $\mathbf{b}=\mathbf{L}^T\boldsymbol{\omega}\mathbf{U}$, $\mathbf{L}=\left(\mathbf{L}_1{}^T,\mathbf{L}_2{}^T,\cdots,\mathbf{L}_n{}^T\right)^T$, $\mathbf{L}_i=\left(\Delta x_i,\Delta y_i,\frac{\Delta x_i^2}{2},\frac{\Delta y_i^2}{2},\Delta x_i\Delta y_i\right)$,

$\boldsymbol{\omega}=diag\left(\omega_1^2,\omega_2^2,\cdots,\omega_n^2\right)$, $\mathbf{U}=\left(u_1-u_i,u_2-u_i,\cdots,u_n-u_i\right)^T$.

Then, Eq. (6) is rewritten as:

$$\mathbf{D}_u=\left(u_{x,i},u_{y,i},u_{xx,i},u_{yy,i},u_{xy,i}\right)^T=\mathbf{A}^{-1}\mathbf{b}=\mathbf{A}^{-1}\mathbf{L}^T\boldsymbol{\omega}\mathbf{U}=\mathbf{M}\mathbf{U} \tag{7}$$

where $\mathbf{M}=\mathbf{A}^{-1}\mathbf{L}^T\boldsymbol{\omega}$.

For the convenience of representation, the elements of the matrix $\mathbf{M}$ are denoted as $m_{kj}^i$, where the superscript *i* indicates that the matrix M is calculated at node *i*. Then the generalized finite difference expressions of spatial derivatives of the unknown function at node *i* are:

$$\begin{aligned} &\frac{\partial u}{\partial x}\Big|_i=\sum_{j=1}^{n_i}m_{1j}^i\left(u_j-u_i\right),\ \frac{\partial u}{\partial y}\Big|_i=\sum_{j=1}^{n}m_{2j}^i\left(u_j-u_i\right),\ \frac{\partial^2 u}{\partial x^2}\Big|_i=\sum_{j=1}^{n}m_{3j}^i\left(u_j-u_i\right),\\ &\frac{\partial^2 u}{\partial y^2}\Big|_i=\sum_{j=1}^{n}m_{4j}^i\left(u_j-u_i\right),\ \frac{\partial^2 u}{\partial x\partial y}\Big|_i=\sum_{j=1}^{n}m_{5j}^i\left(u_j-u_i\right)\end{aligned} \tag{8}$$

## 2.2 Virtual construction of node control domains

Although GFDM has been widely used in various scientific and engineering problems, when there are singular Dirac function terms such as the singular point source in the equation, because no integration treatment is used in GFDM, it is difficult to integrate the singular Dirac term like FEM and FVM to eliminate the singularity in the governing equation. In the traditional FDM for reservoir numerical simulation, each node has a corresponding grid, so the strength of the singular point source can be approximated as the average source strength of the unit grid volume. It can be predicted that this processing method has good convergence, because the smaller the grid volume is, the closer the average source strength of the unit grid volume in FDM is to the real strength of the singular point source. However, as a meshless method, GFDM only has the concept of node but does not have the concept of grid volume. It is difficult to approximate the singular point source to the nonsingular grid average source term as FDM. Therefore, it is difficult to apply GFDM to the reservoir flow problem with point sources. In reservoir numerical simulation, because there will be injection and production wells whose scale is much smaller than the reservoir scale in the reservoir calculation domain, singular point sources often exist. To solve this problem, this paper attempts to build a theory about the node control domain and its volume in GFDM to effectively handle the point source term and form a node control domain-based meshless method.

***Def.* 1**: Each node is defined to have its control domain in addition to its influence domain, and the control domains of different nodes do not intersect each other, which meets:

$$\Omega = \bigcup_i \Omega_i, \ \Omega_i \bigcap \Omega_j = \varnothing, \ \sum_{i=1}^{n_p} V_i = V_\Omega \tag{9}$$

where $\Omega_i$ is the control domain of node *i*, $V_i$ is the volume of the domain $\Omega_i$, named the control volume of node *i*, $V_\Omega$ is the volume of the calculation domain $\Omega$, and $n_p$ is the total number of nodes in the calculation domain. Of course, due to the possible complex node distribution of the point cloud discretizing the computational domain, the geometric shape of the node control domain $\Omega_i$ can not be described explicitly, but this paper only cares about the volume $V_i$ of the corresponding node control domain, without figuring out the specific shape of the node control domain, this is why it is called virtual construction.

***Def.* 2**: The node cloud of the computational domain is called connectable means that for any two nodes in the point cloud (may as well suppose as node *i* and node *j*), if there is node *j* in the influence domain of node *i*, there must also be node *i* in the influence domain of node *j*.

***Def.* 3**: In the case of the connectable point cloud of the computational domain, if *j* is in the influence domain of node *i*, we denote that *i* is neighboring to *j*, and *j* is neighboring to *i*, and (*i*, *j*) is a neighboring pair.

Therefore, the connectable point cloud can be understood as adding the information of *i-j* pairs to the point cloud of the calculation domain, and the connectable point cloud is easy to construct. For example, to set the size of the influence domain of each node to be the same. Of course, the methods of constructing connectable point clouds can be diverse. In Section 2.4, a triangulation-based method of constructing a connectable point cloud will be developed.

Assuming the point cloud of the computational domain is connectable, for any node *i*, by integrating $\Delta u$ at node *i* in the node control domain $\Omega_i$, we can approximately obtain

$$\int_{\Omega_i} \Delta u \mathrm{d}\Omega \approx V_i \Delta u \approx V_i \sum_{j=1}^{n} \left( m_{3j}^i + m_{4j}^i \right) \left( u_j - u_i \right) \tag{10}$$

According to the divergence theorem, it is obtained that

$$\int_{\Omega_i} \Delta u \mathrm{d}\Omega \approx \int_{\Omega_i} \nabla \cdot \left( \nabla u \right) \mathrm{d}\Omega \approx \int_{\partial\Omega_i} \nabla u \cdot \vec{n} \mathrm{d}\Omega \approx \sum_{j=1}^{n} \int_{\partial\Omega_i \bigcap \partial\Omega_j} \nabla u \cdot \vec{n} \mathrm{d}\Omega \tag{11}$$

where $\partial\Omega_i$ and $\partial\Omega_j$ are the boundaries of $\Omega_i$ and $\Omega_j$ respectively, $\partial\Omega_i \bigcap \partial\Omega_j$ are the common boundary of $\Omega_i$ and $\Omega_j$, $\vec{n}$ is the external normal vectors of pairs, then we think that they are only related to $u_i$ and $u_j$. Therefore, by an analogy between Eq. (10) and Eq. (11), it is obtained that:

$$\int_{\partial\Omega_i \bigcap \partial\Omega_j} \nabla u \cdot \vec{n} \mathrm{d}\Omega = V_i \left( m_{3j}^i + m_{4j}^i \right) \left( u_j - u_i \right) \tag{12}$$

Because the point cloud is connectable, node *i* is also in the influence domain of node *j*. Similarly, when taking node *j* as the central node, it can be obtained that

$$\int_{\partial\Omega_i \bigcap \partial\Omega_j} \nabla u \cdot \vec{n}' \mathrm{d}\Omega = V_j \left( m_{3i}^j + m_{4i}^j \right) \left( u_i - u_j \right) \tag{13}$$

There is $\vec{n}' = -\vec{n}$, so we obtain that

$$V_i \left( m_{3j}^i + m_{4j}^i \right) \left( u_j - u_i \right) = V_j \left( m_{3i}^j + m_{4i}^j \right) \left( u_j - u_i \right) \Leftrightarrow V_i \left( m_{3j}^i + m_{4j}^i \right) = V_j \left( m_{3i}^j + m_{4i}^j \right) \tag{14}$$

It should be noted that if node *i* or node *j* is a boundary node, Eq. (14) is often not satisfied. This is because there are nodes only on one side of the boundary, which results in the center of the local point cloud in the boundary node's influence domain being far away from that boundary node [27, 28, 44], i.e. the local point cloud is of low quality, reducing the accuracy of the generalized difference expressions in Eq. (8) for approximating the spatial derivatives of the unknown function. Thus, for the boundary node, $m_{3j}^i + m_{4j}^i < 0$ generally exists, while for the inner node, it can generally meet $m_{3i}^j + m_{4i}^j \geq 0$, Eq. (14) is valid unless the control volume of the boundary node is negative, which is unreasonable. For example, as shown in Fig. 2 (a), node 3 is a boundary node, and the radius of the influence domain is selected as $\sqrt{2} + 0.1$. At this time, if only nodes 2, 3, 4, 7, 8, and 9 are included in the influence domain of node 3, the linear equations in Eq. (6) are unsolvable such that the generalized difference expressions in Eq. (8) cannot be obtained. This singular situation brings trouble to the construction of a generic calculation method of node control volumes. As shown in Fig. 2 (b), when the influence domain of node 3 is further expanded to 2.1 to include nodes 1, 5, and 13, Eq. (15) can be calculated by solving corresponding Eq. (6):

$$\Delta u = 6.1019\times10^{-6}\left(u_1-u_3\right)+1.0714\left(u_2-u_3\right)+1.0714\left(u_4-u_3\right)+6.1019\times10^{-6}\left(u_5-u_3\right)+ \\ -0.0714\left(u_7-u_3\right)-1.8573\left(u_8-u_3\right)-0.0714\left(u_9-u_3\right)+1.0000\left(u_{13}-u_3\right) \tag{15}$$

Thus, when the boundary node 3 is the central node, $m_{38}^3+m_{48}^3=-1.8573<0$. This example proves the above assertion, which makes it difficult for us to calculate the node control volume by using Eq. (14).

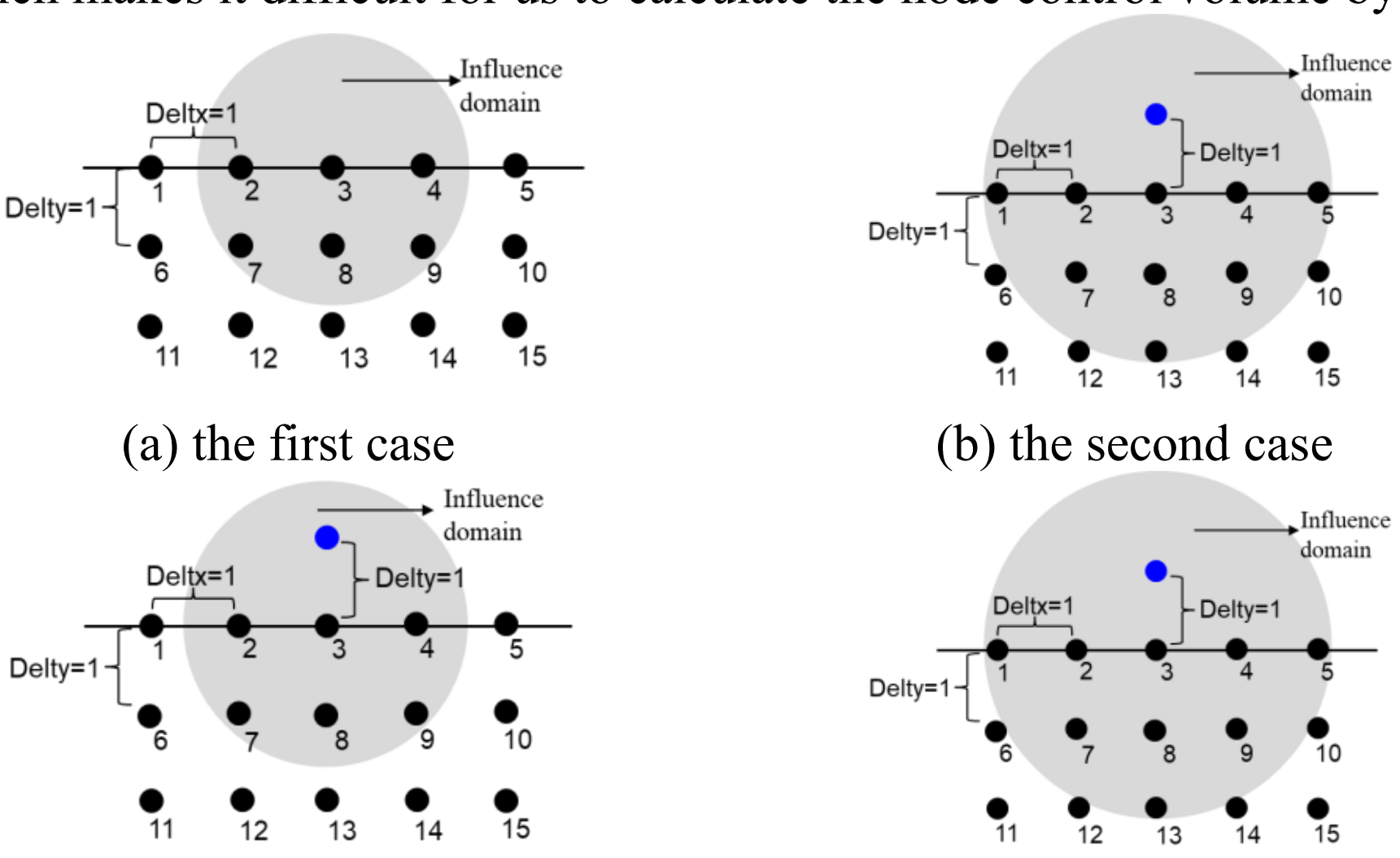

(a) the first case (b) the second case

(c) the first case with a virtual node (d) the second case with a virtual node

Fig. 2 The illustration of adding virtual nodes

We know that when dealing with derivative boundary conditions, GFDM can add a virtual node outside the calculation domain of derivative boundary nodes to improve the approximation accuracy of the generalized finite difference expressions of the spatial derivatives at boundary nodes in Eq. (8) [27, 28]. It inspires us that when calculating the node control volume, to eliminate the problems caused by boundary nodes, virtual nodes corresponding to boundary nodes can be added outside the boundary of the calculation domain, so that for any boundary node *i*, $\left(m_{3j}^i+m_{4j}^i\right)\geq 0$ holds, and the problem that the linear equations in Eq. (6) are unsolvable and Eq. (8) cannot be obtained can be overcome. Still taking Fig. 2 as an example, we only add a virtual node *M* (marked in blue) in the influence domain of node 3 to obtain Fig. 2 (c) (in fact, two virtual nodes also need to be added outside boundary nodes 2 and 4. The illustration here is used to show that even if only one virtual node is added outside the concerned central node 3, it is sufficient to meet the requirement). In this case, using Eq. (6) it can be calculated as:

$$\Delta u = 1.0000\left(u_M-u_3\right)+0.9998\left(u_2-u_3\right)+0.9998\left(u_4-u_3\right)+ \\ 1.7590\times10^{-4}\left(u_7-u_3\right)+0.9996\left(u_8-u_3\right)+1.7590\times10^{-4}\left(u_9-u_3\right) \tag{16}$$

It can be seen that linear equations in Eq. (6) after adding virtual nodes become solvable, and $m_{38}^3+m_{48}^3=0.9996>0$ is obtained.

Similarly, Fig. 2 (b) is added the virtual node *M* to become Fig. 2 (d), then it can be calculated that

$$\Delta u = 1.0000\left(u_M-u_3\right)+4.8824\times10^{-6}\left(u_1-u_3\right)+0.8572\left(u_2-u_3\right)+1.0714\left(u_4-u_3\right) \\ +4.8824\times10^{-6}\left(u_5-u_3\right)++0.1427\left(u_7-u_3\right)+0.7145\left(u_8-u_3\right)+0.1427\left(u_9-u_3\right) \\ +4.4760\times10^{-6}\left(u_{13}-u_3\right) \tag{17}$$

As seen from Eq. (17), after adding virtual nodes outside the boundary node, it is obtained that $m_{38}^3+m_{48}^3=0.7145>0$.

In addition, it should also be noted that, after adding a virtual node outside the boundary node, the control domain in Eq. (14) of the boundary node will not be completely included in the calculation domain. At this time, Eq. (14) needs to be modified as follows:

$$V_i\left(m_{3j}^i+m_{4j}^i\right)=V_j\left(m_{3i}^j+m_{4i}^j\right) \tag{18}$$

where $V_i$ is the control volume of each node after adding virtual nodes on the boundary, rather than the node control volume $\overline{V_i}$ only inside the calculation domain.

***Def.* 4**: Then, we introduce the concept of feature angle at each node to establish the relationship between $V_i$ and $\bar{V}_i$. The node characteristic angle introduced here is the same as that defined at the node in the boundary element method (BEM). For the smooth node on the boundary shown in Fig. 2 (a), the characteristic angle is generally $\pi$, for the corner node on the boundary shown in Fig. 2 (b), the characteristic angle is the included angle $\theta$ of tangents on both sides, and for the internal nodes in the calculation domain, the characteristic angle of the node is $2\pi$. Therefore, denote $\theta_i$ as the characteristic angle of node *i*, it is obtained that:

$$\bar{V}_i = \frac{\theta_i}{2\pi} V_i \tag{19}$$

Thus, Eq. (9) is rewritten as:

$$\sum_{i=1}^{n_p} \frac{\theta_i}{2\pi} V_i = V_\Omega \tag{20}$$

Finally, if there are $n_c$ pairs of *i-j*, there is an Eq. (18) for each pair of *i-j*. combined with Eq. (20), $n_c$ + 1 linear equations can be obtained, and there are $n_p$ unknown node control volumes, and $n_c > n_p$ generally holds. Therefore, Eq. (18) and Eq. (20) form an overdetermined linear equation system, which can be solved via the least square method, and then calculated $\bar{V}_i$ by using Eq. (19). Note that the above derivation process has a premise assumption, that is, the point cloud of the computational domain is connectable. In addition, another detail needs to explain. Since Eq. (18) and Eq. (20) form an overdetermined system of equations, which is solved by minimizing the residue in the 2-norm of the equations, Eq. (20) may not be strictly satisfied. However, the sum of the control volume of each node is equal to the total volume of the calculation domain, which is the key to the mass conservation of the calculation domain. Therefore, Eq. (20) is required to be strictly satisfied. Therefore, this paper uses the "large number method" to multiply both sides of Eq. (20) by a larger number *G*. That is, Eq. (20) is modified as:

$$G\sum_{i=1}^{n_p} \frac{\theta_i}{2\pi} V_i = G V_\Omega \tag{21}$$

In Section 2.3, an empirical method that can calculate reliable node control volumes will be given, and detailed verification will be illustrated in Section 3 “Results of numerical analysis”. Here, it is illustrated that the reasonability of the above calculation method of node control volumes with a ‘naïve’ case where simple Cartesian node collocation is used. As shown in Fig. 3, it is a rectangular calculation domain [0, 20] × [0, 20], using 3×3 uniformly distributed nodes and virtual nodes are added outside the boundary of the calculation domain. May as well select the radius of the influence domain as $10\sqrt{2}+0.1$, and then it can be obtained that the node distribution in the influence domain of each internal or boundary node is the same and symmetrical. Therefore, for any *i-j* pair, $\left(m_{3j}^i + m_{4j}^i\right) = \left(m_{3i}^j + m_{4i}^j\right)$ holds. Then according to Eq. (17), $\forall i, j, V_i = V_j$ holds. For the nodes at the four corners of the calculation domain, there is $\theta_i = \frac{\pi}{2}$ ($i = 1,3,7,9$), For the nodes at the midpoint of the four edges of the calculation domain, there is $\theta_i = \pi$ ($i = 2,4,6,8$), and for the internal nodes, there is $\theta_i = 2\pi$ ($i = 5$), to obtain:

$$\begin{aligned}&\sum_{i=1}^{n_p} \frac{\theta_i}{2\pi} V_i = 4\times\frac{1}{4}V_5 + 4\times\frac{1}{2}V_5 + V_5 = 20\times 20 \Rightarrow V_i = 100 \Rightarrow \bar{V}_5 = 100\\ &\Rightarrow \bar{V}_1 = \bar{V}_3 = \bar{V}_7 = \bar{V}_9 = 25, \bar{V}_2 = \bar{V}_4 = \bar{V}_6 = \bar{V}_8 = 50\end{aligned} \tag{22}$$

The above-calculated node control volumes are matched with the grid volumes in the Cartesian grid, which shows the reasonability of the above definitions and equations about the node control volume.

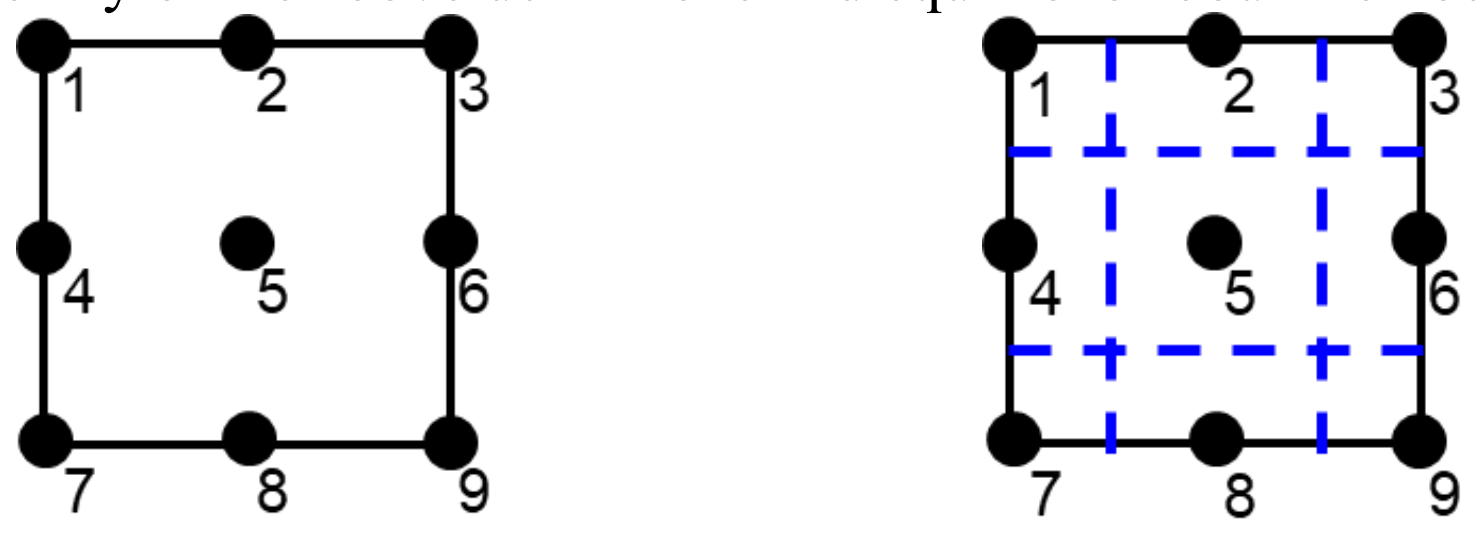

(a) uniform node collocation (b) corresponding orthogonal mesh

Fig. 3 A sketch of the naive case

From the above process, we can predict that the introduced node control volume in GFDM is related not only to the node coordinates of the point cloud but also to the size of the influence domain (i.e. *i-j* pairs between nodes), which is expressed in the form of an abstract function as follows:

$$V_i = f\left(X,C\right) \tag{23}$$

where $X$ is the information of nodes' coordinates and $C$ is information of the *i-j* pairs. This is also easy to understand because the *X* and *C* in GFDM can correspond to the mesh center and mesh neighboring relationship in the mesh division of the calculation domain, while the mesh center and mesh neighboring relationship can uniquely determine the mesh volume. Correspondingly, in the meshless system, the node control volume is uniquely determined by *X* and *C*. In numerical test cases of Section 3, we will show that different influence domain radii will yield different node control volumes.

It should be noted that if building a connectable point cloud by setting each node the same radius of the node influence domain, the given radius value cannot be too large. Because if the given radius value is large, there will be a large number of nodes in the node influence domain, which will bring three problems:

(i) The dimension of the overdetermined equations is large, reducing the calculation efficiency.

(ii) Rao et al. [27, 28] pointed out that, for convection-dominated heat and mass transfer and multiphase flow problems, the large radius of the influence domain will reduce the quality of node distribution in the influence domain of the nodes at or near the boundary, causing the decrease of the approximation accuracy of the generalized difference operator. In addition, the large radius of the influence domain will also increase the dissipation error of the single point upwind (SPU) scheme used for discretization of the convection term.

(iii) The *i-j* pairs are complex, resulting in a large bandwidth of the Jacobian matrix in the Newton iteration process, increasing the computational cost of the linear solver.

2.3 An empirical method of calculating node control volumes

For the narration convenience, the quartic spline weight function in Eq. (3) is denoted as $w_1$, and the negative cubic function in Eq. (4) is denoted as $w_2$.

First, based on implementation experience, we give the method of obtaining reliable node control volumes:

***Step*** **1**: if you use the method of setting the same node influence domain radius to build a connectable point cloud, first give the influence domain radius $r_m$, and build the corresponding connectable point cloud (i.e. *i-j* pairs).

***Step*** **2**: use the function in Eq. (4) as the weight function, i.e. $w_2$, to obtain the overdetermined equations composed of Eq. (18) and Eq. (21). In Section 3, we will explain in detail through specific examples that $w_2$ can obtain a higher-accuracy node control volume distribution than $w_1$.

***Step*** **3**: considering that the solution of overdetermined equations generally adopts the least square method, but the result in the sense of least square is not necessarily a reasonable node control volume. On the one hand, consider a common case, *i* is a node in the area with high node density, and *j* is a node in the area with normal node density. In the connectable point cloud, if *i* is neighboring with to *j*, because the node density in the area where node *i* is located is large, the local point cloud of node *i* contains more nodes, and because node *j* is close to the edge of the influence domain of node *i*, that is, the distance between node *j* and node *i* $d_{ij}$ is close to the radius of the influence domain of node *i*, Then the weight of node *j* in the local point cloud of node *i* will be very small, but the local point cloud of node *j* does not contain a large number of nodes, and the distance between node *j* and its adjacent nodes is similar to $d_{ij}$. Node *i* has a large weight in the local point cloud of node *j*, so it will lead to a great difference between $\left(m_{3j}^{i}+m_{4j}^{i}\right)$ and $\left(m_{3i}^{j}+m_{4i}^{j}\right)$. According to implementation experience, in this case, the Eq. (18) corresponding to the *i-j* pair will affect the accuracy of calculating node control volumes by solving the overdetermined equations in ***Step*** 2. On the other hand, taking the Cartesian collocation in Section 3.1 as an example, when the radius of the influence domain becomes larger, while the node control volume profile is controlled by the complexity of the *i-j* pairs, it can be expected that in Cartesian collocation, for two neighboring nodes far away from the boundary of the calculation domain (may as well denoted as node $i$ and node $j$), the local node clouds of these two nodes is the same. Therefore, there should be:

$$\left(m_{3j}^{i}+m_{4j}^{i}\right)=\left(m_{3i}^{j}+m_{4i}^{j}\right) \tag{24}$$

So,

$$V_i=V_j \tag{25}$$

In the sense of least squares, Eq. (25) may not be strictly satisfied, but it can be expected that the accurate node control volume distribution should meet $V_i \approx V_j$.

Based on the analysis of the above two aspects, it is expected that the satisfying priority of the Eq. (18) obtained when $\left(m_{3j}^{i}+m_{4j}^{i}\right)$ and $\left(m_{3i}^{j}+m_{4i}^{j}\right)$ are close in solving the overdetermined equations is higher than that of the Eq. (18) obtained when $\left(m_{3j}^{i}+m_{4j}^{i}\right)$ and $\left(m_{3i}^{j}+m_{4i}^{j}\right)$ are significantly different. This can be achieved by empirically weighting the equation.

$$\omega_{ij}\left[V_i\left(m_{3j}^{i}+m_{4j}^{i}\right)-V_j\left(m_{3i}^{j}+m_{4i}^{j}\right)\right]=0\,,\quad \omega_{ij}=\frac{\min\left(m_{3j}^{i}+m_{4j}^{i},m_{3i}^{j}+m_{4i}^{j}\right)}{\max\left(m_{3j}^{i}+m_{4j}^{i},m_{3i}^{j}+m_{4i}^{j}\right)} \tag{26}$$

***Step* 4**: Node control volumes are calculated by solving the overdetermined equations composed of Eq. (26) and (21) with the least square method.

In Section 3 "results of numerical analysis", we will explain in detail the proposed method outperforms other methods of calculating node control volumes.

2.4 A triangulation-based method of determining the connectable point cloud

In Section 2.2, we give a method to determine the connectable point cloud by setting the same-radius influence domain for all nodes. However, in the case of point clouds with uneven node density, the same-radius node influence domain will increase the complexity of local point clouds of the nodes in the high-node-density domain. Implementation experience tells us that this case will reduce the calculation accuracy of node control volumes. Therefore, how to reasonably determine the connectable point cloud (the connectable point cloud is the point cloud with the information of $i$-$j$ pairs) is an important issue. In this section, a method to determine the connectable point cloud that can obtain high-accuracy node control volumes will be given for the point cloud composed of the mesh vertices of the Delaunay triangulation of the computational domain:

***Step* 1**: for the triangular mesh, first determine the preliminary local point cloud of nodes according to the principle of "if there is an edge between two nodes, then they are neighboring";

***Step* 2**: for the node inside the computational domain, may as well be denoted as node $i$. if the number of nodes adjacent to node $i$ denoted by $n_i$ is less than 5, which means that the generalized difference operator of the local point cloud of node $i$ cannot be closed solved by Eq. (6), add 5-$n_i$ nodes (denoted as $j$, $j=1,\cdots,n_i$) closest to node $i$ from the nodes in the computational domain except the nodes contained in the local point cloud of node $i$ to the local point cloud of node $i$. At the same time, add node $i$ to the local point cloud of node $j$ ( $j=1,\cdots,n_i$ ). When using the weight function to calculate the weight in Eq. (2), the radius of the influence domain $r_{mi}$ of node $i$ selects 1.5 times the maximum value of the distance between node $i$ and the nodes in its local point cloud. In fact, when using the weight function in Eq. (4) (i.e., $w_2$), no matter what the radius value is, the calculated relative weight between the nodes in the local point cloud of node $i$ is the same, resulting in the same generalized finite difference expressions. Therefore, when the local point cloud is known, the radius of the node influence domain in $w_2$ does not affect the calculation results.

***Step* 3**: for the nodes on the boundary of the computational domain, which may also as well denoted as node $i$, $n_i$ is generally less than 5. There are virtual nodes that need to be added outside the computational domain, therefore, for the boundary node $i$, 5-$n_i$ nodes from the virtual nodes closest to node $i$ need to be added to the local point cloud of node $i$. The radius of the node influence domain used for weight calculation is determined in the same way as in ***Step* 2**.

In Section 3 "results of numerical analysis", we will explain in detail that the proposed connectable point cloud determination method in this section can obtain a high-accuracy node control volume profile than the method giving the same-radius node influence domain, also obtain higher-accuracy simulation results.

2.5 NCDMM discrete scheme of governing equations

This paper takes the two-phase porous flow equations in Eq. (27) as an example to illustrate the specific NCDMM-based discrete scheme.

$$\nabla \cdot \left( \frac{kk_{ra}}{B_a \mu_a} \nabla p_a \right) + q_{a,well} = \frac{\partial}{\partial t} \left( \frac{\phi S_a}{B_a} \right) \tag{27}$$

where $k$ is the permeability, $a = o$ or $w$, which denotes oil phase or water phase; $k_{ra} = k_{ra}(S_w)$ are the relative permeability of phase $a$; $B_a$ is the volume coefficient of the phase $a$; $\mu_a$ is the viscosity of phase $a$; $p_a$ is the pressure of phase $a$; $t$ is time; $S_a$ is the saturation of phase $a$, which meets $S_o + S_w = 1$; $\phi = \phi(p)$ is the porosity which is a function of the pressure. $q_{a,well}$ is the strength of the source or sink of phase $a$ from production or injection wells.

Assuming that the node pressure and water saturation at time $t$ are known, integrate both sides of Eq. (27) on the node control volume of node $i$, and take the implicit scheme for each physical quantity, then the integration of the first term on the left side of Eq. (27) is approximated as:

$$\int_{V_i} \nabla \cdot \left( \frac{kk_{ra}}{B_a \mu_a} \nabla p_a \right) d\Omega = \int_{V_i} \sum_{j=1}^{n_i} \left[ \frac{k_{ij}^{t+\Delta t} k_{ra,ij}^{t+\Delta t}}{B_{a,ij}^{t+\Delta t} \mu_{a,ij}^{t+\Delta t}} \left( m_{3,j}^i + m_{4,j}^i \right) \left( p_{a,(i,j)}^{t+\Delta t} - p_{a,i}^{t+\Delta t} \right) \right] d\Omega \tag{28}$$

in which, it should be emphasized that $n_i$ is the number of real nodes in the influence domain of node $i$, $(i, j)$ refers to the sequence number in all nodes of the $j$-th real node in the influence domain of node $i$ except node $i$ itself. The single-point upstream scheme in Eq. (29) is used for relative permeability to handle the convection term stably [27, 28], the harmonic average scheme in Eq. (30) is used for permeability, and the arithmetic average scheme in Eq. (31) is used for phase volume coefficient and viscosity.

$$k_{ra,ij}^{t+\Delta t} = \begin{cases} k_{ra,i}\left(S_{w,i}^{t+\Delta t}\right) & if\ p_{a,(i,j)}^{t+\Delta t} \geq p_{a,i}^{t+\Delta t} \\ k_{ra,i}\left(S_{w,i}^{t+\Delta t}\right) & if\ p_{a,(i,j)}^{t+\Delta t} < p_{a,i}^{t+\Delta t} \end{cases}, \ a = o \text{ or } w \tag{29}$$

$$k_{ij}^{t+\Delta t} = \frac{2}{1/k_i^{t+\Delta t} + 1/k_j^{t+\Delta t}} \tag{30}$$

$$B_{a,ij}^{t+\Delta t} = \frac{B_{a,i}^{t+\Delta t} + B_{a,j}^{t+\Delta t}}{2}, \quad \mu_{a,ij}^{t+\Delta t} = \frac{\mu_{a,i}^{t+\Delta t} + \mu_{a,j}^{t+\Delta t}}{2}, \ a = o \text{ or } w \tag{31}$$

Eq. (28) can be further approximated as:

$$\int_{V_i} \nabla \cdot \left( \frac{kk_{ra}}{B_a \mu_a} \nabla p_a \right) d\Omega = \sum_{j=1}^{n_i} \left[ \frac{k_{ij}^{t+\Delta t} k_{ra,ij}^{t+\Delta t}}{B_{a,ij}^{t+\Delta t} \mu_{a,ij}^{t+\Delta t}} V \left( m_{3,j}^i + m_{4,j}^i \right) \left( p_{a,(i,j)}^{t+\Delta t} - p_{a,i}^{t+\Delta t} \right) \right], \ a = o \text{ or } w \tag{32}$$

For the second term of the left side of Eq. (27), if there is no singular source or sink term (i.e. no production well or injection well) at node $i$, then:

$$\int_{V_i} q_{a,well} d\Omega = Q_{a,well} = 0, \ a = o \text{ or } w \tag{33}$$

If there are production wells or injection wells at node $i$, since the proposed NCDMM is essentially a meshless method, although the node control volume value has been obtained in Section 2.3, it is difficult to figure out the specific geometric shape of the node control domain. Therefore, in NCDMM, the well index is calculated as:

$$WI_i = \frac{2\pi k_i h_i}{\ln\left(r_{e,i}/r_w\right) + s}, \quad r_{e,i} = 0.14\sqrt{2V_i/h_i} \tag{34}$$

where $r_{e,i}$ is the equivalent radius, $r_w$ is the well radius, and $s$ is the skin factor.

If the well is a production well, then

$$\int_{V_i} q_{a,well} d\Omega = Q_{a,well} = \frac{k_{ra,if}^{t+\Delta t}}{B_{a,if}^{t+\Delta t} \mu_{a,if}^{t+\Delta t}} WI_i \left( p_{wf}^{t+\Delta t} - p_{a,i}^{t+\Delta t} \right), \ a = o \text{ or } w \tag{35}$$

If the well is a water injection well, then

$$\int_{V_i} q_{o,well} d\Omega = Q_{o,well} = 0 \tag{36}$$

$$\int_{V_i} q_{w,well} d\Omega = Q_{w,well} = \frac{1}{B_{w,if}{}^{t+\Delta t}} \left( \frac{k_{ro,if}{}^{t+\Delta t}}{\mu_{o,if}{}^{t+\Delta t}} + \frac{k_{rw,if}{}^{t+\Delta t}}{\mu_{w,if}{}^{t+\Delta t}} \right) WI_i \left( p_{wf}{}^{t+\Delta t} - p_{w,i}{}^{t+\Delta t} \right) \tag{37}$$

The integral of the term on the right side of Eq. (27) is approximated as:

$$\int_{V_i} \frac{\partial(\phi S_a)}{\partial t} d\Omega = \frac{V_i}{\Delta t} \left[ \left( \frac{\phi_i^{t+\Delta t} S_{a,i}{}^{t+\Delta t}}{B_{a,i}{}^{t+\Delta t}} \right) - \left( \frac{\phi_i^t S_{a,i}{}^t}{B_{a,i}{}^t} \right) \right], \; a = o \text{ or } w \tag{38}$$

Based on the above results, the approximate expressions of Eq. (27) can be obtained as follows:

$$\sum_{j=1}^{n_i} \left[ \frac{k_{ij}{}^{t+\Delta t} k_{ra,ij}{}^{t+\Delta t}}{B_{a,ij}{}^{t+\Delta t} \mu_{a,ij}{}^{t+\Delta t}} V_i \left( m_{3,j}^i + m_{4,j}^i \right) \left( p_{a,(i,j)}{}^{t+\Delta t} - p_{a,i}{}^{t+\Delta t} \right) \right] + Q_{a,well} = \frac{V_i}{\Delta t} \left[ \left( \frac{\phi_i^{t+\Delta t} S_{a,i}{}^{t+\Delta t}}{B_{a,i}{}^{t+\Delta t}} \right) - \left( \frac{\phi_i^t S_{a,i}{}^t}{B_{a,i}{}^t} \right) \right] \tag{39}$$

To make direct use of the existing nonlinear solver in current FVM-based reservoir numerical simulators, the concept of transmissibility similar to that in the FVM-based discrete scheme can be extracted from the first item on the left side of Eq. (39), that is, when node $i$ is considered (i.e. consider the influence domain of node $i$), transmissibility between node $i$ and node $j$ (if more rigorous, $j$ should be replaced by $(i, j)$, because $(i, j)$ is the sequential number in all nodes of the $j$-th node in the influence domain of node $i$. Here denoted as $j$ for simplicity of expression) is

$$T_{i,ij}{}^{t+\Delta t} = \left[ k_{ij}{}^{t+\Delta t} V_i \left( m_{3,j}^i + m_{4,j}^i \right) \right] \tag{40}$$

where $i$ in the subscript $i,ij$ means the transmissibility between node $i$ and node $j$ is obtained when the influence domain of node $i$ is considered.

Similarly, when considering node $j$, we can obtain

$$T_{j,ij}{}^{t+\Delta t} = \left[ k_{ji}{}^{t+\Delta t} V_j \left( m_{3,i}^j + m_{4,i}^j \right) \right] \tag{41}$$

Combined with Eq. (18) and Eq. (30), it is obtained that

$$T_{i,ij}{}^{t+\Delta t} = T_{j,ij}{}^{t+\Delta t} \tag{42}$$

Eq. (42) shows that the proposed NCDMM satisfies the local mass conservation [45, 46, 47, 48], which is important in flow simulation, then the transmissibility can be uniformly written as:

$$T_{ij}{}^{t+\Delta t} = T_{i,ij}{}^{t+\Delta t} = T_{j,ij}{}^{t+\Delta t} = k_{ij}{}^{t+\Delta t} V_i \left( m_{3,j}^i + m_{4,j}^i \right) = k_{ji}{}^{t+\Delta t} V_j \left( m_{3,i}^j + m_{4,i}^j \right) \tag{43}$$

Then Eq. (39) can be rewritten as

$$\sum_{j=1}^{n_i} \left[ \frac{k_{ra,ij}{}^{t+\Delta t}}{B_{a,ij}{}^{t+\Delta t} \mu_{a,ij}{}^{t+\Delta t}} T_{ij}{}^{t+\Delta t} \left( p_{a,(i,j)}{}^{t+\Delta t} - p_{a,i}{}^{t+\Delta t} \right) \right] + Q_{a,well} = \frac{V_i}{\Delta t} \left[ \left( \frac{\phi_i^{t+\Delta t} S_{a,i}{}^{t+\Delta t}}{B_{a,i}{}^{t+\Delta t}} \right) - \left( \frac{\phi_i^t S_{a,i}{}^t}{B_{a,i}{}^t} \right) \right] \tag{44}$$

It can be seen that the proposed NCDMM firstly introduces the basic concepts of the node control volume and the connectable point cloud to GFDM, then obtains the NCDMM discrete scheme satisfying the local mass conservation by integrating the GFDM discrete scheme on the node control volume.

## 2.6 Treatment of boundary conditions and construction of global equations

As can be seen from Section 2.5, the proposed NCDMM is essentially the integral form of GFDM, but this integral form refers to the integral of the discretization scheme of the governing equation. For the boundary conditions, the treatment of NCDMM is the same as that of GFDM, so compared with the traditional FVM, NCDMM can directly handle various boundary conditions. When NCDMM deals with derivative boundary conditions, the virtual nodes added at the boundary node need to be used to improve the approximation accuracy of generalized finite difference operators for the spatial derivatives of unknown functions at the boundary node [27, 28, 44]. Fig. 4 is a sketch of the virtual node added at the derivative boundary node from the literature [27]. It can be seen that, when the blue virtual node is not added, the black nodes in the influence domain of the red boundary node are on the same side of the tangent at the boundary node. At this time, the position of the center of gravity of all nodes in the boundary-node influence domain deviates much from that of the considered boundary node, such that the accuracy of the spatial derivatives of the unknown function given by Eq. (8) is low. After adding a blue virtual node in the influence domain of the boundary node, the accuracy of Eq. (8) for spatial derivatives will be significantly improved. Therefore, the virtual node corresponding to the derivative boundary node will also participate in the construction of the global equations. For details, please refer to the literature [27].

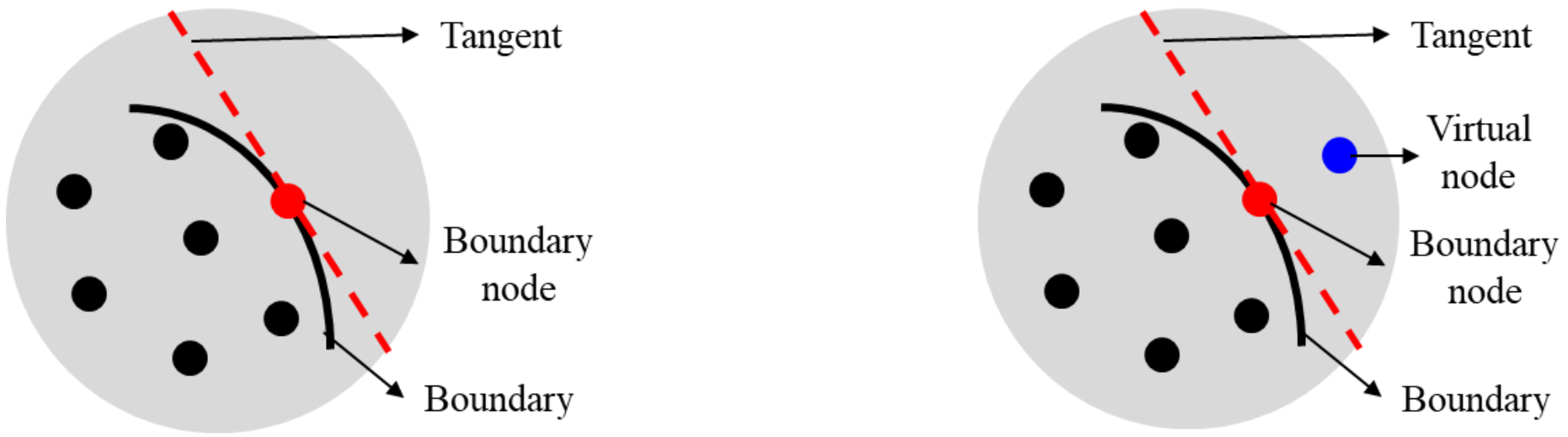

(a) the original case (b) the case after adding the virtual node

Fig. 4 Sketch of using the virtual node to handle the derivative boundary condition [27]

If a boundary node *i* satisfies the first-type (Dirichlet) boundary conditions in Eq. (45), the discrete equation at node *i* is shown in Eq. (46):

$$p|_{\Gamma} = f\left(x,y\right),\quad S_w|_{\Gamma} = g\left(x,y\right) \tag{45}$$

$$p_{o,i} = f\left(x_i,y_i\right),\quad S_{w,i} = g\left(x_i,y_i\right) \tag{46}$$

where $x_i$ and $y_i$ are the *x* coordinate and *y* coordinate of node *i* respectively.

Suppose a boundary node *i* satisfies the second-type (Neumann) boundary conditions, and the virtual node added for the treatment of the derivative boundary condition at node *i* is node *j*, then in the global equations, the equations at node *i* is the discrete schemes of the governing equations in Eq. (44). Corresponding to the Neumann boundary conditions in Eq. (47), the equation at node *j* is Eq. (48).

$$\frac{\partial p}{\partial \mathbf{n}}|_{\Gamma} = f\left(x,y\right) \tag{47}$$

$$n_x\sum_{l=1}^{n_i} m_{1l}^i\left(p_{o,(i,l)} - p_{o,i}\right) + n_y\sum_{l=1}^{n_i} m_{2l}^i\left(p_{o,(i,l)} - p_{o,i}\right) = f\left(x_i,y_i\right) \tag{48}$$

where, $n_x$ and $n_y$ are the components of the external normal vector $\mathbf{n}$ at the boundary node, and $n_i$ is the number of nodes in the influence domain of node *i* other than node *i* itself (including virtual nodes). Since virtual node *j* is in the influence domain of node *i*, there is *l* such that $(i,l) = j$.

For the third-type (Robin-type) boundary condition which is also the derivative boundary condition, the treatment method is the same as that for Neumann boundary conditions.

Therefore, suppose there are $n_m$ nodes in the calculation domain, including $n_{1BC}$ nodes that meet the first-type (Dirichlet) boundary condition and $n_{2BC}$ nodes that meet the derivative boundary condition, then the number of the virtual nodes involved in the construction of the global equations are also $n_{2BC}$, so the nodes participating in the construction of the global equations are $n_m$ real nodes and the $n_{2BC}$ virtual nodes corresponding to the nodes which meet the derivative boundary conditions. Among these nodes, for the $n_m$ real nodes, the equation at the interior node or the boundary node that meet the derivative boundary conditions is Eq. (44), and the equation corresponding to the node which meets the Dirichlet boundary condition is Eq. (46). For $n_{2BC}$ virtual nodes, the corresponding equations are Eq. (48).

Therefore, there are only two types of equations in global equations, one is Eq. (44) representing the NCDMM-based discrete governing equations, and the other is Eq. (46) and Eq. (48) representing the boundary conditions. For the first type, the node control volumes and transmissibilities of *i-j* pairs in NCDMM are the required parameters in the nonlinear solver for the FVM-based reservoir numerical simulator. Therefore, the NCDMM discrete governing equations shown in Eq. (44) at the $n_m - n_{1BC} - n_{2BC}$ interior nodes and $n_{2BC}$ derivative-boundary-condition nodes can be quickly constructed like the traditional FVM. Then, Eq. (48) representing the derivative boundary conditions are filled at the virtual nodes corresponding to the $n_{2BC}$ derivative-boundary-condition nodes, and finally, the equations shown in Eq. (46) are added to the Dirichlet-boundary-condition nodes. Of course, if the bottom hole pressure (BHP) needs to be calculated, the equation between well production rates and BHP needs to be added. The nonlinear solver in the FVM-based simulator can be employed to solve the formed global equations. It can be seen that, Eqs. (46) and (48) representing the boundary conditions are generally linear. Therefore, compared with the traditional FVM, the proposed NCDMM does not increase the nonlinearity of the equations, such that the computational efficiency of the proposed NCDMM is comparable to that of traditional FVM. The numerical examples in Section 3 will prove this assertion with specific data. Note that, when it is the most common closed boundary condition in reservoir simulation, the discrete equation of the boundary condition is not required, but only the connection table composed of the connection between the nodes (i.e., *i-j* pairs) and the corresponding transmissibility is constructed without considering the connection between the boundary node and virtual nodes in the first item on the left of Eq. (44) at each boundary node. In this case, virtual boundary nodes are only used as geometric entities to help calculate the node control volumes and the transmissibility of connections between the boundary node and its neighboring real nodes, and there is not any physical quantity at the virtual nodes that need to be solved.

2.7 Convergence analysis of NCDMM

We know that the traditional FDM has low-order accuracy, that is, the convergence order is less than or equal to 2, as shown in Fig. 4 (a), node 0 is the central node, with $\Delta x_1 = -\Delta x$ and $\Delta y_1 = 0$ for node 1, $\Delta x_2 = \Delta x$ and $\Delta y_2 = 0$ for node 2, $\Delta x_3 = 0$ and $\Delta y_3 = -\Delta y$ for node 3, $\Delta x_4 = 0$ and $\Delta y_4 = \Delta y$ for node 4. Then in the traditional FDM, the second-order difference expressions of spatial derivatives of the unknown function at node 0 are

$$u_{xx,0} = \frac{u_1 + u_2 - 2u_0}{\Delta x^2} + O\left(\Delta x^2\right),\ u_{yy,0} = \frac{u_3 + u_4 - 2u_0}{\Delta y^2} + O\left(\Delta y^2\right),\ u_{x,0} = \frac{u_2 - u_1}{2\Delta x} + O\left(\Delta x^2\right),$$
$$u_{y,0} = \frac{u_4 - u_3}{2\Delta y} + O\left(\Delta y^2\right) \tag{49}$$

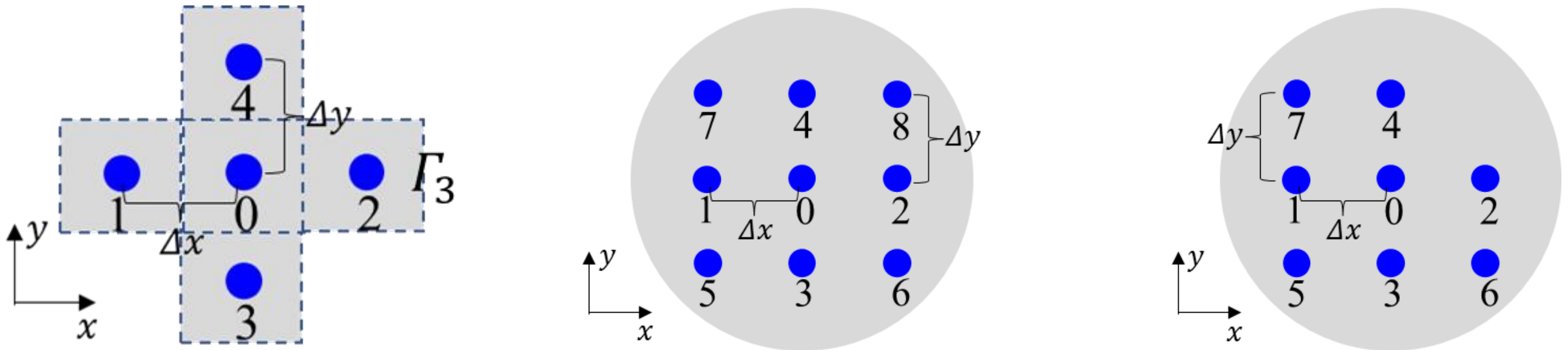

(a) Cartesian grid (b) a symmetric local point cloud (c) an asymmetric local point cloud

Fig. 5 Sketches of the Cartesian grid, a symmetric local point cloud, and an asymmetric local point cloud

Fu et al. [49] have proved the consistency of GFDM schemes. The key lies in the convergence order of GFDM. Firstly, according to the generalized finite difference expressions in Section 2.1 derived by ignoring third-order and higher-order derivative terms of Taylor expansion, it can be inferred that the convergence order of GFDM is also less than or equal to 2. Here, the generalized finite difference expressions are calculated based on two different local point clouds in Fig. 5 (b) and (c) respectively, to show that GFDM has the second-order accuracy when the local point cloud symmetry is good. In the case of local point clouds with poor symmetry, it generally has only the first-order accuracy, which is the cost that FDM needs to pay when it is extended to GFDM in the case of arbitrary point clouds.

(*i*) Example 1: In the local point cloud shown in Fig. 5 (b), the central node is node 0. May as well suppose $\Delta x = \Delta y$, and take the radius of the influence domain of node 0 as $1.8\Delta x$, so there are other 8 nodes (i.e. from node 1 to node 8) around node 0 in its influence domain.

Then, using the GFDM theory in Section 2.1 may as well choose Eq. (3) as the weight function, then the generalized difference expressions of spatial derivatives at node 0 can be calculated. Taking $u_{xx0}$ as an example, it is obtained that

$$u_{xx0} \approx \sum_{j=1}^{8} m_{3j}\left(u_j - u_0\right) \tag{50}$$

in which

$$\left(m_{31},m_{32},m_{33},m_{34},m_{35},m_{36},m_{37},m_{38}\right)=$$

$$\begin{pmatrix} \dfrac{9.6308\times10^{-1}}{\Delta x^2},\dfrac{9.6308\times10^{-1}}{\Delta x^2},-\dfrac{3.6917\times10^{-2}}{\Delta x^2},-\dfrac{3.6917\times10^{-2}}{\Delta x^2}, \\ \dfrac{1.8459\times10^{-2}}{\Delta x^2},\dfrac{2.4918\times10^{-2}}{\Delta x^2},\dfrac{1.8459\times10^{-2}}{\Delta x^2},\dfrac{1.8459\times10^{-2}}{\Delta x^2} \end{pmatrix} \tag{51}$$

According to Taylor's expansion of the unknown function at node 0, we can get:

$$\begin{aligned}
u_1&=u_0-\Delta x u_{x0}+\frac{1}{2}\Delta x^2 u_{xx0}-\frac{1}{6}\Delta x^3 u_{xxx0}+O\left(r^4\right),\\
u_2&=u_0+\Delta x u_{x0}+\frac{1}{2}\Delta x^2 u_{xx0}+\frac{1}{6}\Delta x^3 u_{xxx0}+O\left(r^4\right),\\
u_3&=u_0-\Delta y u_{y0}+\frac{1}{2}\Delta y^2 u_{yy0}-\frac{1}{6}\Delta y^3 u_{yyy0}+O\left(r^4\right),\\
u_4&=u_0+\Delta y u_{y0}+\frac{1}{2}\Delta y^2 u_{yy0}+\frac{1}{6}\Delta y^3 u_{yyy0}+O\left(r^4\right),\\
u_5&=u_0-\Delta x u_{x0}-\Delta y u_{y0}+\frac{1}{2}\Delta x^2 u_{xx0}+\Delta x\Delta y u_{xy}+\frac{1}{2}\Delta y^2 u_{yy0}-\frac{1}{3!}\sum_{i=0}^{3}C_3^i\Delta x^i\Delta y^{3-i}u_{x^iy^{3-i}}+O\left(r^4\right),\\
u_6&=u_0+\Delta x u_{x0}-\Delta y u_{y0}+\frac{1}{2}\Delta x^2 u_{xx0}-\Delta x\Delta y u_{xy}+\frac{1}{2}\Delta y^2 u_{yy0}+\frac{1}{3!}\sum_{i=0}^{3}\left[\left(-1\right)^{3-i}C_3^i\Delta x^i\Delta y^{3-i}u_{x^iy^{3-i}}\right]+O\left(r^4\right),\\
u_7&=u_0-\Delta x u_{x0}+\Delta y u_{y0}+\frac{1}{2}\Delta x^2 u_{xx0}-\Delta x\Delta y u_{xy}+\frac{1}{2}\Delta y^2 u_{yy0}+\frac{1}{3!}\sum_{i=0}^{3}\left[\left(-1\right)^{i}C_3^i\Delta x^i\Delta y^{3-i}u_{x^iy^{3-i}}\right]+O\left(r^4\right),\\
u_8&=u_0+\Delta x u_{x0}+\Delta y u_{y0}+\frac{1}{2}\Delta x^2 u_{xx0}+\Delta x\Delta y u_{xy}\frac{1}{2}\Delta y^2 u_{yy0}+\frac{1}{3!}\sum_{i=0}^{3}C_3^i\Delta x^i\Delta y^{3-i}u_{x^iy^{3-i}}+O\left(r^4\right)
\end{aligned} \tag{52}$$

By introducing Eq. (51) and Eq. (52) into Eq. (50), Eq. (50) can be rewritten as:

$$\begin{aligned}
\sum_{j=1}^{8}m_{3j}\left(u_j-u_0\right)&=\left(m_{32}+m_{36}+m_{38}-m_{31}-m_{35}-m_{37}\right)\Delta x u_{x0}\\
&+\left(m_{34}+m_{37}+m_{38}-m_{33}-m_{35}-m_{36}\right)\Delta y u_{y0}\\
&+\frac{1}{2}\left(m_{31}+m_{32}+m_{35}+m_{36}+m_{37}+m_{38}\right)\Delta x^2 u_{xx0}\\
&+\frac{1}{2}\left(m_{33}+m_{34}+m_{35}+m_{36}+m_{37}+m_{38}\right)\Delta y^2 u_{yy0}+O\left(r^4\right)\\
&=u_{xx0}+O\left(r^2\right)
\end{aligned} \tag{53}$$

Then obtain:

$$u_{xx,0}=\sum_{j=1}^{8}m_{3j}\left(u_j-u_0\right)+O\left(r^2\right) \tag{54}$$

Similarly, it is obtained that

$$u_{yy,0}=\sum_{j=1}^{8}m_{4j}\left(u_j-u_0\right)+O\left(r^2\right) \tag{55}$$

For the first-order spatial derivative, the generalized difference expression of $u_{x0}$ is

$$u_{x0}=\sum_{j=1}^{8}m_{1j}\left(u_j-u_0\right) \tag{56}$$

where

$$\begin{aligned}&\left(m_{11},m_{12},m_{13},m_{14},m_{15},m_{16},m_{17},m_{18}\right)\\ &=\left(\frac{-0.4808}{\Delta x},\frac{0.4808}{\Delta x},0,0,\frac{-0.0096}{\Delta x},\frac{0.0096}{\Delta x},\frac{-0.0096}{\Delta x},\frac{0.0096}{\Delta x}\right)\end{aligned} \tag{57}$$

By introducing Eq. (52) and Eq. (57) into Eq. (56), Eq. (56) can be rewritten as:

$$\sum_{j=1}^{8} m_{1j}\left(u_j - u_0\right) = u_{x0} + O\left(\Delta x^2\right) \tag{58}$$

As seen from Eqs. (54), (55), and (58), in the case of the local point cloud shown in Fig. 5 (b), GFDM has the same second-order accuracy as the traditional FDM for the first-order and second-order spatial derivatives.

(*ii*) Example 2: Remove node 8 in Fig. 4 (c) to form a local point cloud in Fig. 5 (c) with poor symmetry. At this time, it can be calculated as follows:

$$\sum_{j=1}^{7} m_{3j}\left(u_j - u_0\right) \tag{59}$$

in which

$$\left(m_{31}, m_{32}, m_{33}, m_{34}, m_{35}, m_{36}, m_{37}\right) = \left(\begin{array}{l} \dfrac{9.6262\times10^{-1}}{\Delta x^2}, \dfrac{9.8754\times10^{-1}}{\Delta x^2}, -\dfrac{3.7377\times10^{-2}}{\Delta x^2}, \\ -\dfrac{1.2459\times10^{-2}}{\Delta x^2}, \dfrac{2.4918\times10^{-2}}{\Delta x^2}, \dfrac{1.2459\times10^{-2}}{\Delta x^2}, \dfrac{1.2459\times10^{-2}}{\Delta x^2} \end{array}\right) \tag{60}$$

Combining the Taylor expansion in Eq. (52), we can get:

$$\sum_{j=1}^{7} m_{3j}\left(u_j - u_0\right) = u_{xx,0} + \frac{1}{3!}\frac{2.4918\times10^{-2}}{\Delta x^2}\left(3\Delta x^2 \Delta y u_{xxy,0} + 3\Delta x \Delta y^2 u_{xyy,0}\right) + O\left(r^2\right) \tag{61}$$

Since Eq. (59) and Eq. (60) are obtained when $\Delta x = \Delta y$, rewrite Eq. (61) as:

$$\sum_{j=1}^{7} m_{3j}\left(u_j - u_0\right) = u_{xx,0} + \frac{2.4918\times10^{-2}}{3!}\left(3\Delta y + 3\frac{\Delta y}{\Delta x}\Delta y\right) + O\left(r^2\right) = u_{xx,0} + O\left(r\right) \tag{62}$$

So the generalized finite difference expression in Eq. (62) has only first-order accuracy, not second-order accuracy. Similarly, for the first-order spatial derivative, we can get:

$$u_{x0} = \sum_{j=1}^{7} m_{1j}\left(u_j - u_0\right) \tag{63}$$

in which

$$\left(m_{11}, m_{12}, m_{13}, m_{14}, m_{15}, m_{16}, m_{17}\right) = \left(\begin{array}{l} \dfrac{-4.8107\times10^{-1}}{\Delta x}, \dfrac{4.9353\times10^{-1}}{\Delta x}, \dfrac{-2.3880\times10^{-4}}{\Delta x}, \dfrac{1.2698\times10^{-2}}{\Delta x} \\ , \dfrac{-6.2296\times10^{-3}}{\Delta x}, \dfrac{6.4684\times10^{-3}}{\Delta x}, \dfrac{-1.2698\times10^{-2}}{\Delta x} \end{array}\right) \tag{64}$$

Then Eq. (65) is obtained. In the case of the local point cloud, the generalized difference expression of the first-order spatial derivative already contains the first-order term, but the coefficients of the first-order terms are small, which can be approximately regarded as having second-order accuracy.

$$\sum_{j=1}^{8} m_{1j}\left(u_j - u_0\right) = u_{x0} - 3.4694\times10^{-18} u_{xx0} - 1.7347\times10^{-18} u_{xy0} + O\left(\Delta x^2\right) \tag{65}$$

The difference in the convergence order of the generalized difference expressions between the two local point clouds in Fig. 5 (b) and Fig. 5 (c) shows that the convergence order of GFDM is less than or equal to that of the traditional FDM. When ignoring the influence of the first-order upwind scheme on the convergence order of the discrete scheme, the FVM commonly used in reservoir numerical simulation has second-order accuracy for the diffusion term expressed by the second-order spatial derivative, which can be regarded as the integral form of the traditional FDM to a certain extent. Of course, many high-order FVMs have been developed in computational fluid dynamics, but because the size of the reservoir calculation domain is relatively large, the grid size is generally about 10m or even larger, so high-order methods are rarely applied to reservoir numerical simulation. The proposed NCDMM is essentially an integral form of GFDM, which can be regarded as an extended finite volume method (EFVM). Therefore, NCDMM has the same convergence order of diffusion term as GFDM, that is, the convergence order of NCDMM is also less than or equal to the second order.

## 3. Results of numerical analysis

In this section, four numerical test cases are implemented to analyze the computational performances of the proposed NCDMM.

3.1 Two-phase porous flow in a rectangular heterogeneous domain

As mentioned above, injection and production wells exist in reservoir numerical simulation, and these wells are regarded as the point source or sink of the reservoir domain, Therefore, this example compares the calculation results of NCDMM with those of FVM commonly used in reservoir numerical simulation to illustrate that the proposed NCDMM can handle the point source problem in reservoir numerical simulation with high accuracy. Fig. 6 (a) shows the strongly heterogeneous permeability profile, which is in the base-10 logarithm. As shown in Fig. 6 (b), there are a production well and an injection well in the reservoir with constant flow rates of 60 $m^3$/d.The size of the reservoir domain in this example is 600m×180m, and the thickness is 3m (Although this work focuses on the 2D problem, because the unit of the well flow rate in reservoir flow problems is generally $m^3$/d, it is necessary to give the thickness of the reservoir model in the design of the numerical example, so that the strength of the source or sink term on the 2D domain can be obtained by dividing the well flow rate by the thickness), and the reservoir boundary is closed. The basic physical properties and oil-water two-phase relative permeability data are given in Table 1 and Table. 2, respectively. Fig. 6 (c), (d), and (e) show the Cartesian mesh, Cartesian and irregular node collocations in the reservoir calculation domain respectively, in which the black solid points are the boundary and internal nodes of the calculation domain, The red diamond points are virtual nodes added to calculate the node control volumes and deal with the derivative boundary condition. Based on the Cartesian mesh, the FVM with the two-point flux approximation (TPFA) scheme [40, 41] commonly used in reservoir numerical simulation is employed to provide the reference solution.

For the Cartesian collocation, different radii (10m, 15m, 20m, and 25m) of the node influence domain are used to determine the *i-j* pairs (i.e., connectable point cloud) and calculate corresponding node control volumes. For the irregular collocation, which comes from the mesh vertices of the triangulation of the reservoir domain, two methods are used to construct the connectable point cloud, one is to make the radius of the influence domain of each node 7.5m, and the other is to use the triangulation-based method given in Section 2.4. For the convenience of narration, these two connectable point clouds for the irregular collocation are denoted as $\mathbf{NP}_1$ and $\mathbf{NP}_2$, respectively.

To verify the empirical method of calculating node control volumes (Not multiplied by the domain thickness) in Section 2.3, Fig. 7 compares the calculated profiles of node control volumes when adopting weight function $w_1$ without Step 3 in Section 2.3 (simply denoted by $w_1$ in this paper), weight function $w_2$ without Step 3 in Section 2.3 (simply denoted by $w_2$ in this paper), and weight function $w_2$ with Step 3 in Section 2.3 (i.e., the method in Section 2.3, simply denoted by weighted $w_2$ in this paper), and Fig. 8 compares the colorful distribution maps of node control volumes for different radii of node influence domain when adopting $w_1$, $w_2$, and weighted $w_2$, and Fig. 9 compares the colorful distribution maps of node control volumes for $\mathbf{NP}_1$ and $\mathbf{NP}_2$ when adopting $w_1$, $w_2$, and weighted $w_2$.

As seen in Fig. 7 and Fig. 8, with the increase of the radius of the node influence domain, the distribution of node control volumes calculated by $w_1$ will become more and more uneven. Using $w_2$ will significantly reduce the unevenness of the distribution of node control volumes, and using weighted $w_2$ can achieve a better effect of reducing the unevenness than $w_2$. According to the brief analysis in Section 2.3, for Cartesian coordinates, even when the radius of the node influence domain increases to 25m, the local point cloud of the two neighboring nodes away from the boundary should be the same, so the control volume of the two nodes should be the same, which indicates that the reasonable distribution of the node control volume, in this case, should have quite good uniformity. It can also be seen from Fig. 9 that for the irregular collocation, compared with using $w_1$, using $w_2$ or weighted $w_2$ can obtain significantly more even node control volume distribution in the area with even node distribution density, and can maintain a reasonably small node control volume in the area with dense node density. The subsequent comparison of simulation results will further illustrate this issue. It should be pointed out that in this example, compared with $\mathbf{NP}_1$, $\mathbf{NP}_2$ can get a more even node control volume distribution, and the subsequent simulation results will also verify that the calculation accuracy of $\mathbf{NP}_2$ is higher than that of $\mathbf{NP}_1$. This verifies the analysis in Sections 2.3 and 2.4 that in the case of point clouds with uneven node density, the same-radius node influence domain in $\mathbf{NP}_1$ will much increase the number of nodes in the local point cloud of the node in the area with high node density, resulting in the increase of the difference of $\left(m_{3j}^{i}+m_{4j}^{i}\right)$ and $\left(m_{3i}^{j}+m_{4i}^{j}\right)$, then reducing the calculation accuracy of node control volumes and simulation results, which can be avoided by $\mathbf{NP}_2$ based on the triangulation-based method in Section 2.4.

Figs. 10-13 compare the pressure and oil saturation profiles calculated by using different connectable point clouds, and Fig. 14 compares the corresponding $L_2$ errors, which are calculated via Eq. (66).

$$error_p=\sqrt{\frac{\sum_{i=1}^{n_p}\left(p_i^{cal}-p_i^{ref}\right)^2}{n_p}},\quad error_{S_w}=\sqrt{\frac{\sum_{i=1}^{n_p}\left(S_{w,i}^{cal}-S_{w,i}^{ref}\right)^2}{n_p}} \tag{66}$$

where the superscript *cal* denotes the calculated value, the superscript *ref* denotes the reference solution and $n_p$ denotes the number of nodes participating in the comparison.

Whether from the intuitive comparison of colorful distribution maps or the comparison of specific $L_2$ error data shown in Fig. 14, The following results can be found:

(*i*) Using weighted $w_2$ can improve the calculation accuracy of node control volume distribution, and then improve the accuracy of simulation results. Except that the calculation accuracy of weighted $w_2$ is slightly lower than that of $w_1$ at 'Cartesian colocation, 10m' and 'Cartesian colocation, 15m', in other connectable point clouds, the calculation accuracy of weighted $w_2$ will be significantly higher than that of $w_1$. In the case of 'Cartesian colon, 10m' and 'Cartesian colon, 15m', the figure also shows that compared with the use of $w_1$, the use of weighted $w_2$ has little impact on the calculation of node control volume distribution, which should be the reason why the calculation accuracy of weighted $w_2$ is slightly lower than that of $w_1$ in the case of these two point clouds. We know that the meshless method is more used in complex geometry, so there are more cases of using irregular point clouds, and the subsequent examples in this section will further fully demonstrate that for general connectable point clouds, using weighted $w_2$ can significantly

(*ii*) in the case of the same collocation, the determination method of connectable point cloud significantly affects the calculation accuracy. Fig. 14 shows that under the same Cartesian collocation, the oil saturation calculation error dominated by dissipation error will increase with the increase of the radius of the node influence domain. In addition, in the case of the same irregular collocation, the calculation error corresponding to $\mathbf{NP}_2$ obtained by using the triangulation-based method in Section 2.4 to determine the connectable point cloud is significantly smaller than the calculation error corresponding to $\mathbf{NP}_1$ determined by setting each node to have the same-radius influence domain, which also indicates that in the case of any collocation, finding a suitable construction method of the connectable point cloud to provide each node with high-quality local point cloud is the key issue to improve the computational performance of proposed NCDMM, and it is a very important future work.

Overall, NCDMM using the empirical method of calculating node control volumes (i.e. weighted $w_2$) given in Section 2.3 can achieve satisfactory calculation accuracy in case of different connectable point clouds.

Table 3 shows the relevant parameters in the nonlinear solver in this section. Fig. 15 compares the cumulative Newton iteration steps vs. time in the nonlinear solution process of different cases. As seen from the figure, for Cartesian collocation, with the increase of the radius of the node influence domain, the consuming Newton iteration steps of the NCDMM will decrease, and all lower than those of the FVM, indicating that with the increase of node influence domain, the more abundant *i-j* pairs, the faster the convergence speed of nonlinear solution. This comparison fully shows that for the cases which both NCDMM and FVM can handle, NCDMM can achieve comparable computational efficiency to traditional FVM.

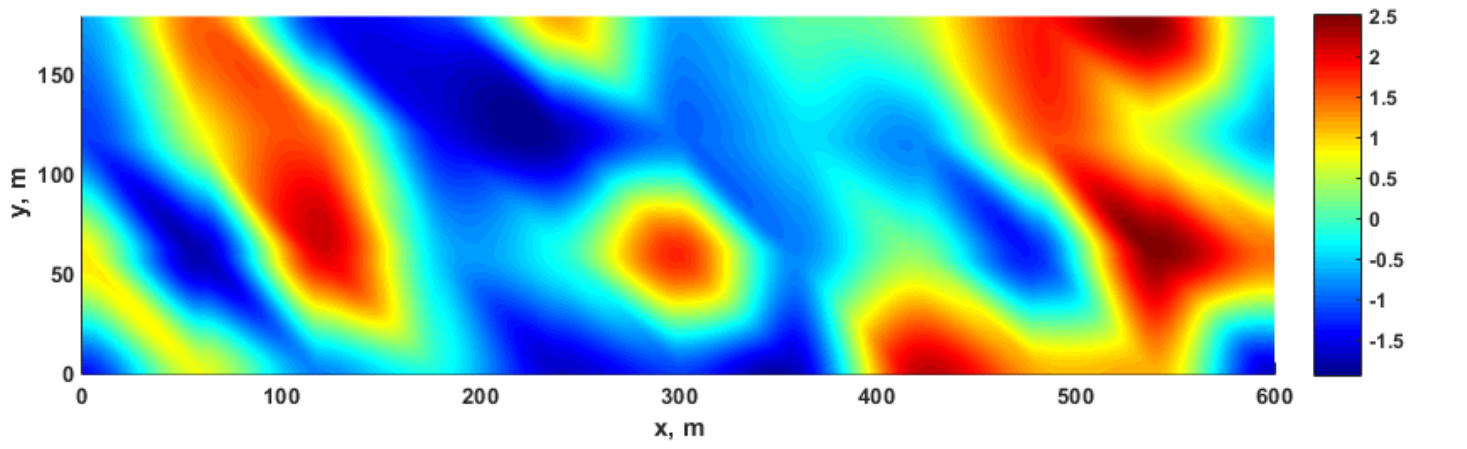

(a) ten-based logarithm of the permeability profile, $\log_{10}(k)$

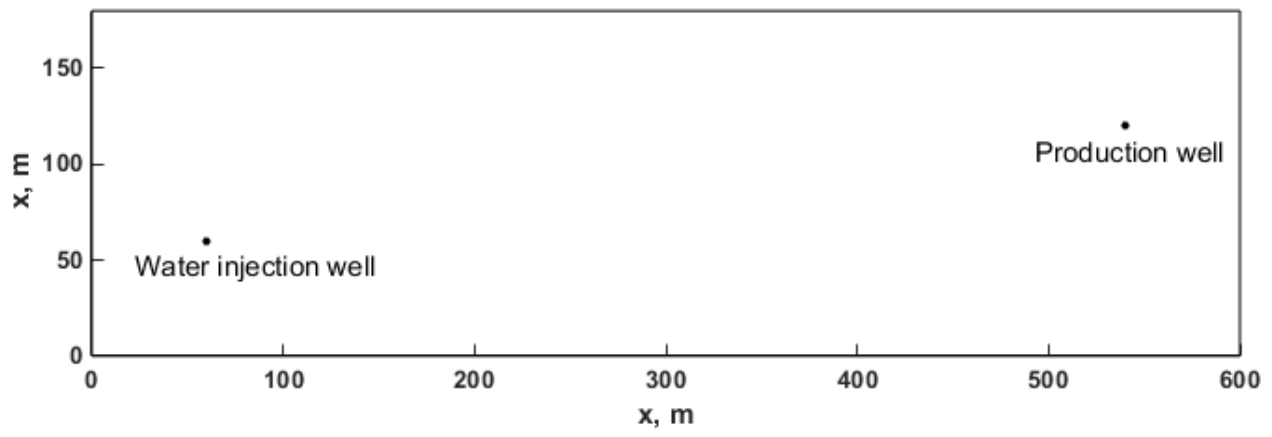

(b) sketch of the reservoir domain

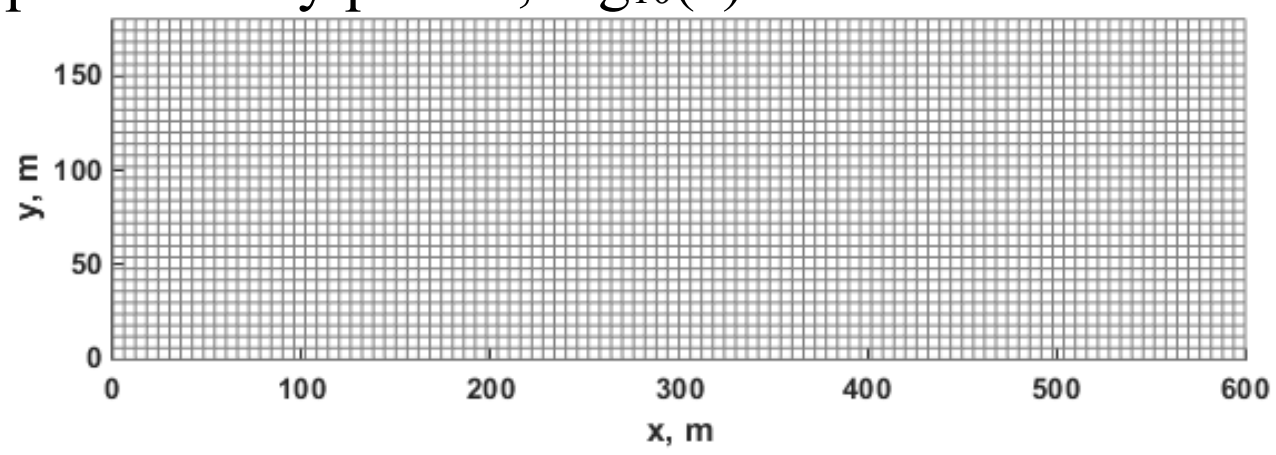

(c) Cartesian mesh

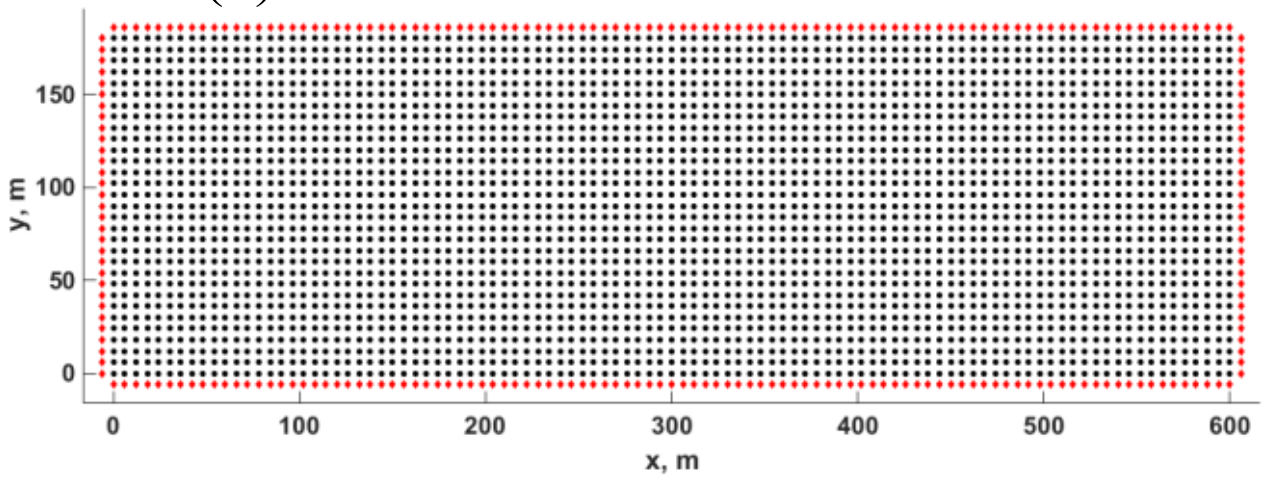

(d) Cartesian collocation

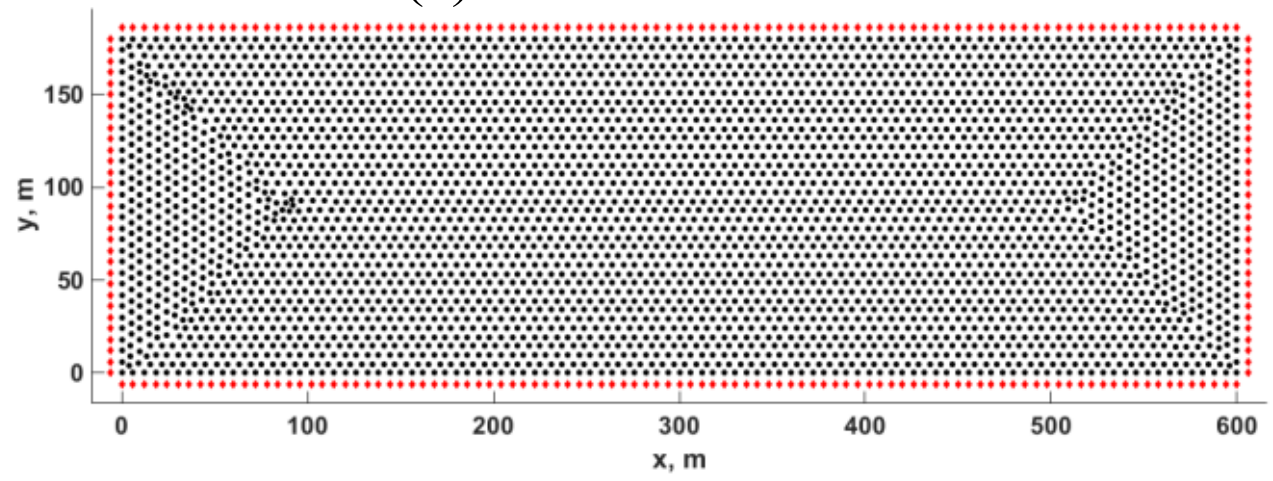

(e) Irregular collocation

Fig. 6 Sketches of the reservoir domain and various discretization

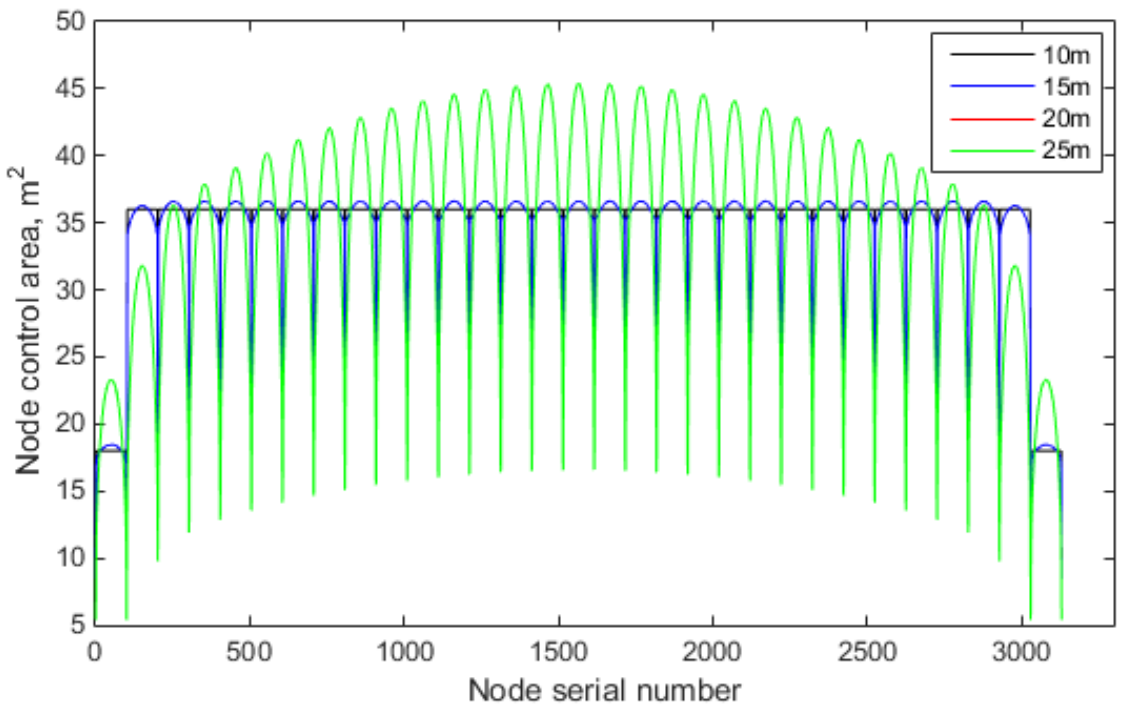

(a) $w_1$

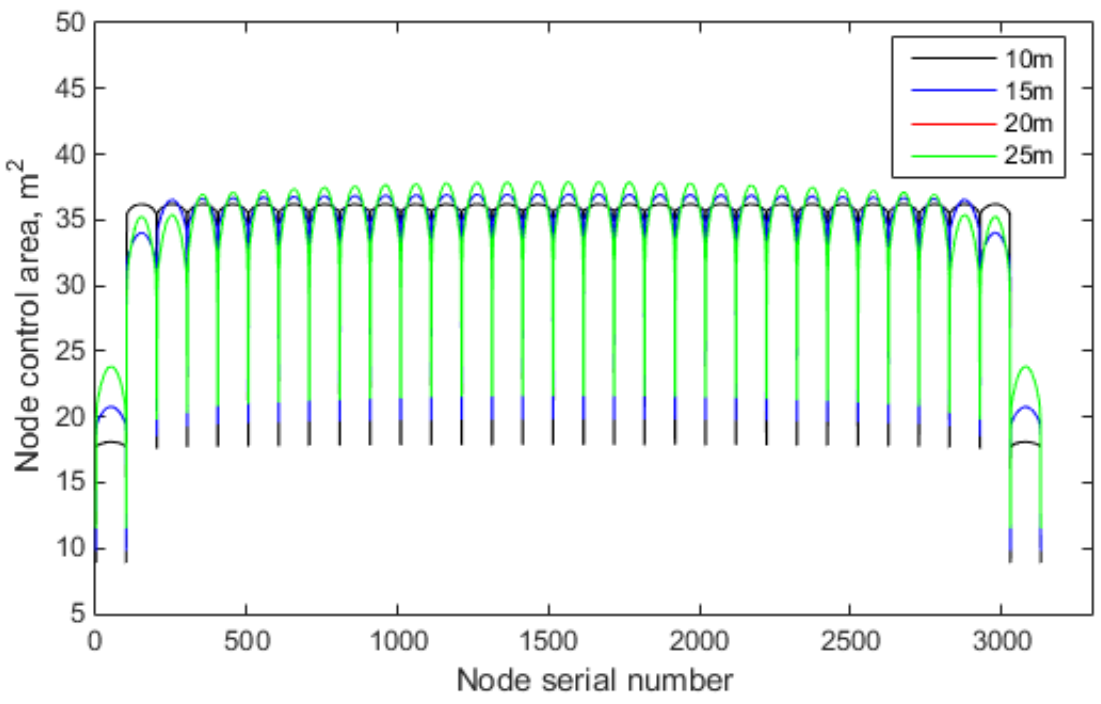

(b) $w_2$

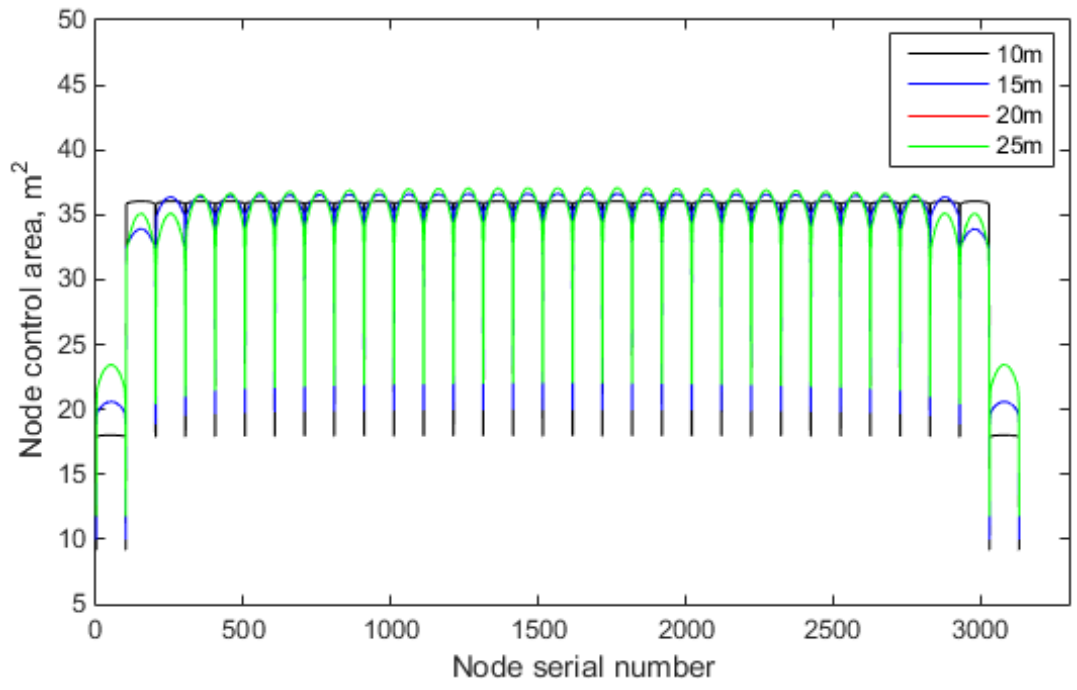

(c) weighted $w_2$

Fig. 7 Comparison of the calculated node control areas in case of different radii of node influence domain, different point clouds, and different weight functions

$w_1$

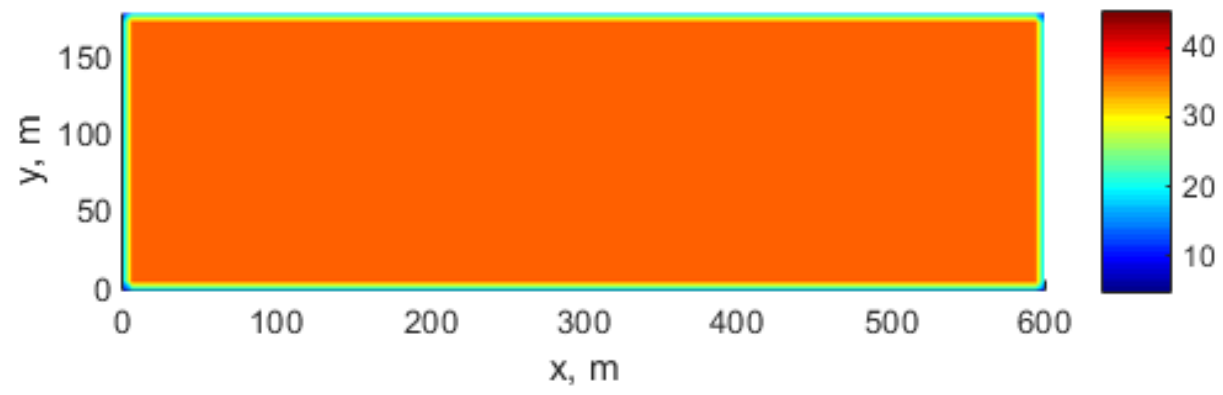

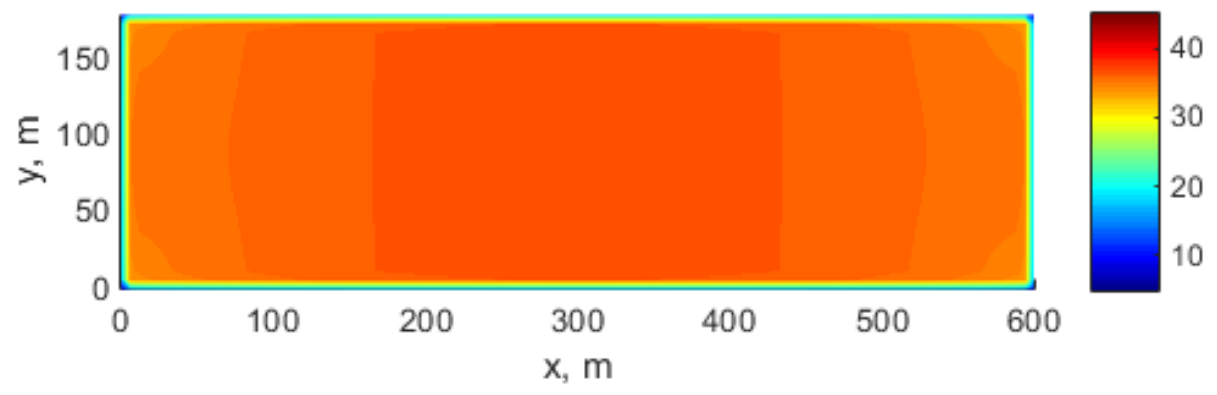

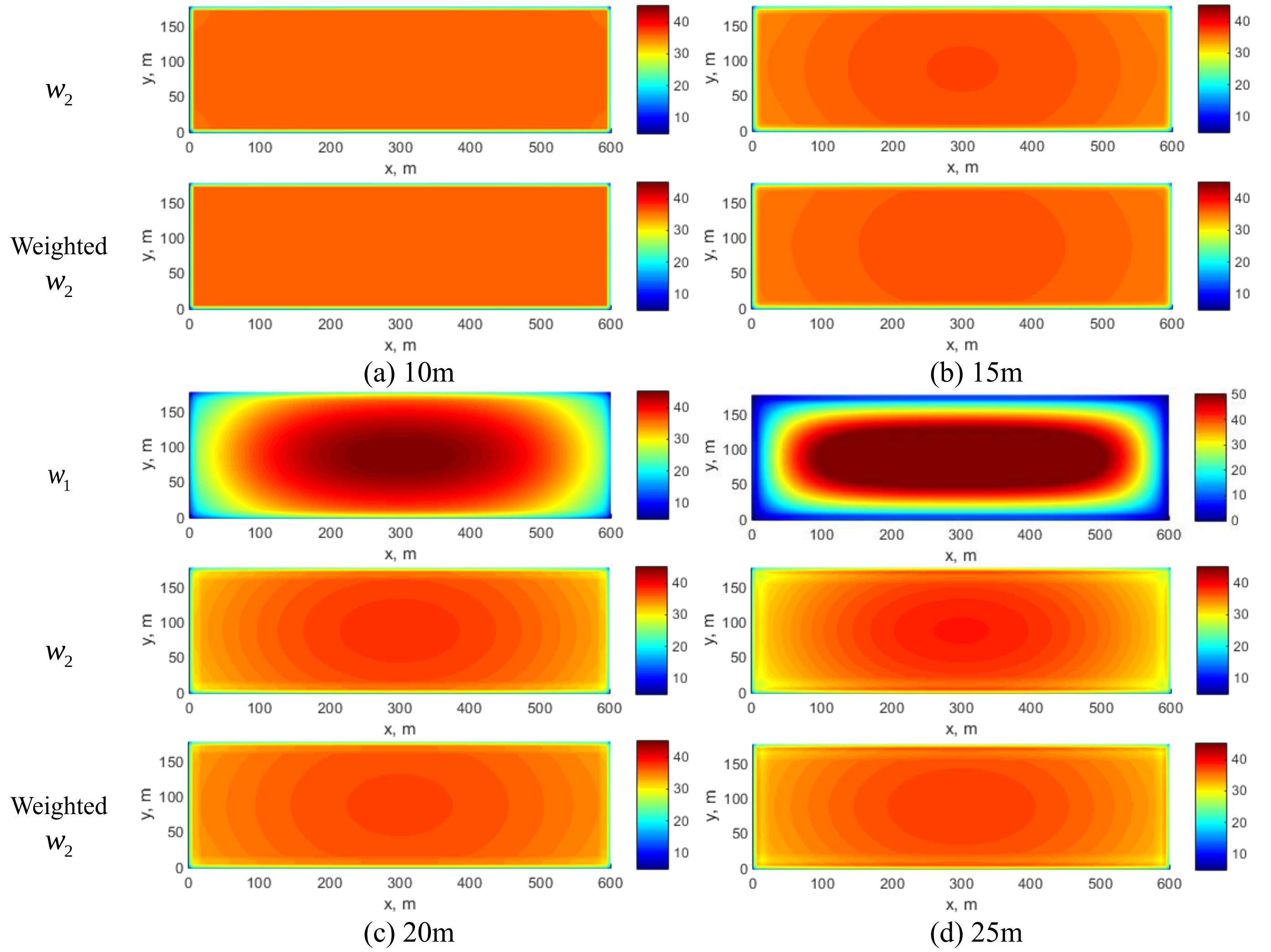

Fig. 8 Comparisons of calculated profiles of node control volumes

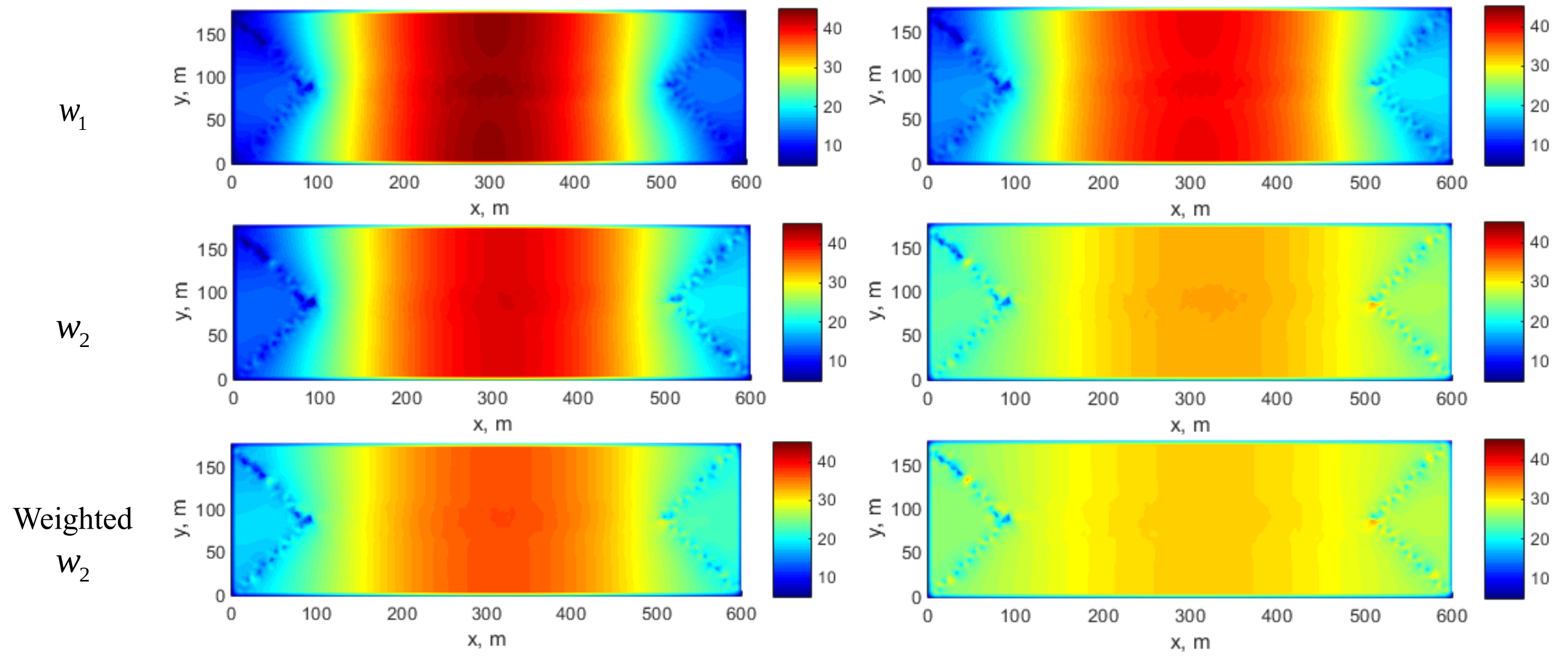

(a) given radius of node influence domain (b) triangulation-based method

Fig. 9 Node-control-volume profiles of the irregular point cloud calculated by different methods

Table. 1 Basic physical properties used in this numerical example

| Properties | Values | Properties | Values |
|---|---|---|---|
| Porosity | 0.2 | Rock compressibility | $1\times10^{-4}$ MPa$^{-1}$ |
| Oil compressibility | $3\times10^{-3}$ MPa$^{-1}$ | Water compressibility | $4\times10^{-4}$ MPa$^{-1}$ |
| Oil viscosity | 2 mPa·s | Water viscosity | 0.6 mPa·s |
| Initial reservoir pressure | 15 MPa | Initial water saturation | 0.20 |
| Oil volume factor | 1.0 | Water volume factor | 1.0 |
| Reservoir thickness | 3 m | Well radius | 0.1 m |
| Skin factor | 0 | | |

Table. 2 Basic physical properties used in this numerical example

| *Sw* | *Krw* | *Kro* | *Sw* | *Krw* | *Kro* |
|---|---|---|---|---|---|
| 0.2 | 0 | 1 | 0.55 | 0.3403 | 0.1736 |
| 0.25 | 0.0069 | 0.8403 | 0.6 | 0.4444 | 0.1111 |
| 0.3 | 0.0278 | 0.6944 | 0.65 | 0.5625 | 0.0625 |
| 0.35 | 0.0625 | 0.5625 | 0.7 | 0.6944 | 0.0278 |
| 0.4 | 0.1111 | 0.4444 | 0.75 | 0.8403 | 0.0069 |
| 0.45 | 0.1736 | 0.3403 | 0.8 | 1 | 0 |
| 0.5 | 0.25 | 0.25 | | | |

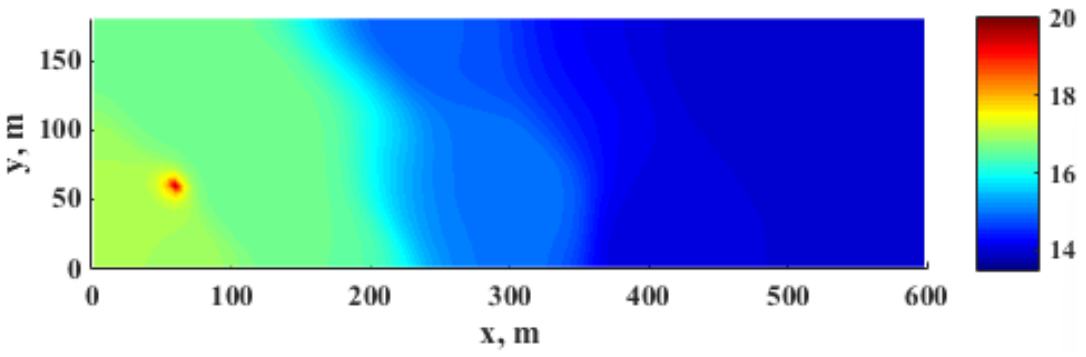

FVM reference solution

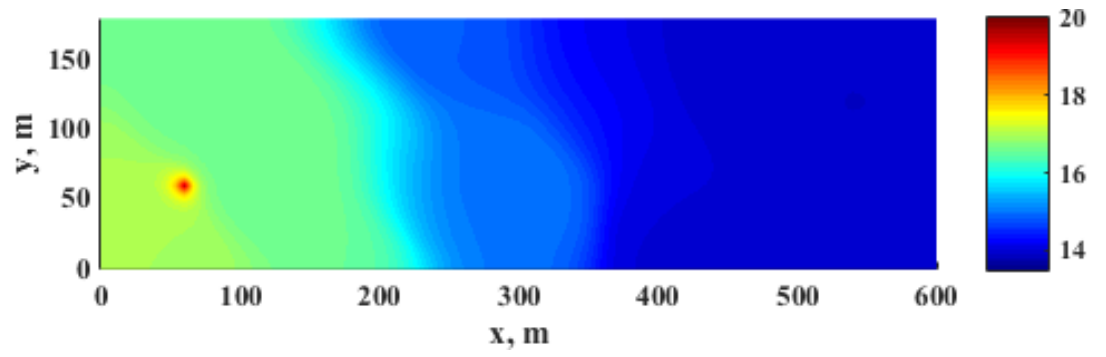

NCDMM, Cartesian, 10m, $w_1$

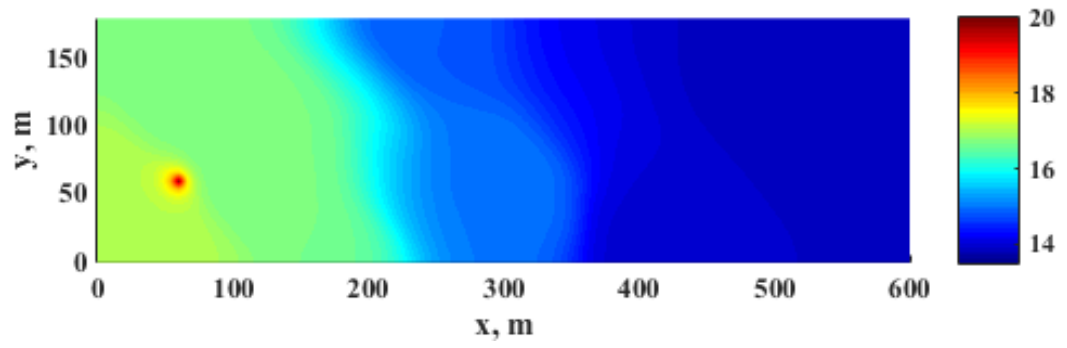

NCDMM, Cartesian, 10m, weighted $w_2$

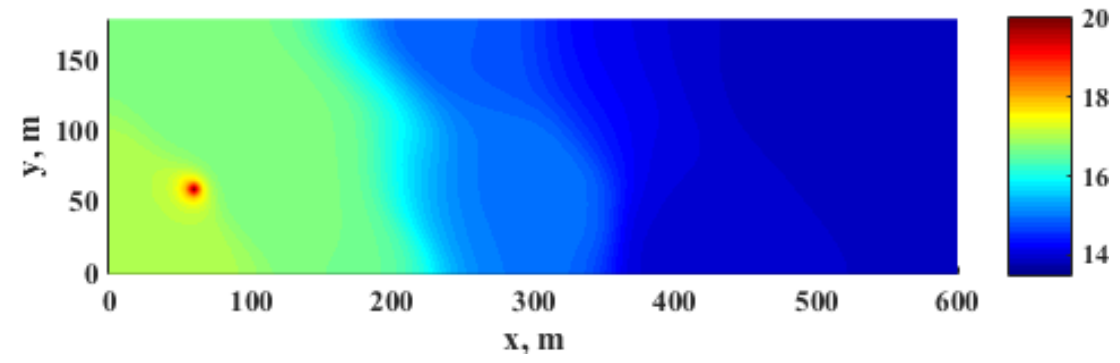

NCDMM, Cartesian, 15m, $w_1$

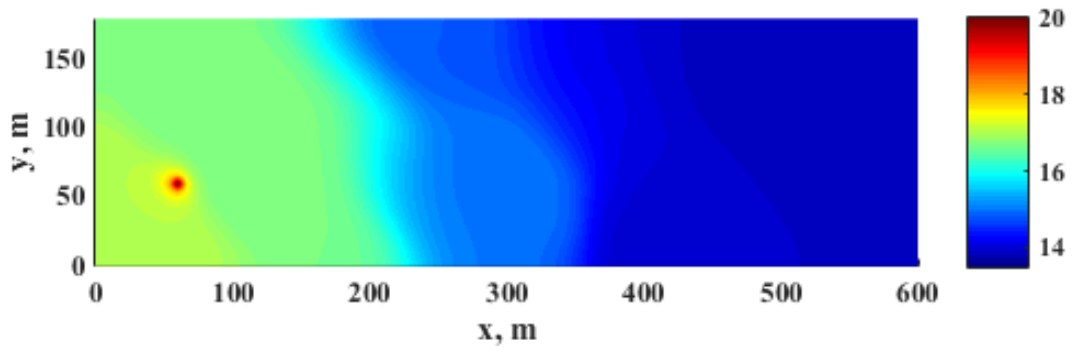

NCDMM, Cartesian, 15m, weighted $w_2$

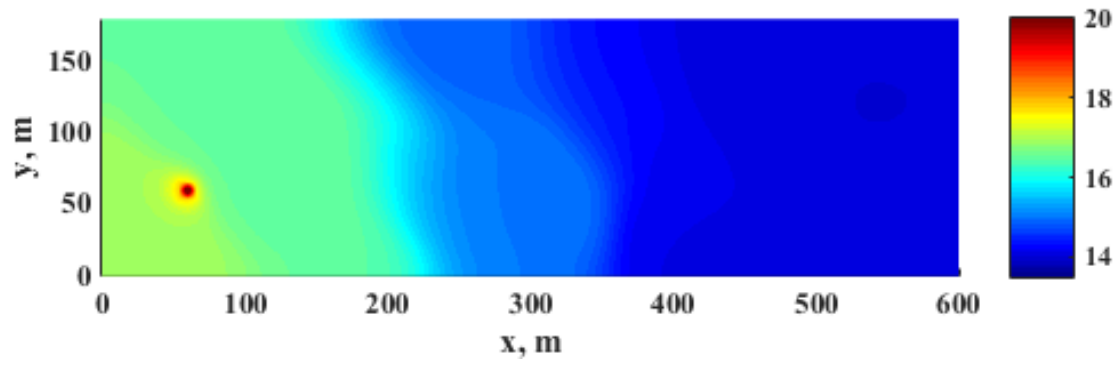

NCDMM, Cartesian, 20m, $w_1$

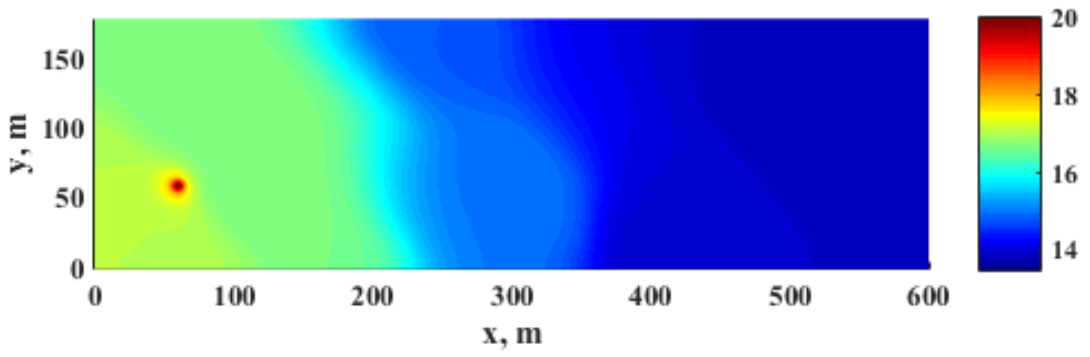

NCDMM, Cartesian, 20m, weighted $w_2$

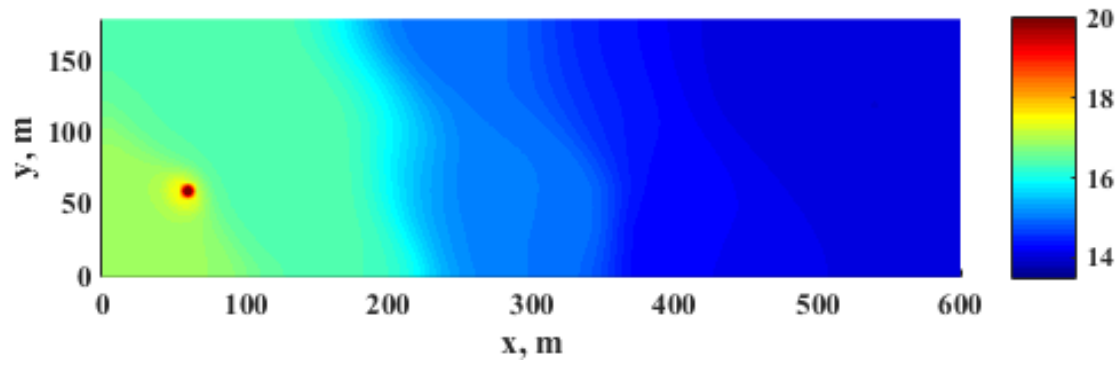

NCDMM, Cartesian, 25m, $w_1$

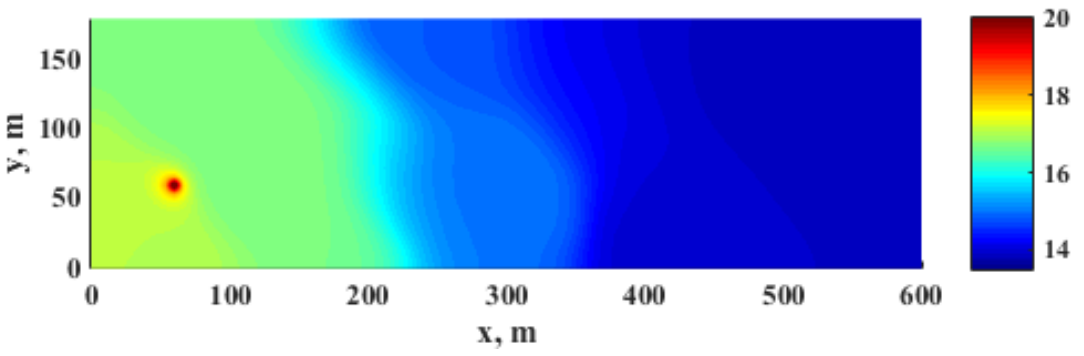

NCDMM, Cartesian, 25m, weighted $w_2$

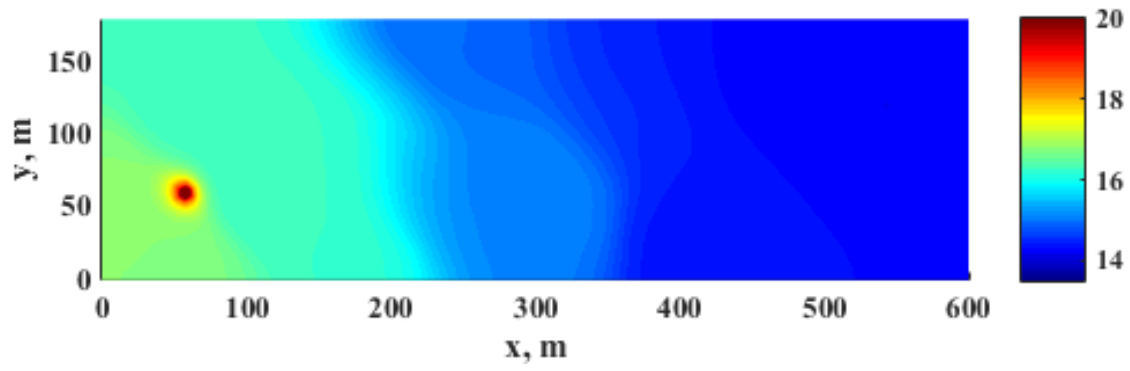

NCDMM, irregular, $\mathbf{NP}_1$, $w_1$

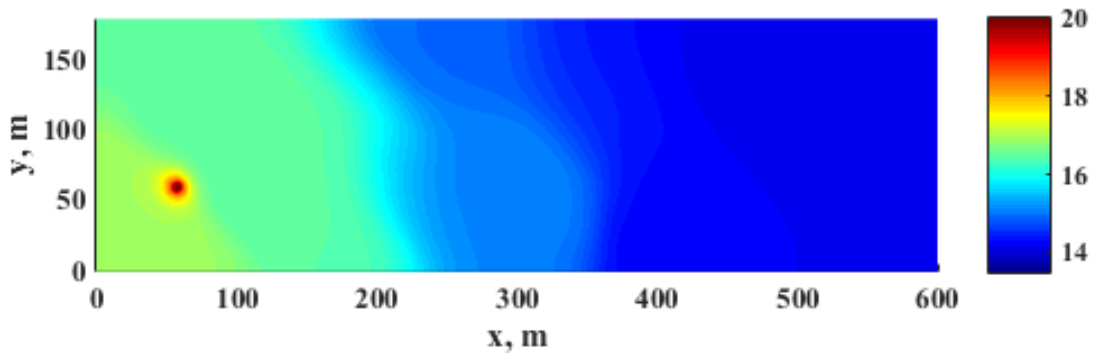

NCDMM, irregular, $\mathbf{NP}_1$, weighted $w_2$

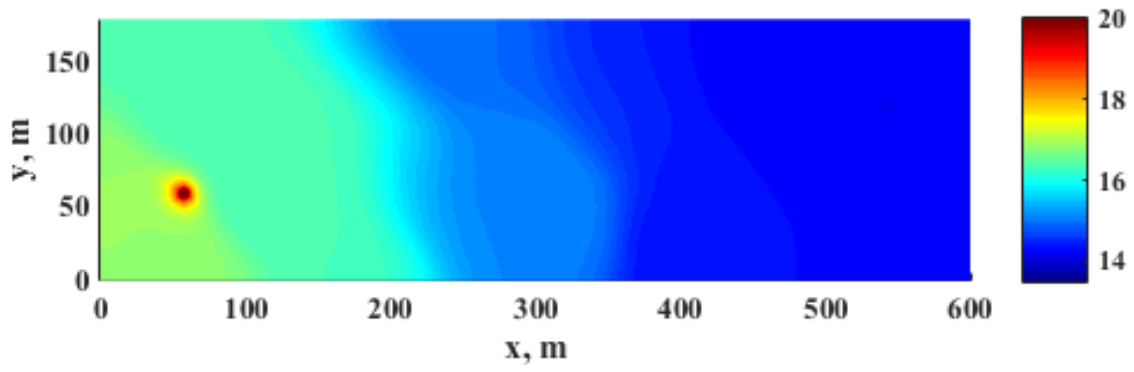

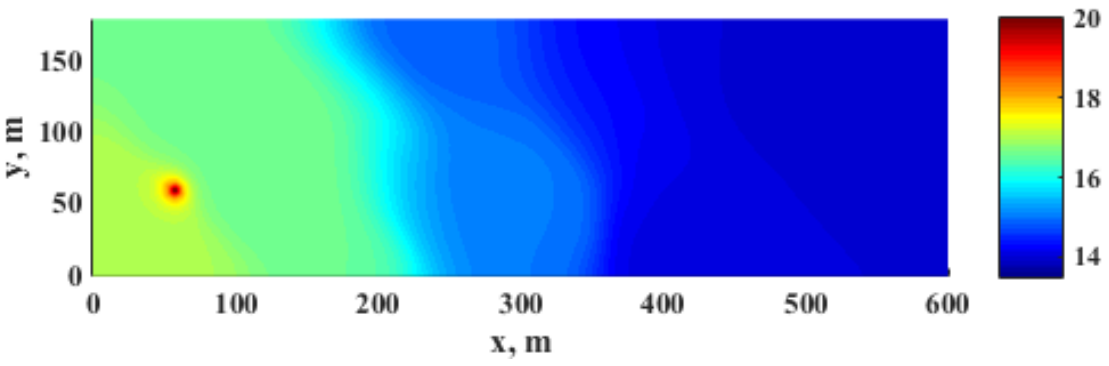

NCDMM, irregular, $\mathbf{NP}_2$, $w_1$ NCDMM, irregular, $\mathbf{NP}_2$, weighted $w_2$

Fig. 10 Calculated pressure profiles at 125 days in different cases

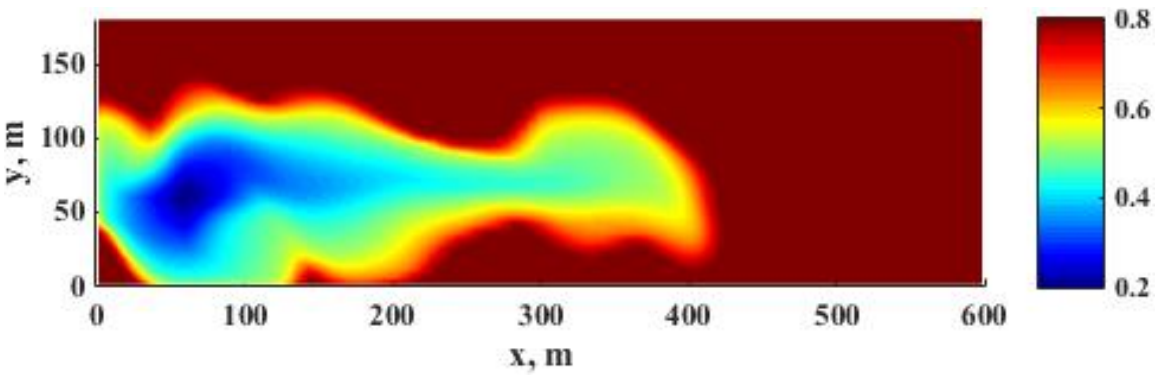

FVM reference solution

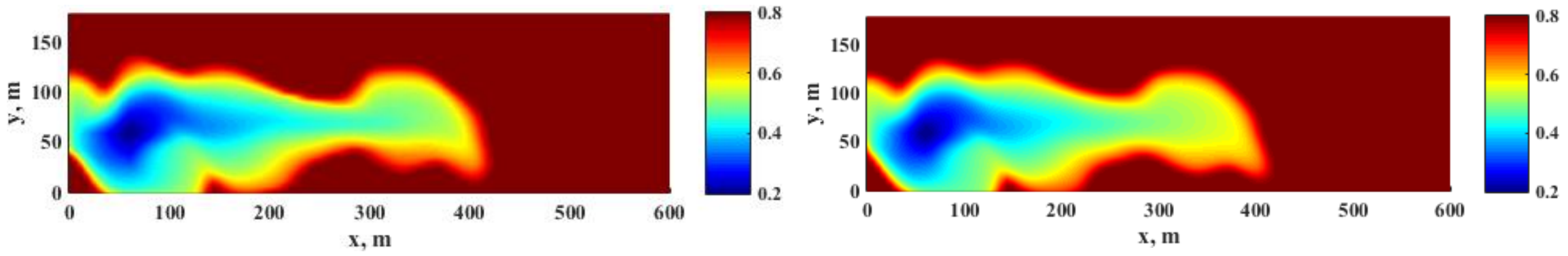

NCDMM, Cartesian, 10m, $w_1$ NCDMM, Cartesian, 10m, weighted $w_2$

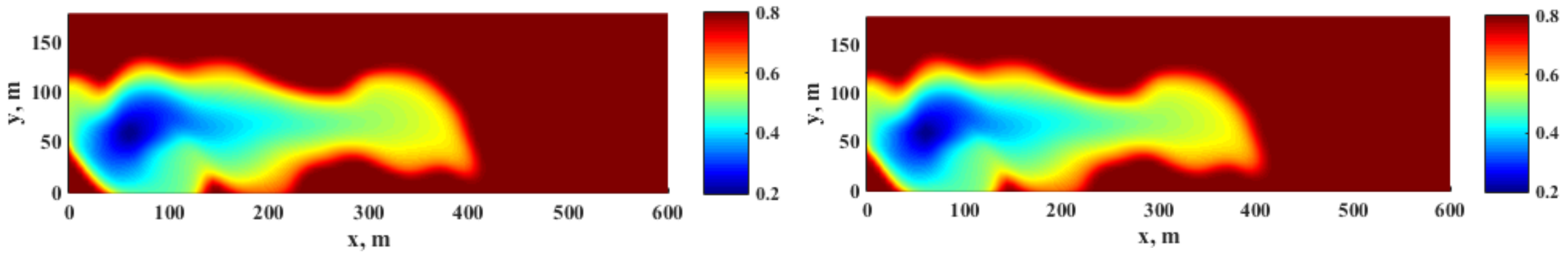

NCDMM, Cartesian, 15m, $w_1$ NCDMM, Cartesian, 15m, weighted $w_2$

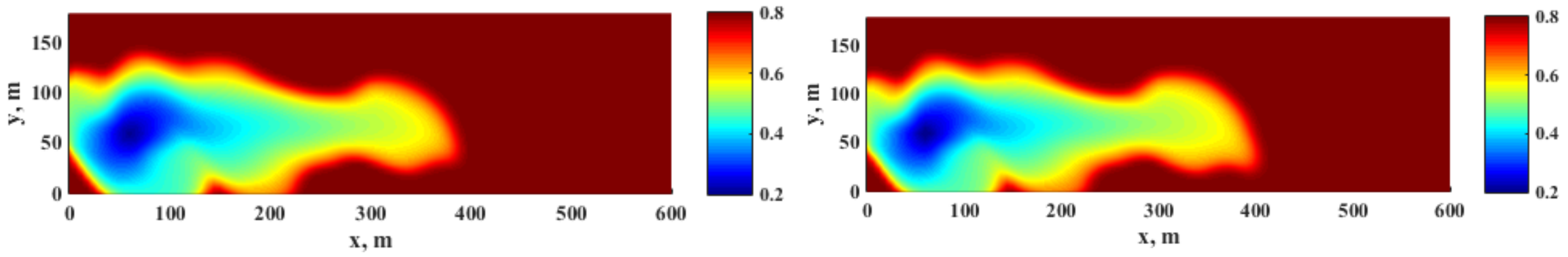

NCDMM, Cartesian, 20m, $w_1$ NCDMM, Cartesian, 20m, weighted $w_2$

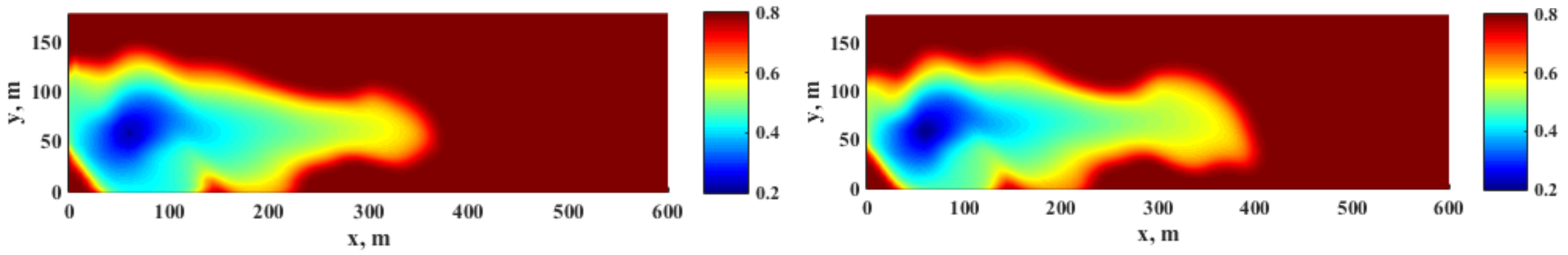

NCDMM, Cartesian, 25m, $w_1$ NCDMM, Cartesian, 25m, weighted $w_2$

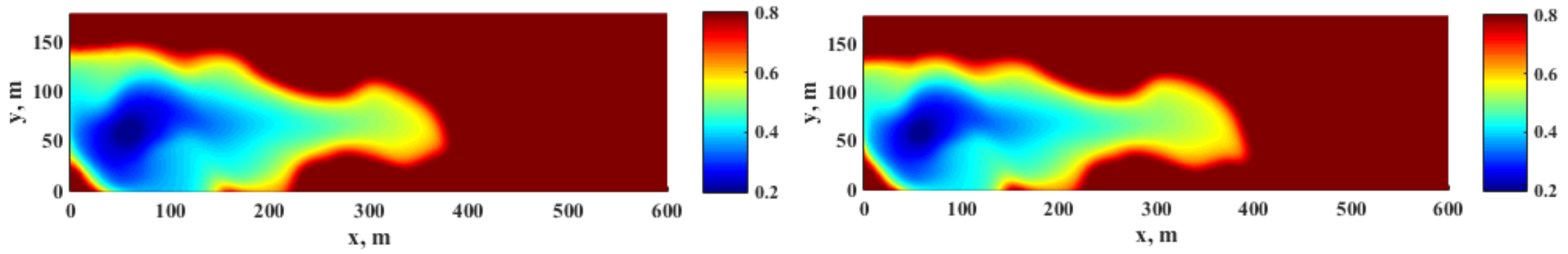

NCDMM, irregular, $\mathbf{NP}_1$, $w_1$ NCDMM, irregular, $\mathbf{NP}_1$, weighted $w_2$

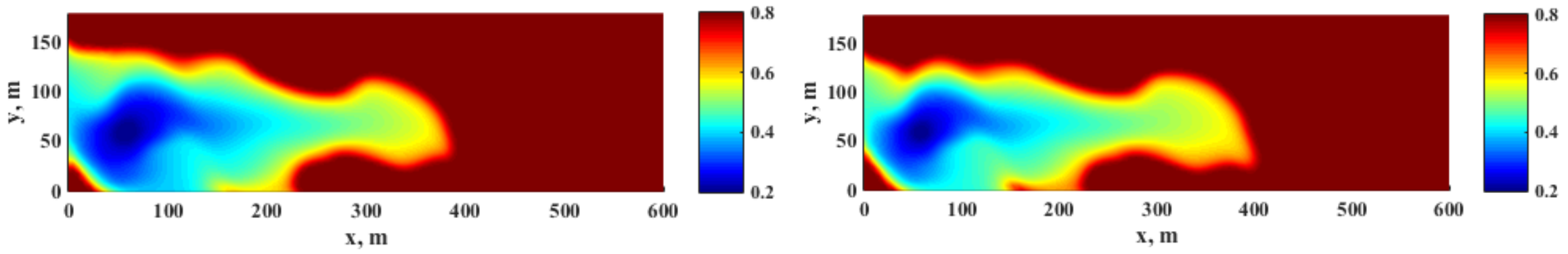

NCDMM, irregular, $\mathbf{NP}_2$, $w_1$ NCDMM, irregular, $\mathbf{NP}_2$, weighted $w_2$

Fig. 11 Calculated oil saturation profiles at 125 days in different cases

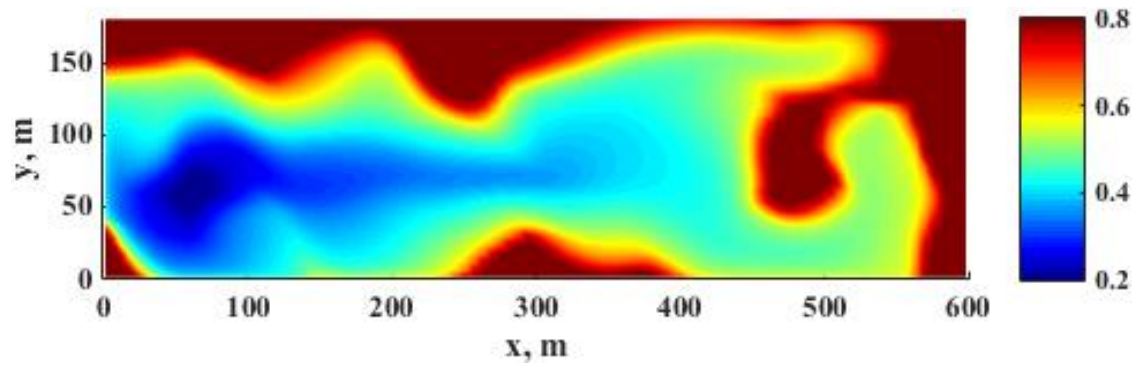

FVM reference solution

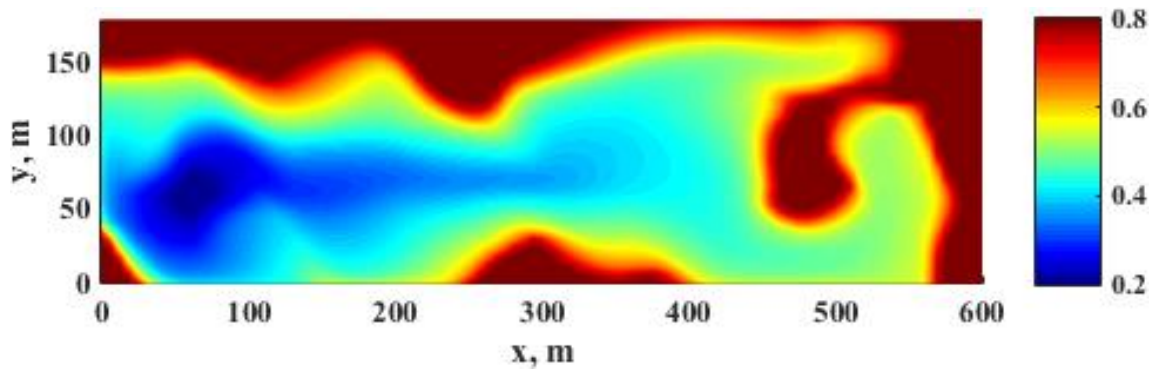

NCDMM, Cartesian, 10m, $w_1$

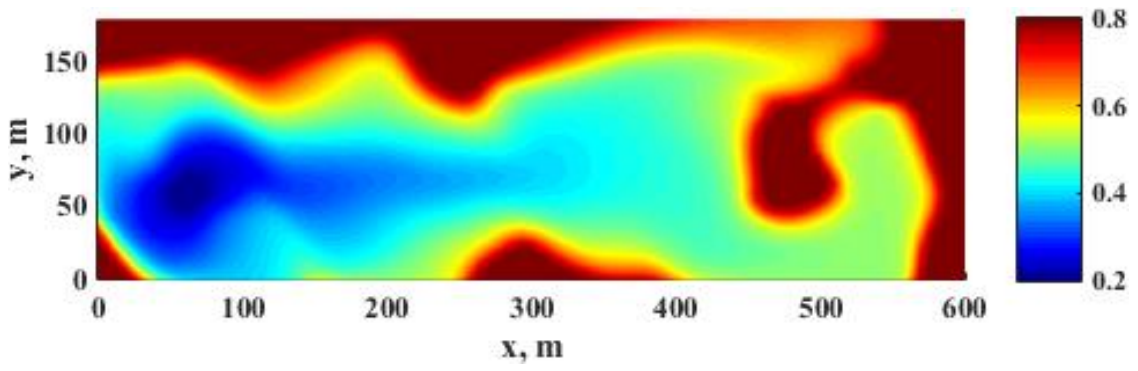

NCDMM, Cartesian, 10m, weighted $w_2$

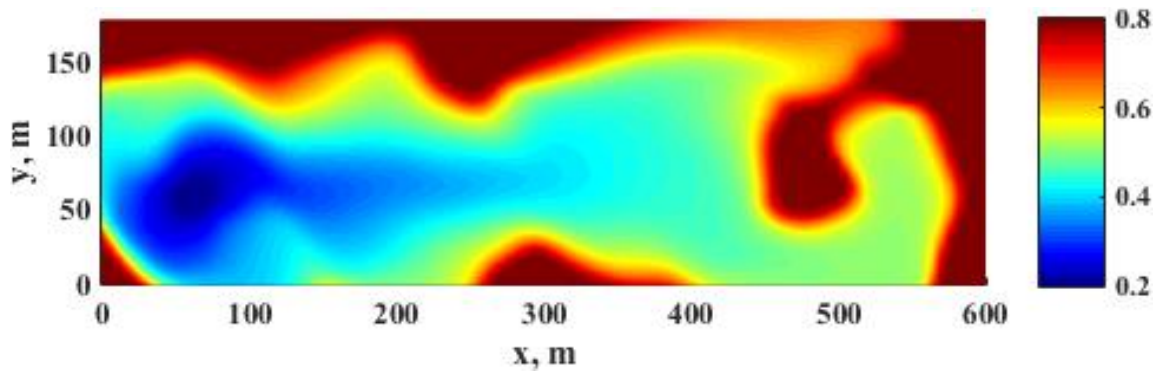

NCDMM, Cartesian, 15m, $w_1$

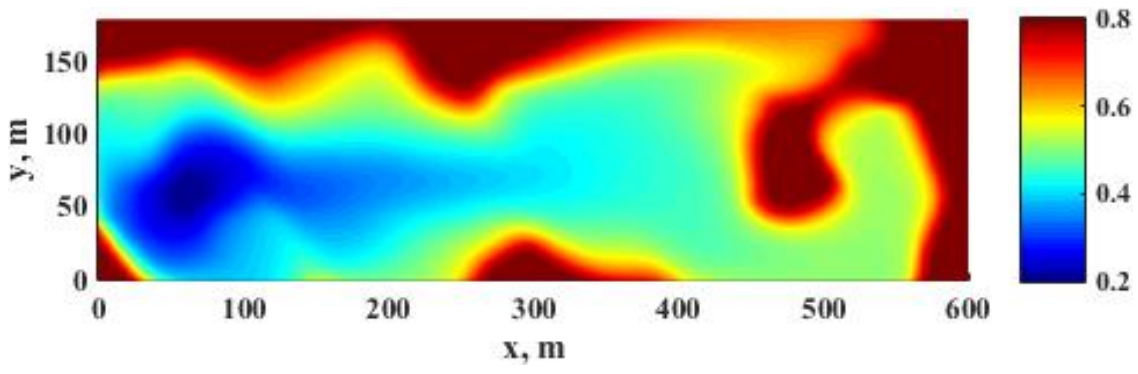

NCDMM, Cartesian, 15m, weighted $w_2$

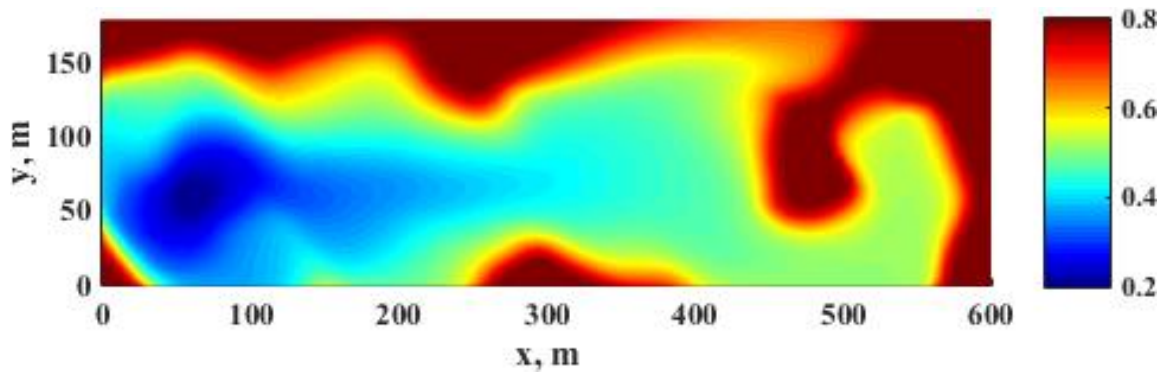

NCDMM, Cartesian, 20m, $w_1$

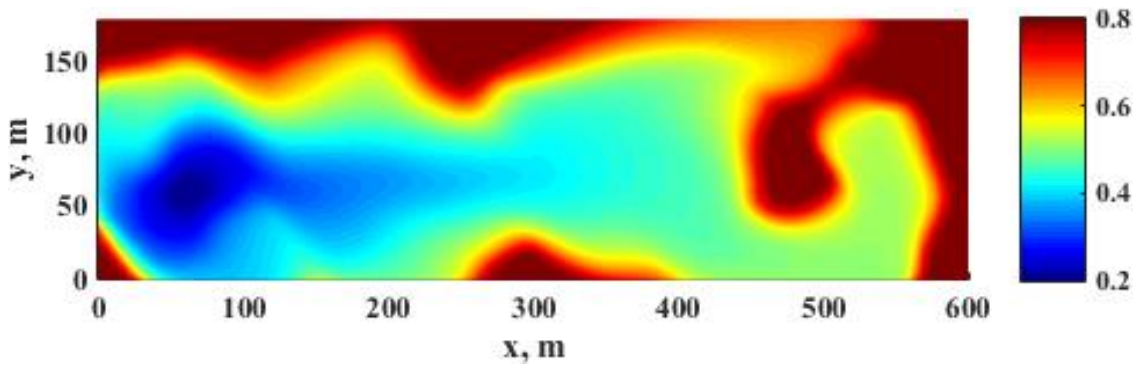

NCDMM, Cartesian, 20m, weighted $w_2$

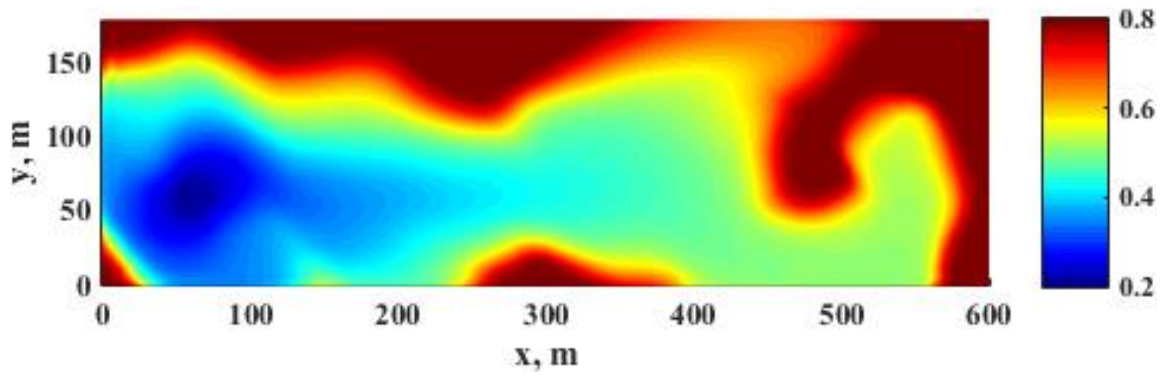

NCDMM, Cartesian, 25m, $w_1$

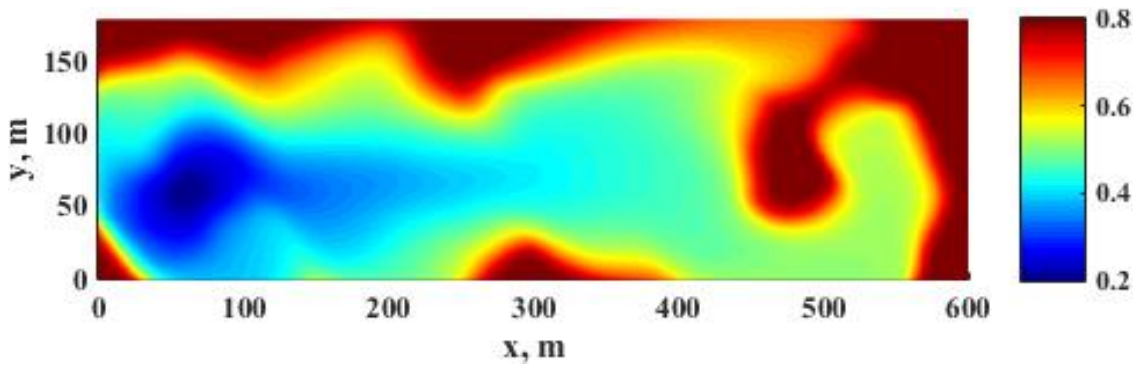

NCDMM, Cartesian, 25m, weighted $w_2$

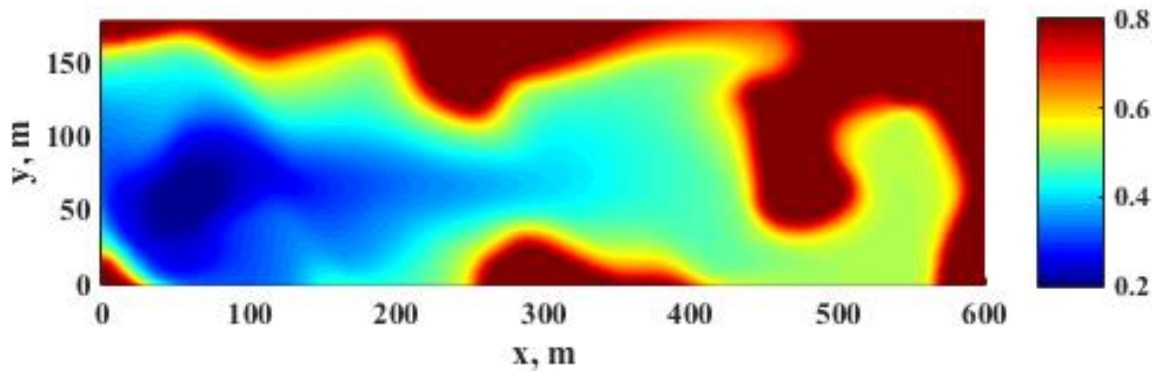

NCDMM, irregular, $\mathbf{NP}_1$, $w_1$

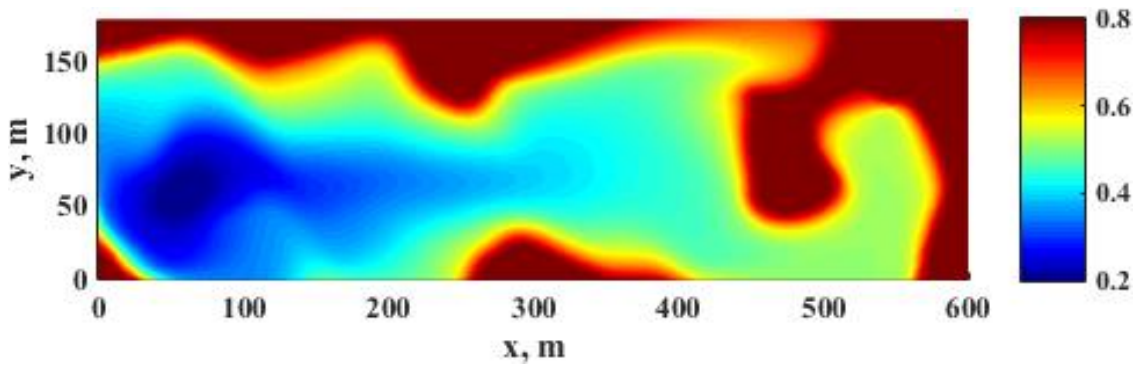

NCDMM, irregular, $\mathbf{NP}_1$, weighted $w_2$

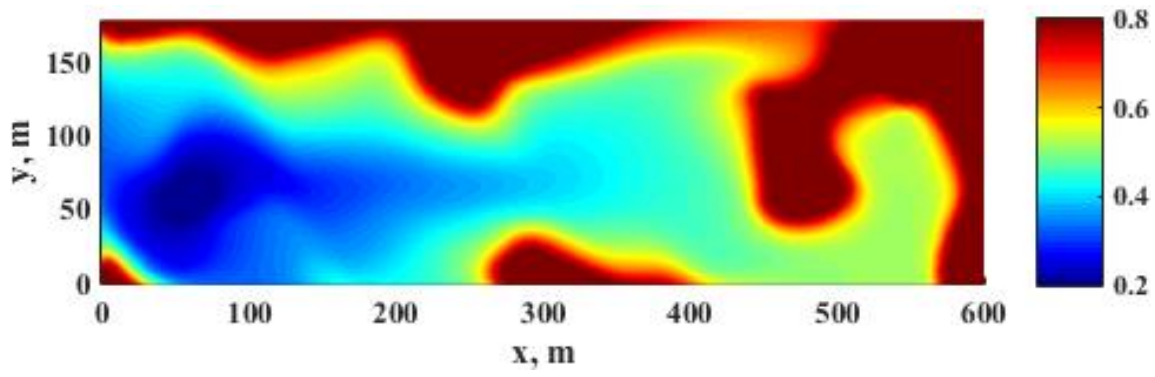

NCDMM, irregular, $\mathbf{NP}_2$, $w_1$

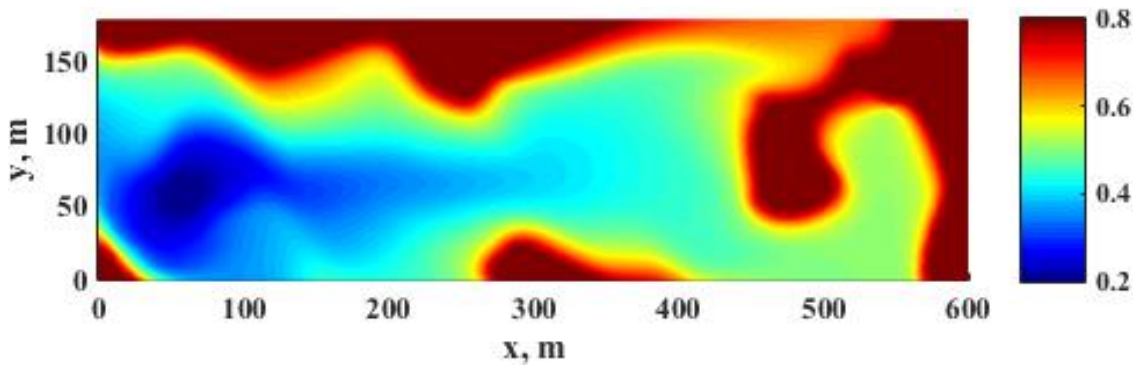

NCDMM, irregular, $\mathbf{NP}_2$, weighted $w_2$

Fig. 12 Calculated oil saturation profiles at 300 days in different cases

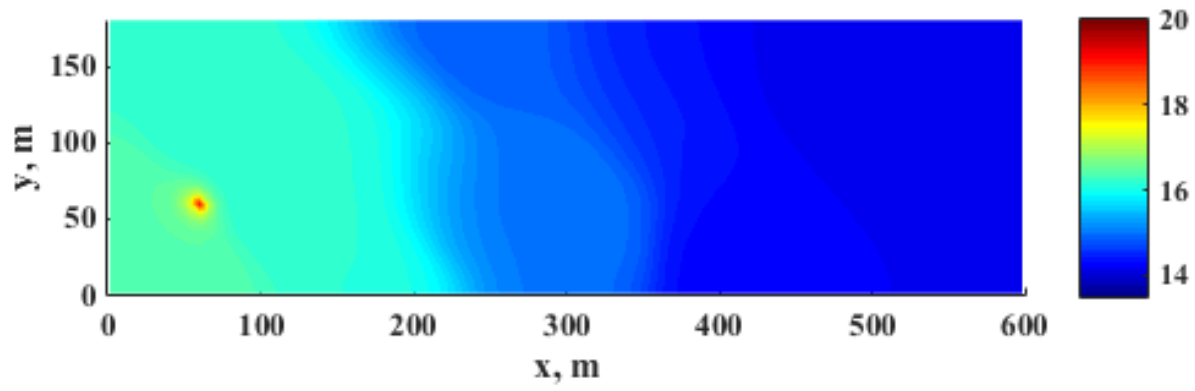

FVM reference solution

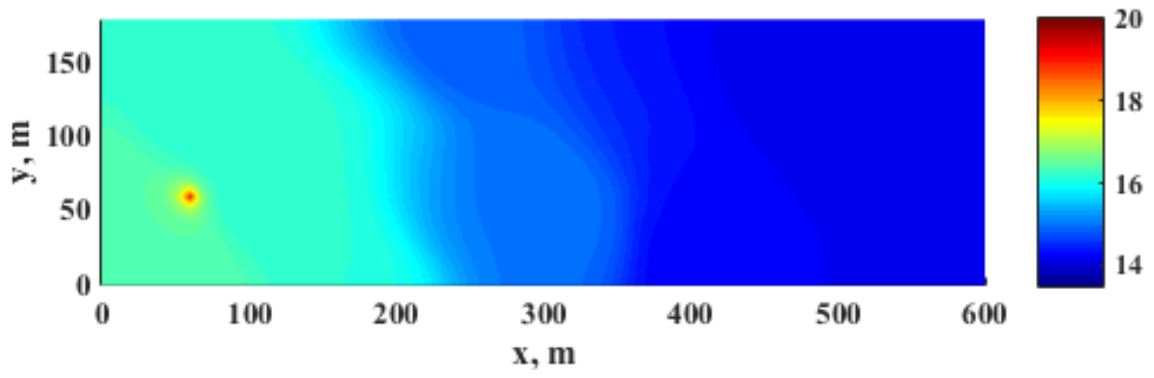

NCDMM, Cartesian, 10m, $w_1$

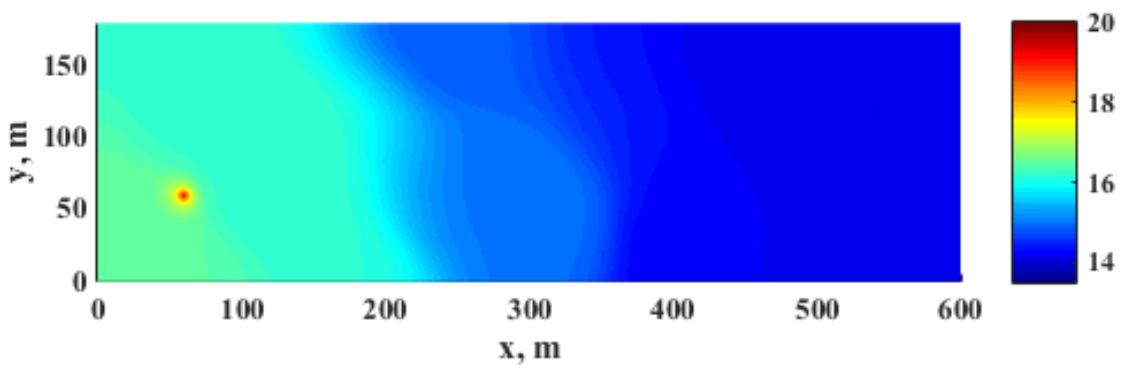

NCDMM, Cartesian, 10m, weighted $w_2$

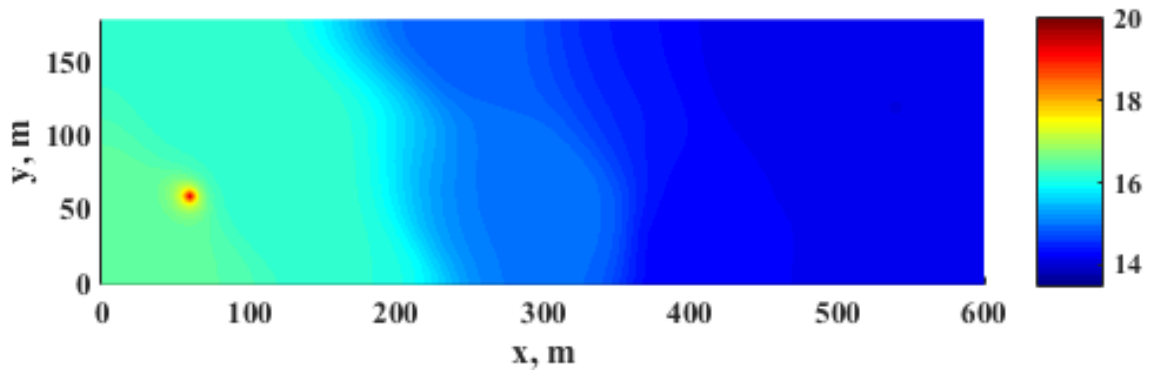

NCDMM, Cartesian, 15m, $w_1$

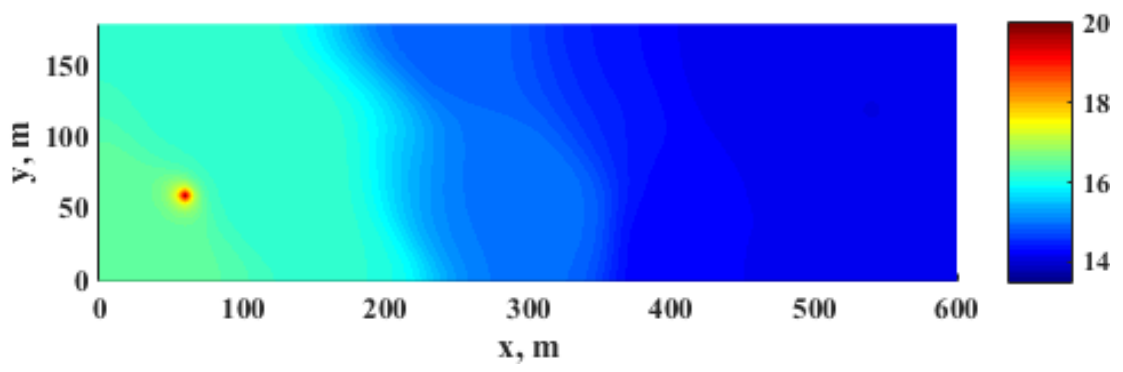

NCDMM, Cartesian, 15m, weighted $w_2$

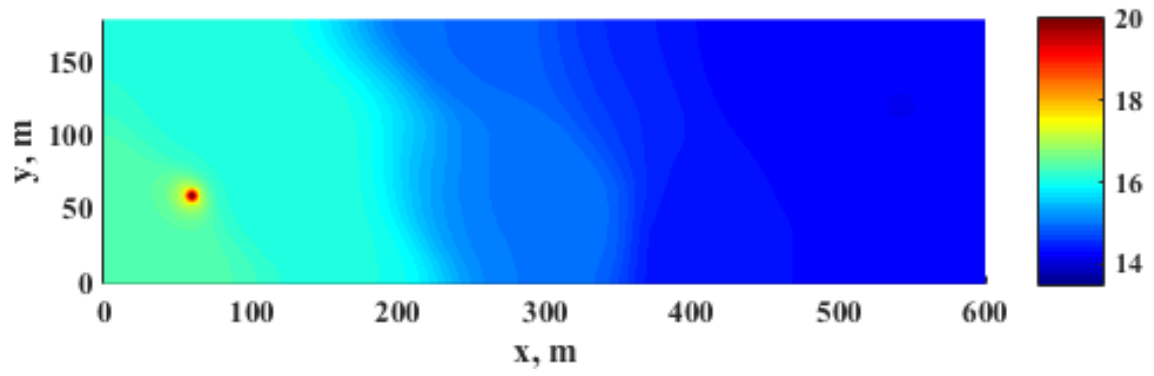

NCDMM, Cartesian, 20m, $w_1$

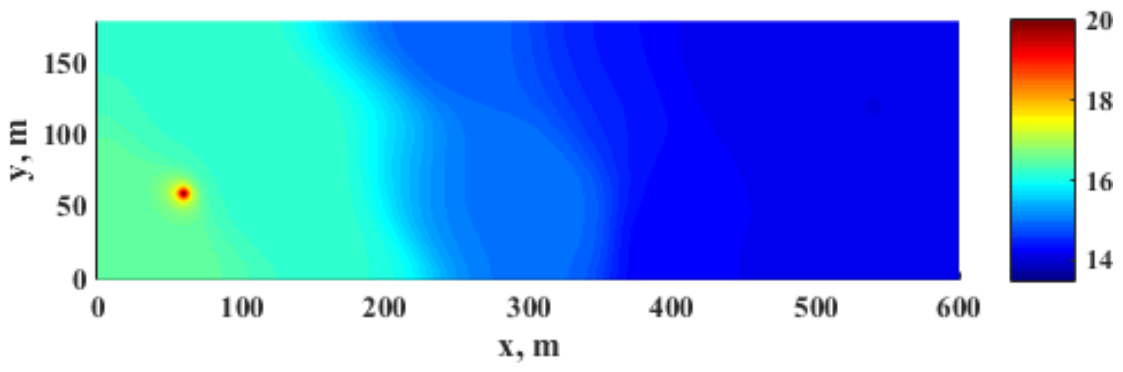

NCDMM, Cartesian, 20m, weighted $w_2$

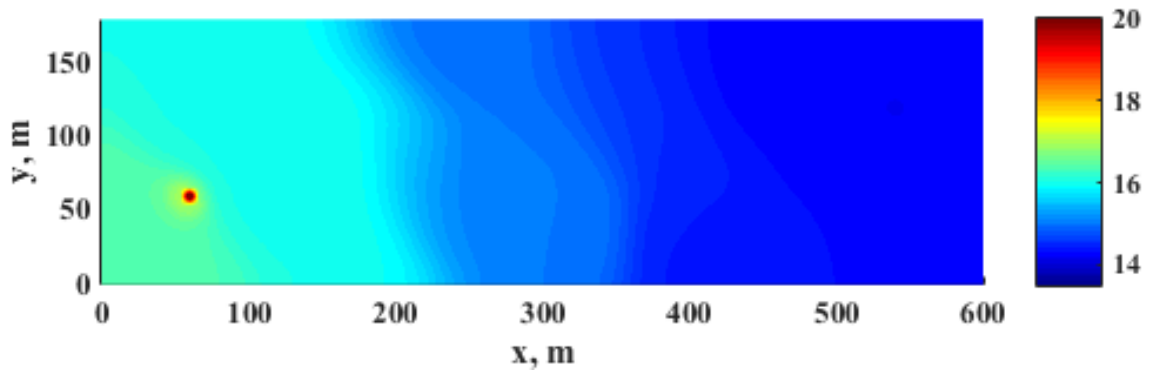

NCDMM, Cartesian, 25m, $w_1$

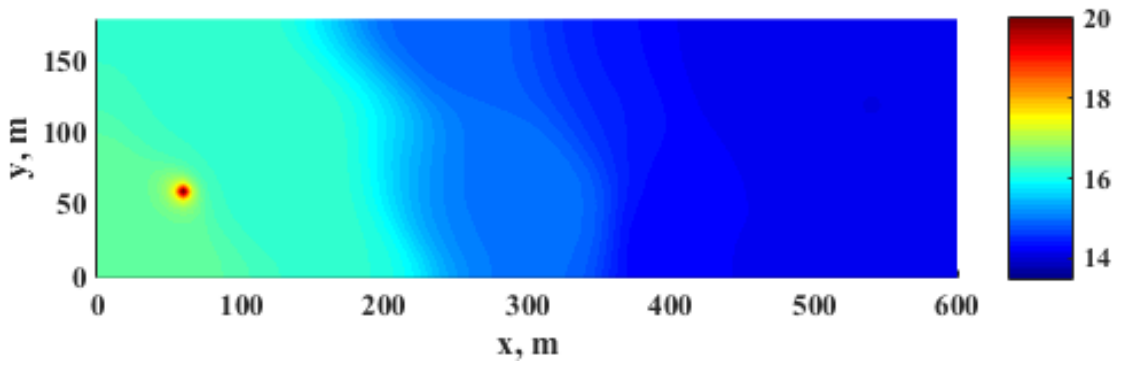

NCDMM, Cartesian, 25m, weighted $w_2$

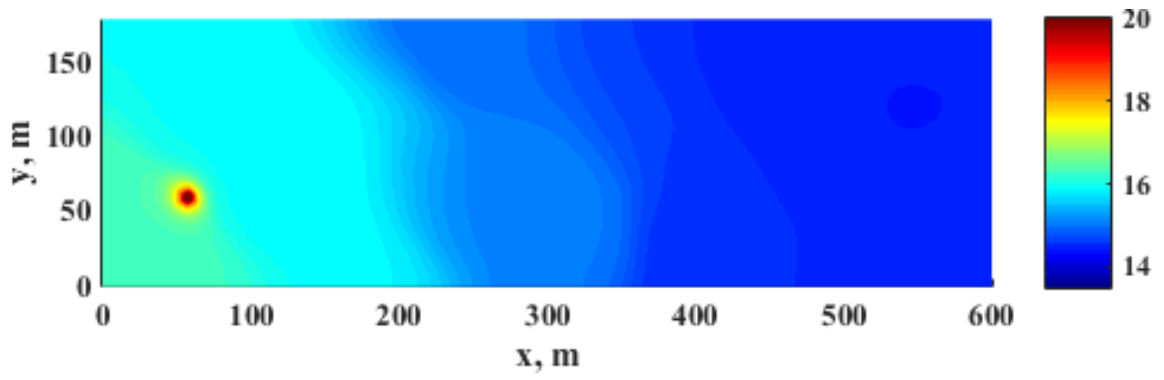

NCDMM, irregular, $\mathbf{NP}_1$, $w_1$

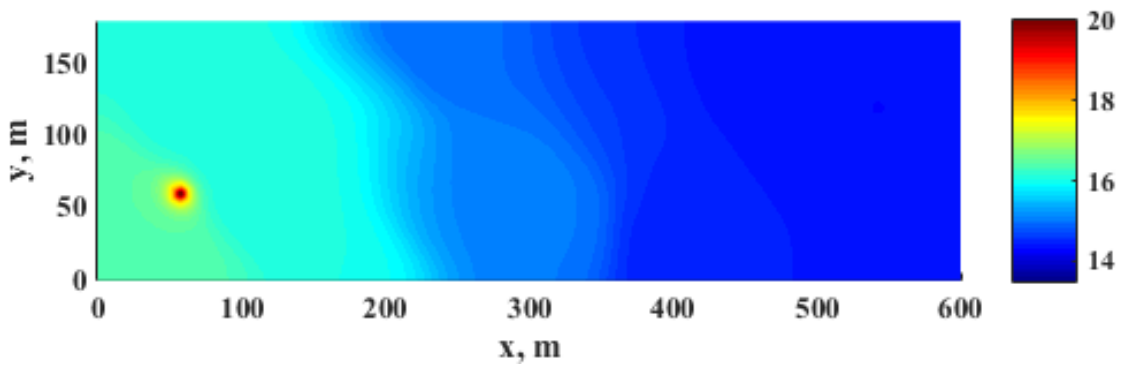

NCDMM, irregular, $\mathbf{NP}_1$, weighted $w_2$

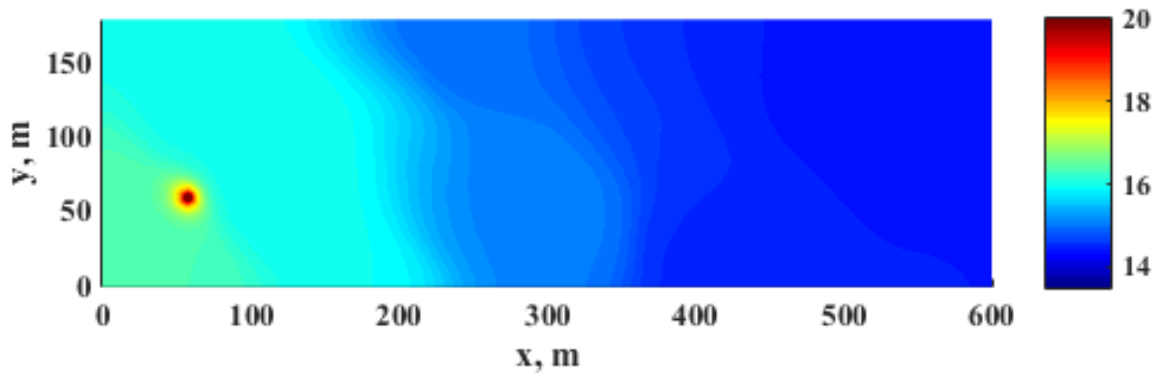

NCDMM, irregular, $\mathbf{NP}_2$, $w_1$

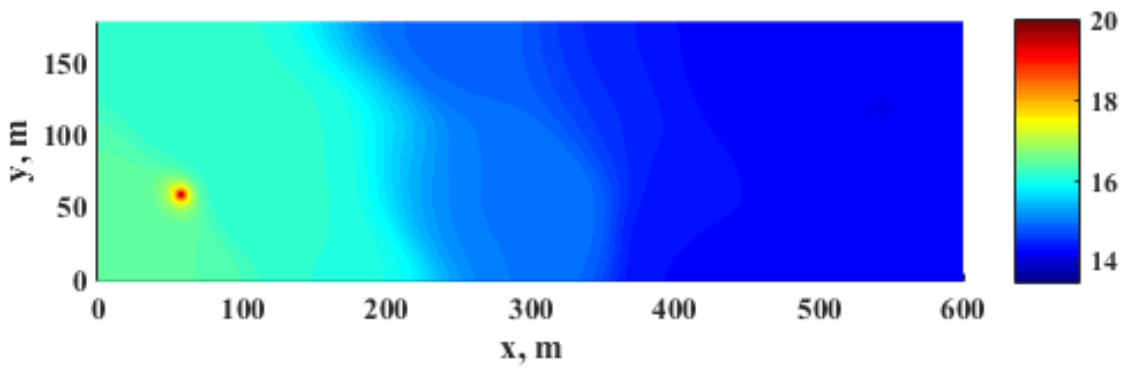

NCDMM, irregular, $\mathbf{NP}_2$, weighted $w_2$

Fig. 13 Calculated pressure profiles at 300 days in different cases

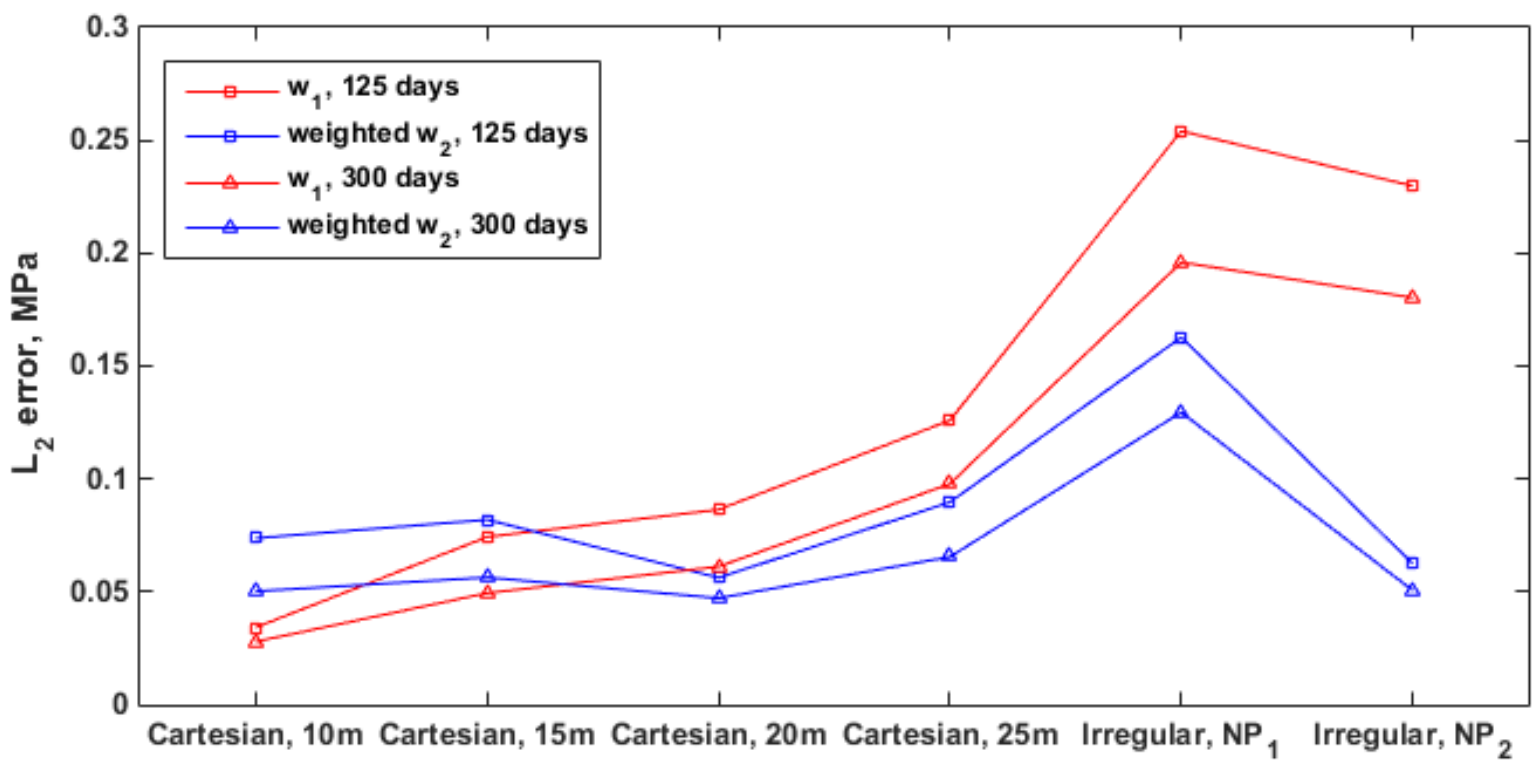

(a) pressure

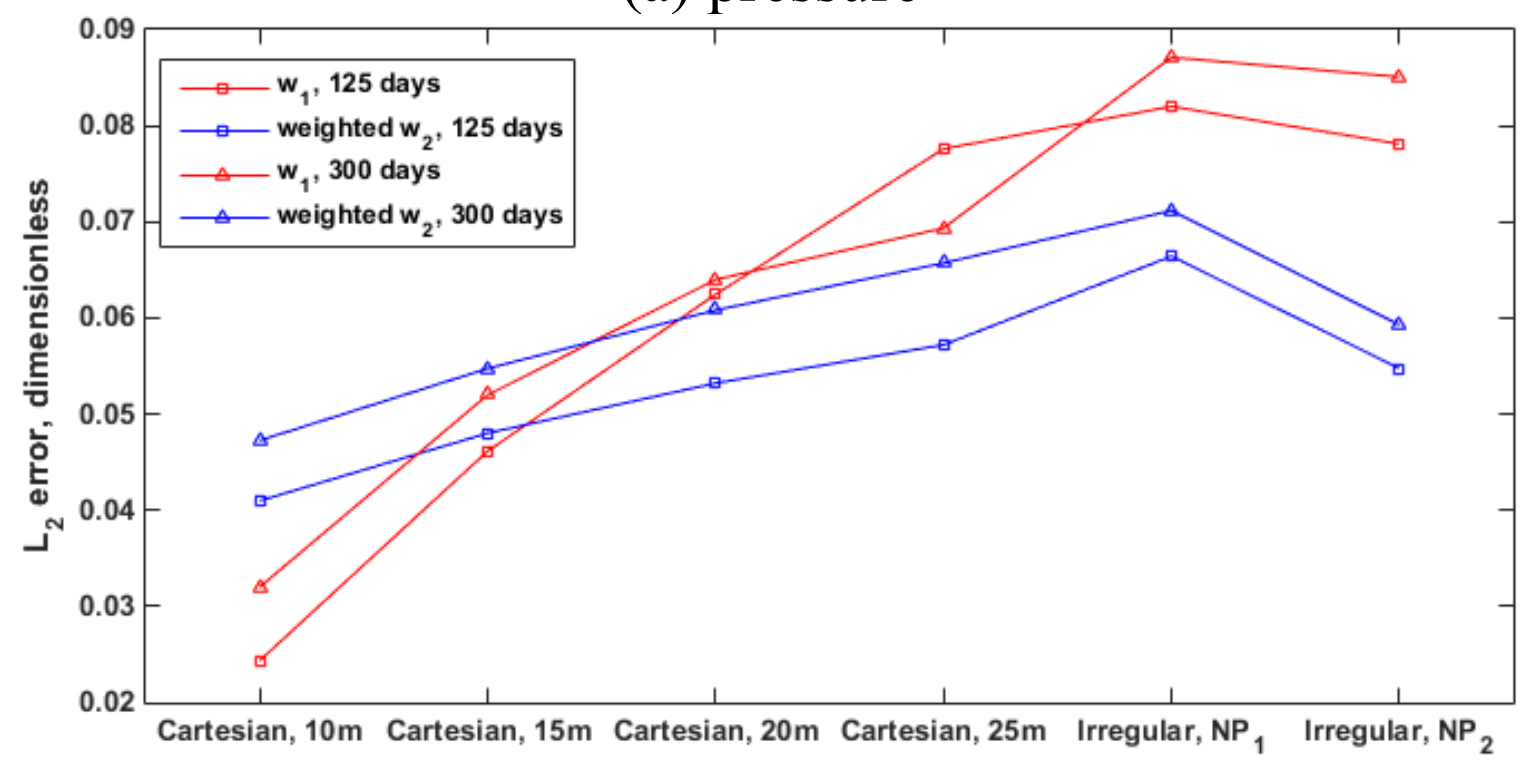

(b) oil saturation

Fig. 14 relative errors of calculated pressure and oil saturation profiles in different cases

Table 3 Parameters in Newton nonlinear solver

| Parameter | Value |
|---|---|
| Maximum time step size | 2 days |
| Minimum time step size | 0.001 day |
| Maximum Newton iterations | 50 |
| Error tolerance | $10^{-6}$ |
| $\eta_p$ | 5 MPa |
| $\eta_{S_w}$ | 0.05 |

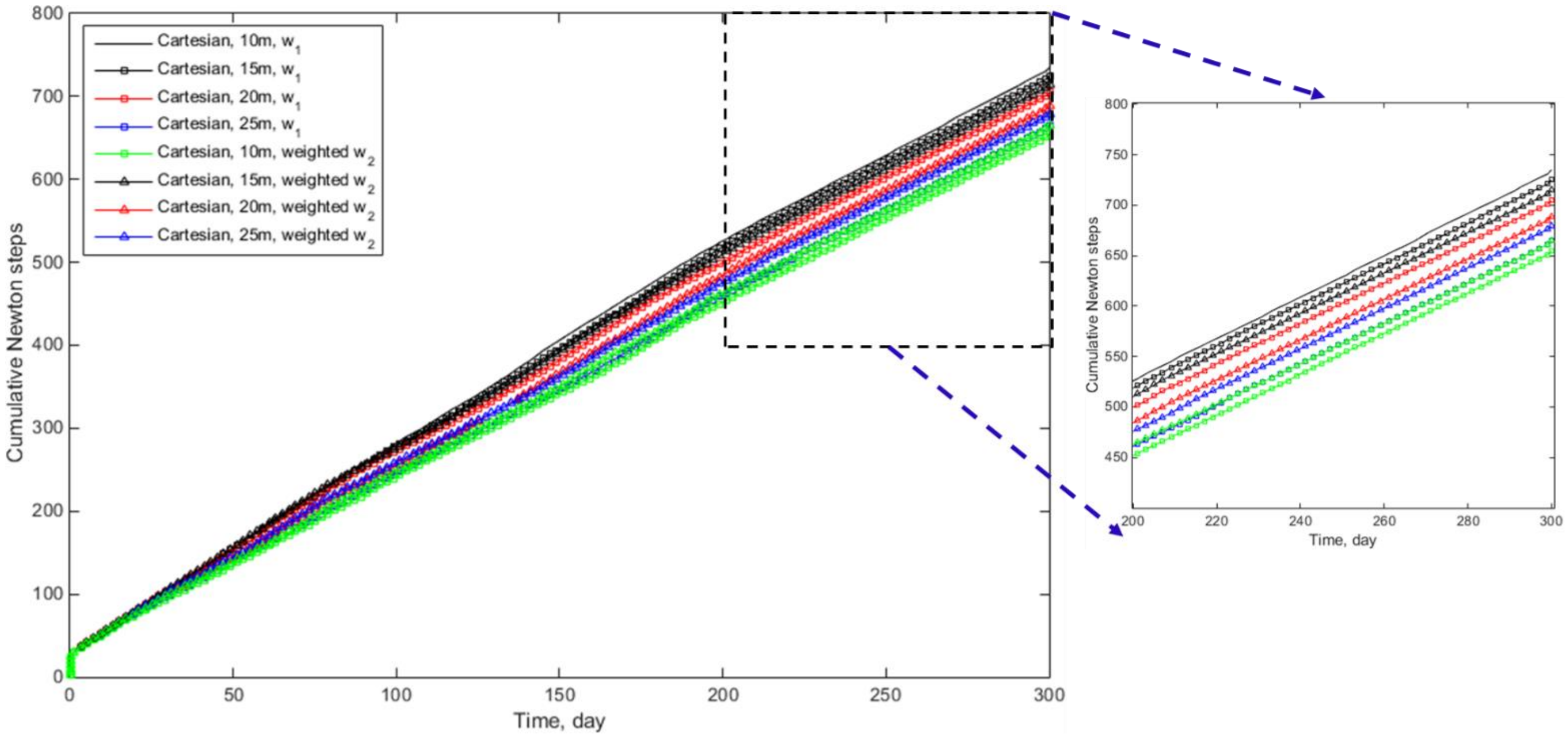

Fig. 15 Cumulative Newton iteration steps versus the production time

## 3.2 Two-phase porous flow in a polygonal domain

In this section, the rectangular reservoir model is changed to the polygonal domain in Fig. 16 (a). The heterogeneous permeability meets Eq. (67), and other physical parameters and the closed boundary condition are the same as those in Section 3.1.

$$k = 100 \times \exp\left[2\left(x/600\right)^2 + 2\left(y/180\right)^2\right] \tag{67}$$

Fig. 16 (b) and (c) show the two point clouds employed in this example, one of which is composed of mesh vertices of the triangulation of the reservoir domain. The other is called the pseudo-Cartesian point cloud, which is composed of Cartesian points inside the calculation domain and points on the boundary of the calculation domain. Its generation steps are:

**Step 1**: Firstly, generate nodes to discretize the boundary of the computational domain, and the set of the nodes is denoted as $\mathbf{I}_1$.

**Step 2**: regardless of the specific shape of the boundary, generate a Cartesian lattice for the larger rectangular area containing the calculation domain (the distance between adjacent nodes in the Cartesian lattice can be taken as the distance between neighboring nodes in $\mathbf{I}_1$);

**Step 3**: leave the nodes in the Cartesian lattice within the boundary (it can be realized by judging whether the points are in the polygon formed by the nodes in $\mathbf{I}_1$, such as the inpolygon function in Matlab), denote this part of the nodes as $\mathbf{I}_2$, and eliminate the nodes in $\mathbf{I}_2$ that are too close to the nodes in $\mathbf{I}_1$ (the distance threshold can generally be selected as half of the average spacing of nodes in the Cartesian lattice). This part of the eliminated nodes is denoted as $\mathbf{I}_3$, then the finally generated point cloud is $(\mathbf{I}_2 - \mathbf{I}_3) \cup \mathbf{I}_1$.

It can be seen that the quasi-Cartesian point cloud is suitable for complex geometries, and the only difficulty in generating this point cloud is to generate nodes on the boundary of the calculation domain in step 1, and the boundary of the calculation domain is one dimension lower than that of the calculation domain. Therefore, the generation difficulty of the quasi-Cartesian point cloud is less than that of directly generating the point cloud of the calculation domain with a complex boundary geometry.

For the first point cloud, this example uses the triangulation-based method given in Section 2.4 to determine the connectable point cloud, which is denoted as $\mathbf{NP}_1$. For the second point cloud, this example constructs the connectable point cloud by setting the influence domain of each node to have the same radius of 9m, which is denoted as $\mathbf{NP}_2$.

Fig. 17 compares the calculated profiles of node control volumes when adopting $w_1$, $w_2$, and weighted $w_2$. It can be seen that the node distribution density of these two point clouds in the middle of the calculation domain is relatively even, especially the quasi-Cartesian point cloud, in which except that the node density near the boundary is different, the node distribution density inside the calculation domain is the same everywhere. Therefore, a reasonable node control volume distribution should be even inside the calculation domain. For example, for the quasi-Cartesian point cloud, the control volume of the inner points far away from the domain boundary should be 25m$^2$. As seen from these different node control volume distribution maps, whether $\mathbf{NP}_1$ or $\mathbf{NP}_2$, when $w_1$ is used, the calculated node control volumes in the middle of the computational domain are significantly higher than the actual node control volume, and the overall node control volume distribution generally shows the characteristics of rapid reduction from the middle of the computational domain to the boundary of the computational domain. With $w_1$ to $w_2$, and then to weighted $w_2$, the calculated node control volumes within the calculation domain are more and more evenly distributed. For example, for the quasi-Cartesian point cloud, the maximum node control volume calculated by $w_1$ reaches nearly 35m$^2$, while the node control volumes calculated by weighted $w_2$ are generally about 25m$^2$. Moreover, it can be seen that when $w_1$ is just replaced by $w_2$, the calculation accuracy of the node control volume distribution has been greatly improved, which also fully shows that the use of $w_2$ is indeed helpful to improve the calculation accuracy of node control volumes.

As shown in Fig. 16 (b), construct a rectangular model that can include the polygonal calculation domain, discretize the rectangular domain with the simple Cartesian mesh, make the mesh permeability outside the polygonal domain 0 mD to play the role of the closed boundary condition, and then use the FVM-based solution of the rectangular model as the reference solution.

Fig. 18, 19, 20, and 21 compare the calculated water saturation and pressure profiles in different cases, and Fig. 22 compares the $L_2$ errors in different cases. It can be seen that whether $\mathbf{NP}_1$ or $\mathbf{NP}_2$ is used, compared with $w_1$, the pressure and oil saturation distribution calculated by weighted $w_2$ are intuitively closer to the reference solution, and the corresponding $L_2$ error is significantly lower. In this example, when using weighted $w_2$, the $L_2$ error of $\mathbf{NP}_2$ is generally slightly lower than that of $\mathbf{NP}_1$, which is the same as the analysis in Section 3.1, and the $L_2$ error of pressure is less than 0.08MPa, and the $L_2$ error of oil saturation is less than 0.035, which shows that the proposed NCDMM can achieve satisfactory calculation accuracy when using different connectable point clouds to discretize the calculation domain with polygonal boundary geometry, and also reflects that the type of connectable point cloud directly affects the calculation accuracy, which is similar to the influence of mesh quality on the calculation performance of mesh-based numerical methods.

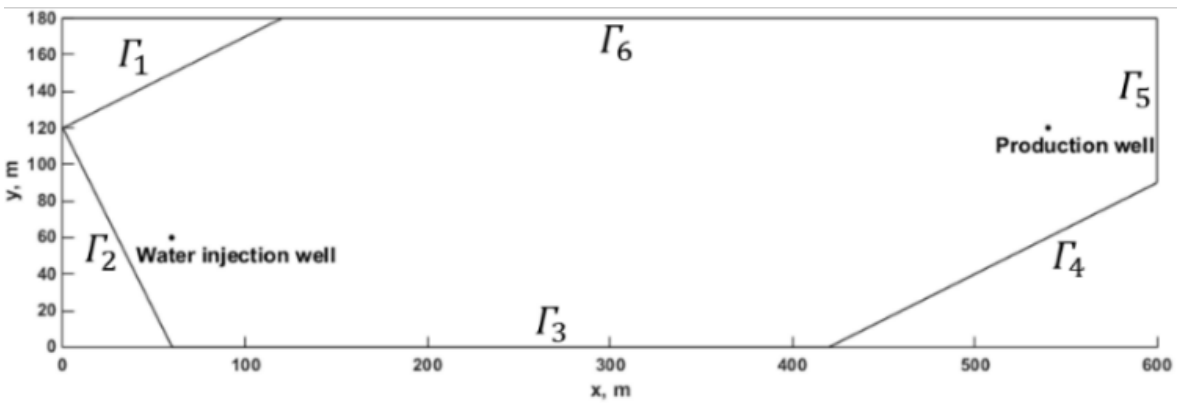

(a) reservoir calculation domain

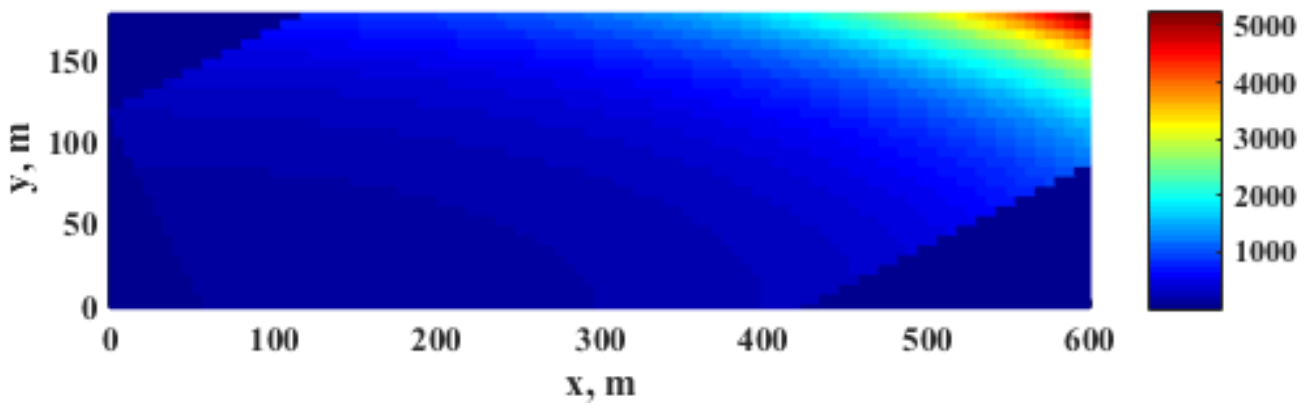

(b) permeability distribution

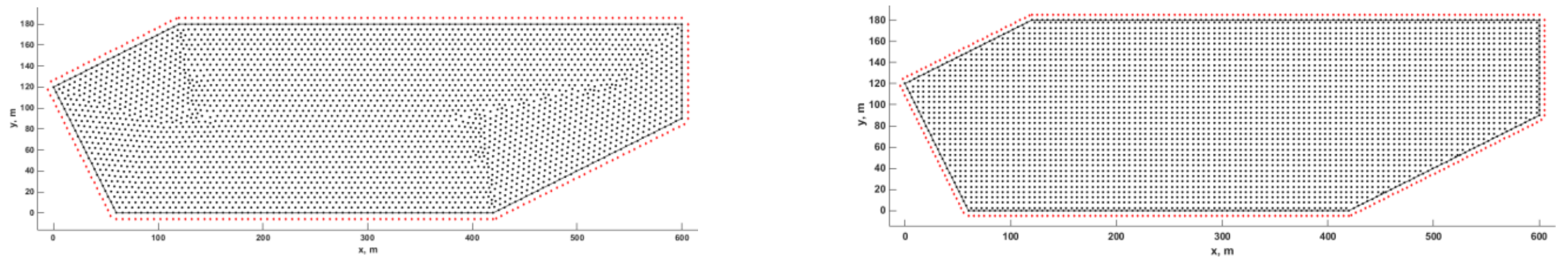

(c) $\mathbf{NP}_1$ (d) $\mathbf{NP}_2$

Fig. 16 Sketches of the reservoir domain and point clouds

$w_1$

$w_2$

Weighted $w_2$

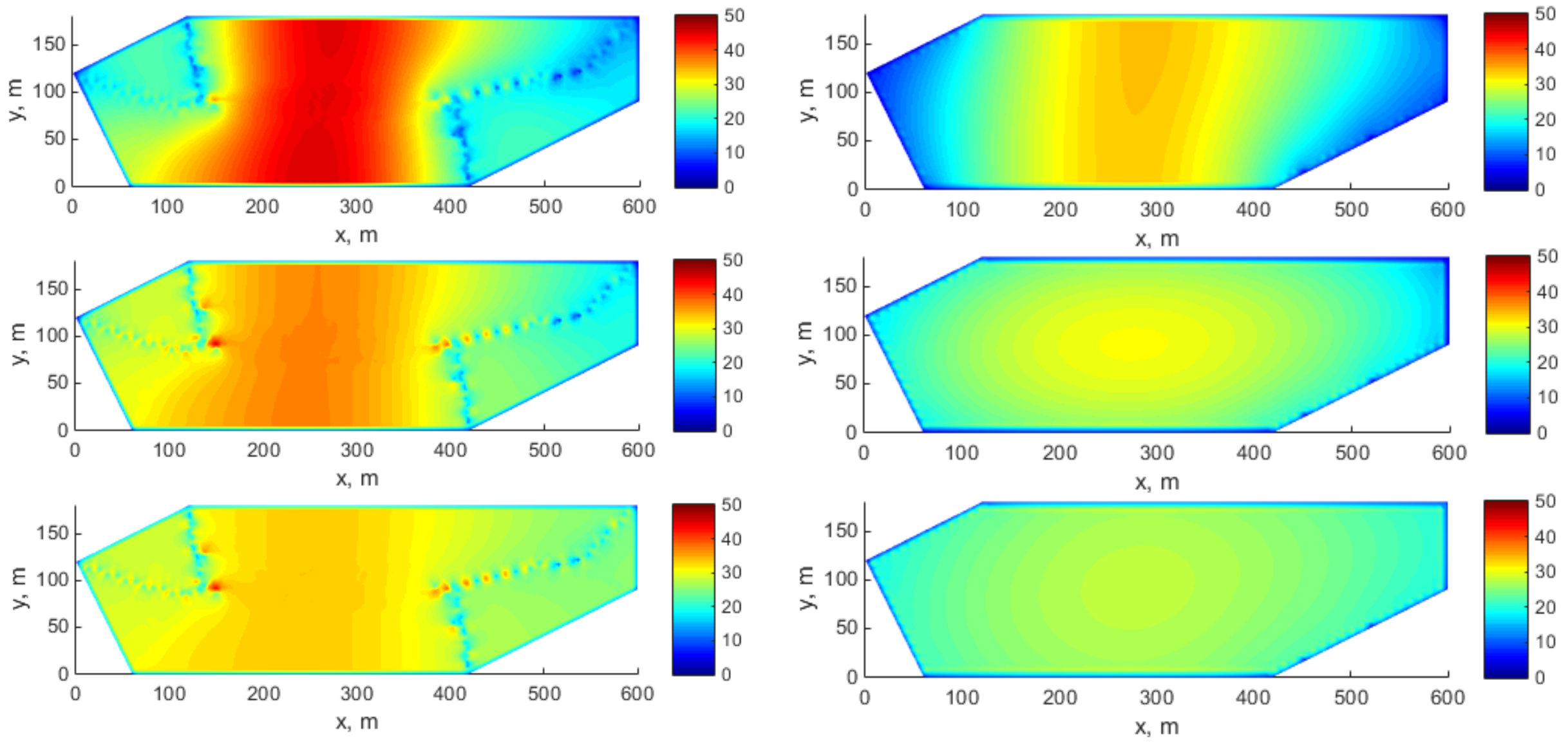

(a) $\mathbf{NP}_1$ (b) $\mathbf{NP}_2$

Fig. 17 Calculated profiles of node control volumes by adopting $w_1$, $w_2$, and weighted $w_2$

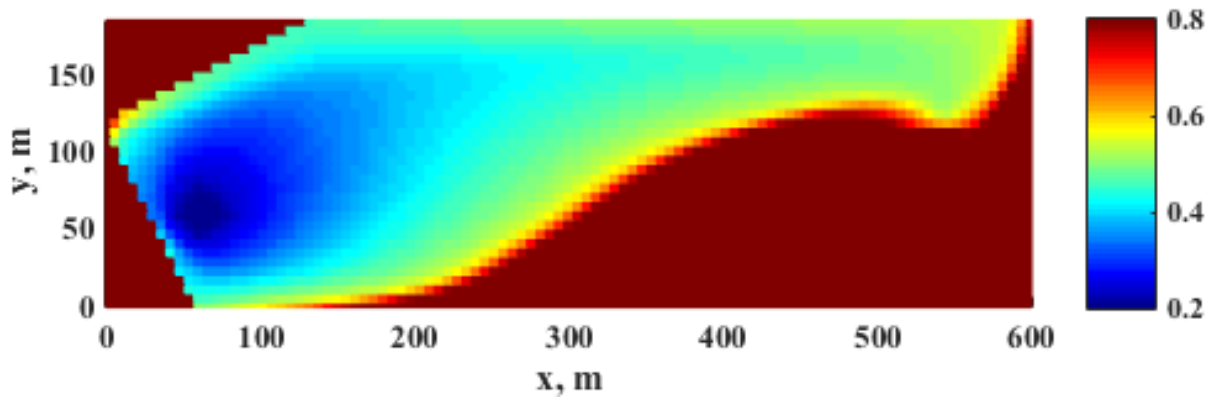

(a) FVM reference solution

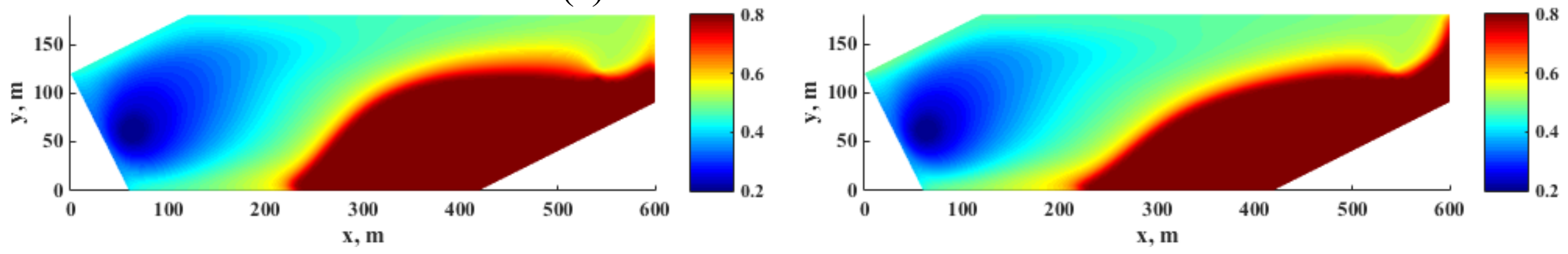

(b) $\mathbf{NP}_1$, $w_1$ (c) $\mathbf{NP}_1$, weighted $w_2$

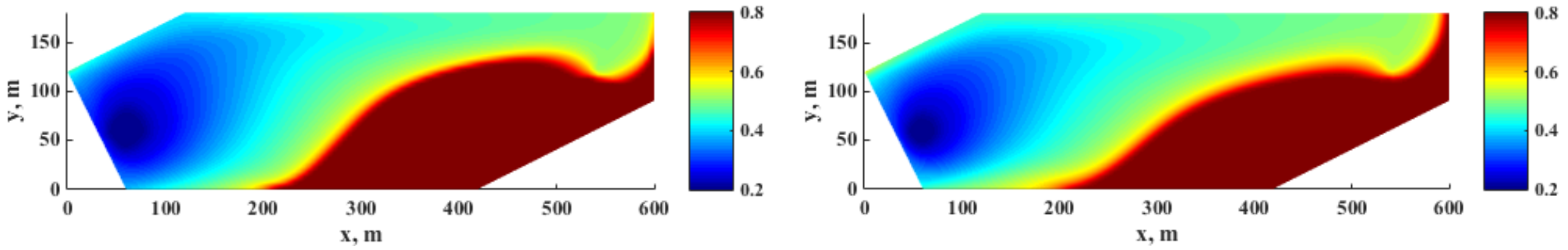

(d) $\mathbf{NP}_2$, $w_1$ (e) $\mathbf{NP}_2$, weighted $w_2$

Fig. 18 Comparison of calculated oil saturation profiles at 250 days

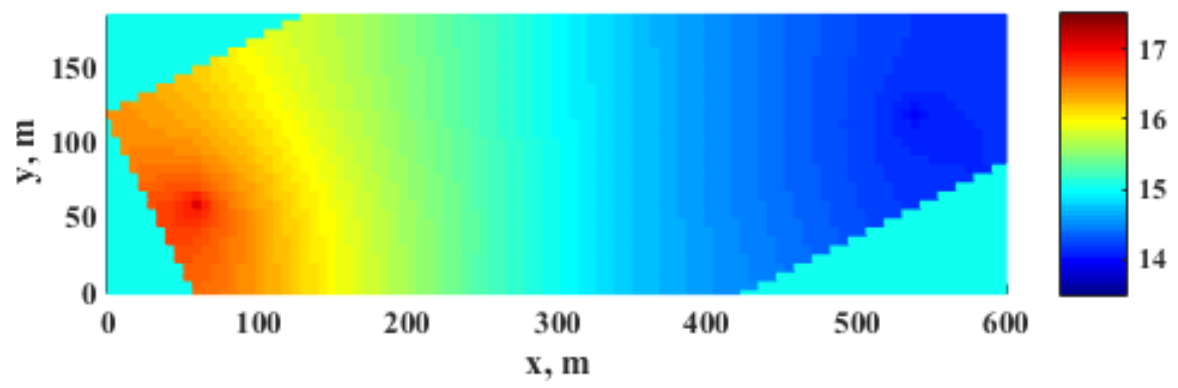

(a) FVM reference solution

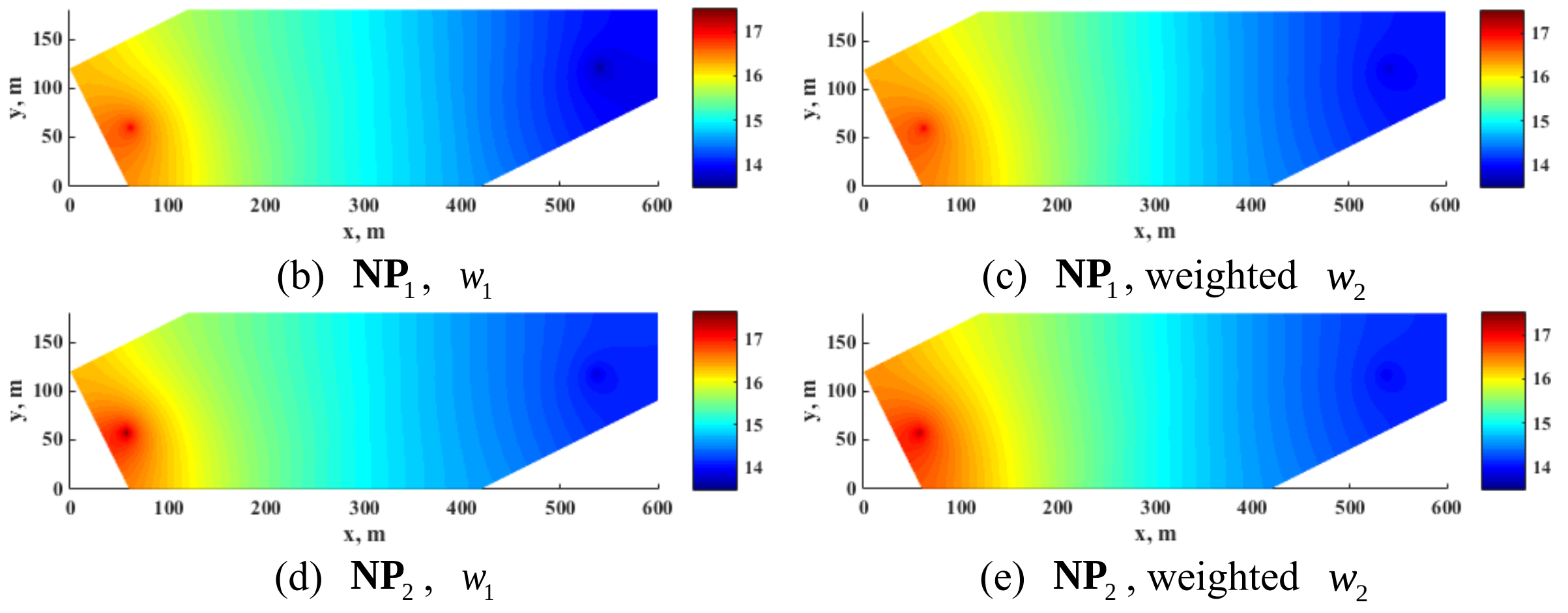

(b) $\mathbf{NP}_1$, $w_1$ (c) $\mathbf{NP}_1$, weighted $w_2$

(d) $\mathbf{NP}_2$, $w_1$ (e) $\mathbf{NP}_2$, weighted $w_2$

Fig. 19 Comparison of calculated pressure profiles at 250 days

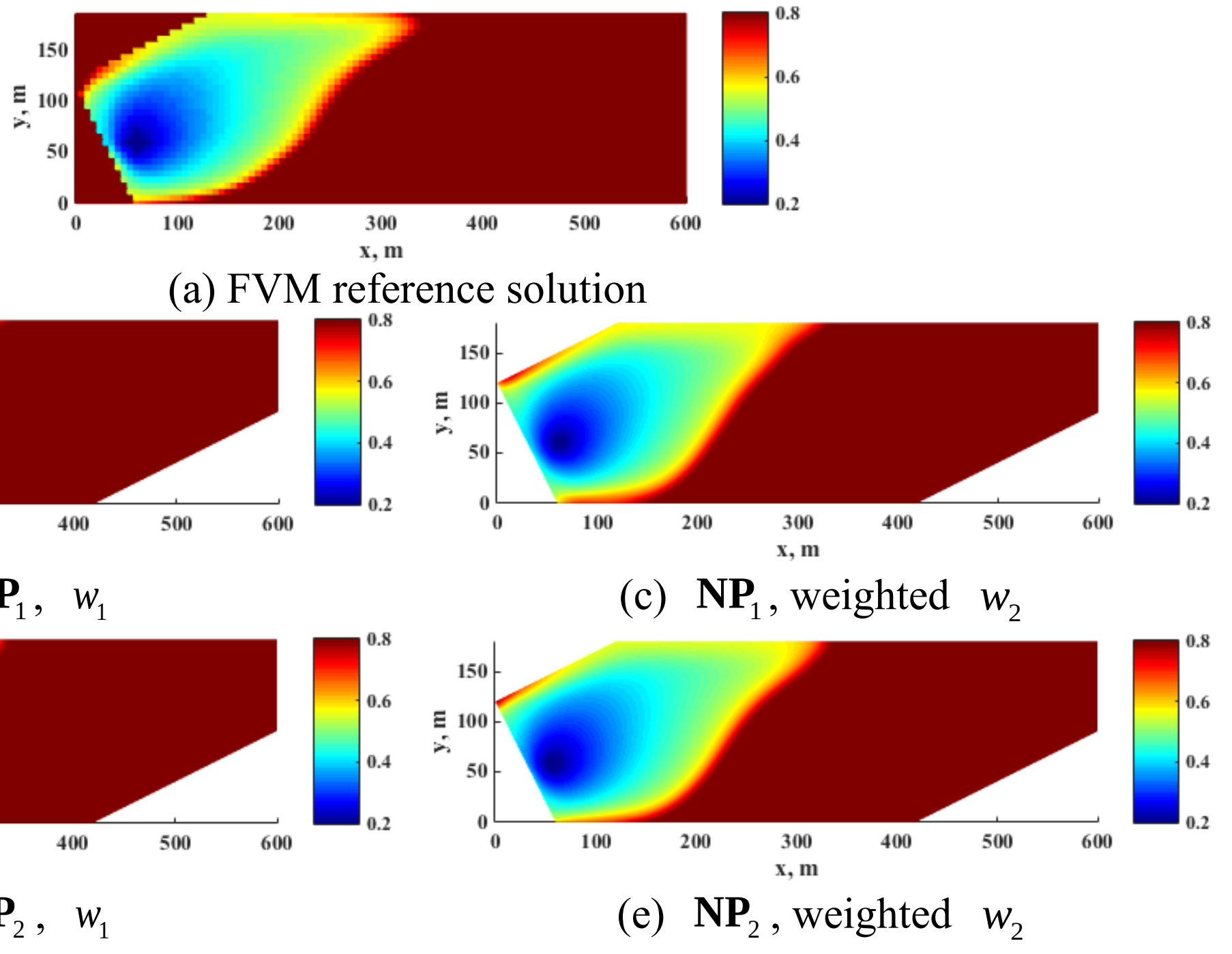

(a) FVM reference solution

(b) $\mathbf{NP}_1$, $w_1$ (c) $\mathbf{NP}_1$, weighted $w_2$

(d) $\mathbf{NP}_2$, $w_1$ (e) $\mathbf{NP}_2$, weighted $w_2$

Fig. 20 Comparison of calculated oil saturation profiles at 125 days

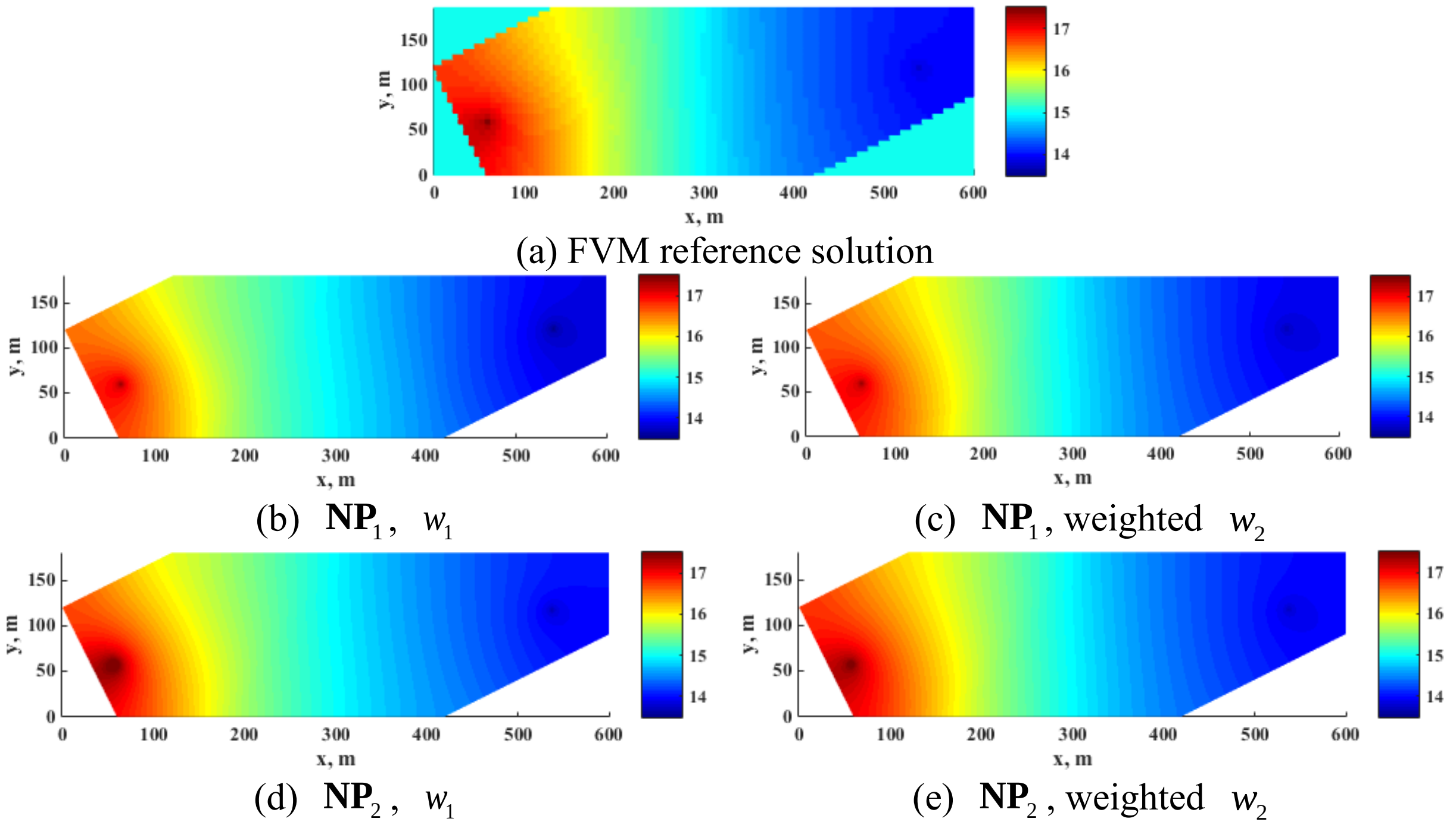

(a) FVM reference solution

(b) $\mathbf{NP}_1$, $w_1$ (c) $\mathbf{NP}_1$, weighted $w_2$

(d) $\mathbf{NP}_2$, $w_1$ (e) $\mathbf{NP}_2$, weighted $w_2$

Fig. 21 Comparison of calculated pressure profiles at 125 days

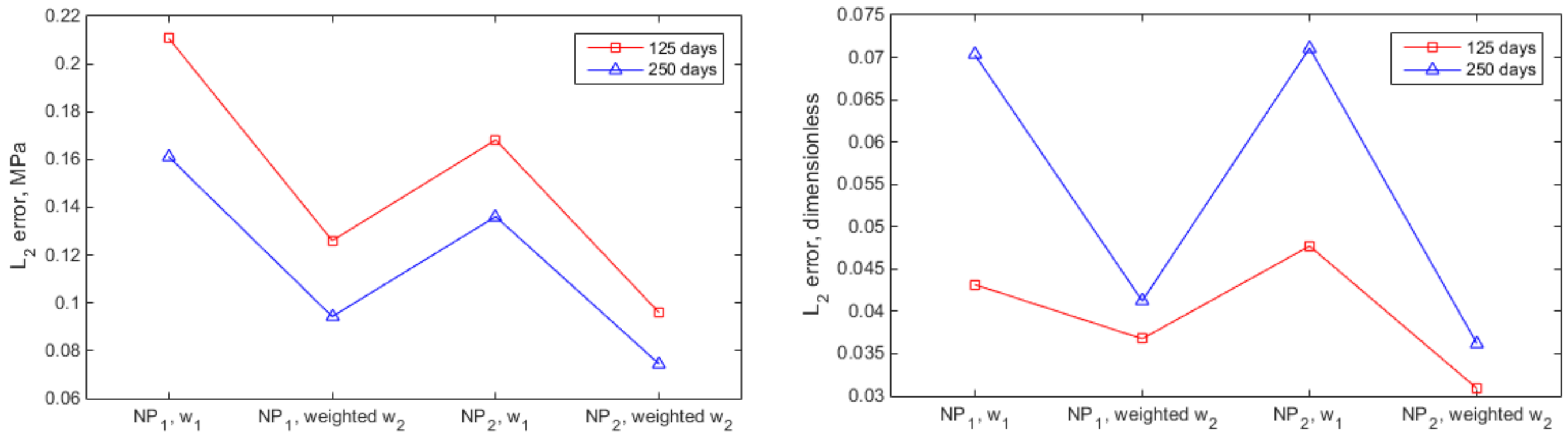

(a) pressure (b) oil saturation

Fig. 22 $L_2$ errors of calculated pressure and oil saturation profiles

3.3 Two-phase porous flow with a Dirichlet boundary condition

In this section, the reservoir model, node collocations, and relevant physical properties in Section 3.2 are still used, while the permeability is set as 100mD. There are two production wells in coordinates (200, 72), (420, 103) with a constant liquid production rate 60m$^3$/d. The initial and boundary conditions are listed in Eq. (68) and Eq. (69) respectively, in which, the boundaries are at the Dirichlet boundary condition, indicating strong edge water. The initial water saturation is the irreducible water saturation.

$$p\big|_{t=0} = 10MPa\,,\quad S_w\big|_{t=0} = 0.2 \tag{68}$$

$$p\big|_{\Gamma} = 10MPa\,,\quad S_w\big|_{\Gamma} = 0.8\,,\quad \Gamma = \Gamma_1 \cup \Gamma_2 \cup \Gamma_3 \cup \Gamma_4 \cup \Gamma_5 \cup \Gamma_6 \tag{69}$$

Fig. 23 shows the point clouds with different average node spacings (3m, 6m, 10m, and 15m) composed of mesh vertices of triangulations. The triangulation-based method of determining the connectable point cloud is adopted. Fig. 24 compares the calculated profiles of node control volumes when adopting $w_1$, $w_2$, and weighted $w_2$. As seen from which, no matter what the average spacing of nodes is, with $w_1$ to $w_2$, and then to weighted $w_2$, except that the node control volume in the area with high node density remains reasonably small, the node control volume distribution in the area with even node density is becoming more and more even (i.e., more accurate). The FVM-based fine-mesh solution in a larger computational domain is used as the reference solution, in which, to impose the Dirichlet boundary condition in Eq. (68) when FVM is used, volume, pressure, and water saturation of the cells outside the real calculation domain are set as 10$^{10}$m$^3$, 10MPa, and 0.8 respectively, and the transmissibility between the cell outside the calculation domain and the cell inside the calculation domain is set large.

Fig. 25 and Fig. 26 compare the pressure and water saturation profiles at 100 days and 300 days calculated by NCDMM with these point clouds in Fig. 23. Fig. 27 shows the relative errors versus the average node spacings. As seen from the comparisons of NCDMM results and FVM reference solutions, when the average node spacing becomes smaller, the NCNMM results gradually converge to the reference solution, and the corresponding computational error gradually decreases. Therefore, although the FVM can also deal with the Dirichlet boundary conditions, it is not straightforward enough and is more complicated to deal with in the case of complex geometry. In contrast, NCDMM can handle Dirichlet boundary conditions very directly (including, of course, the Robin-type boundary conditions that are more difficult to be handled directly by FVM), and the comparison of computational results shows that NCDMM achieves good computational performance, demonstrating one of the computational advantages of NCDMM over FVM.

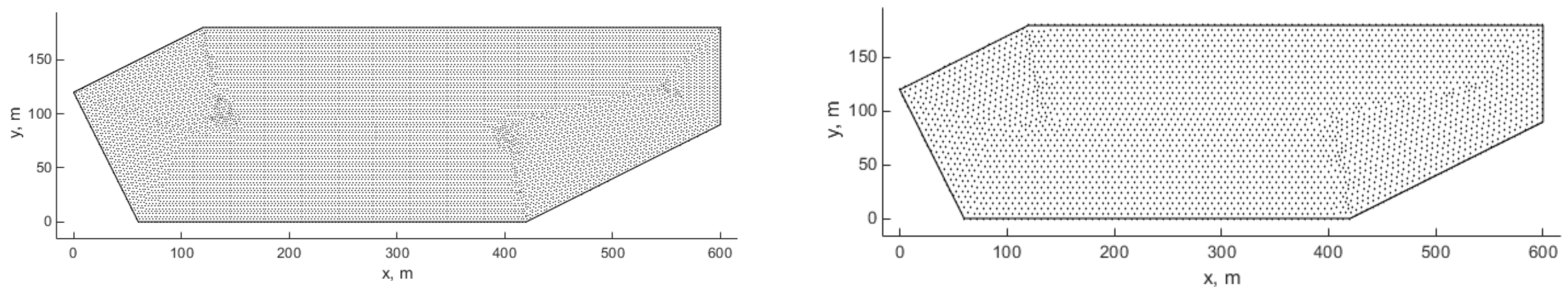

(a) the average node spacing is 3m (b) the average node spacing is 6m

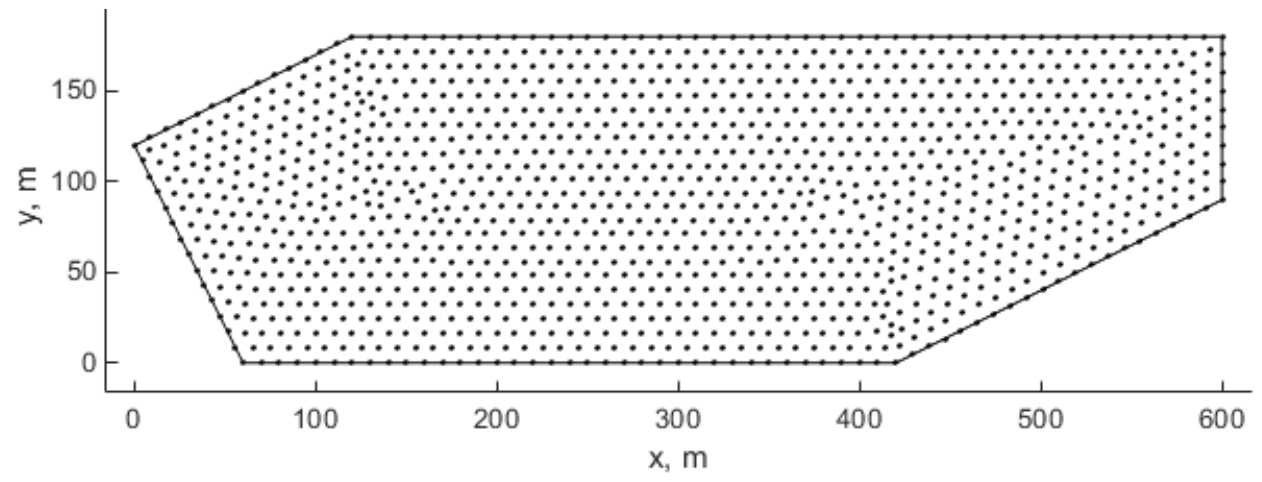

(c) the average node spacing is 10m

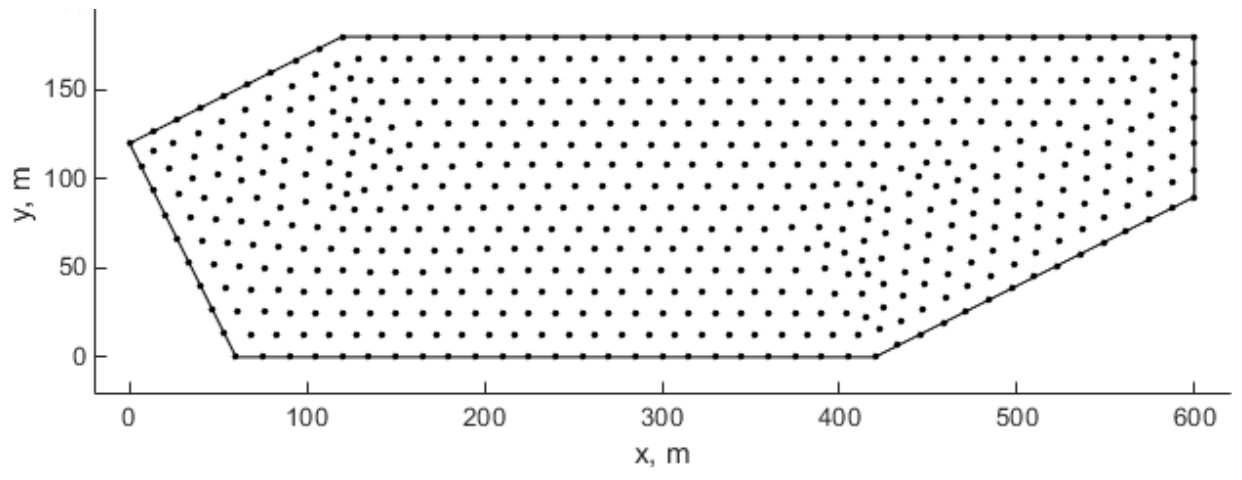

(d) the average node spacing is 15m

Fig. 23 The point clouds used in this example

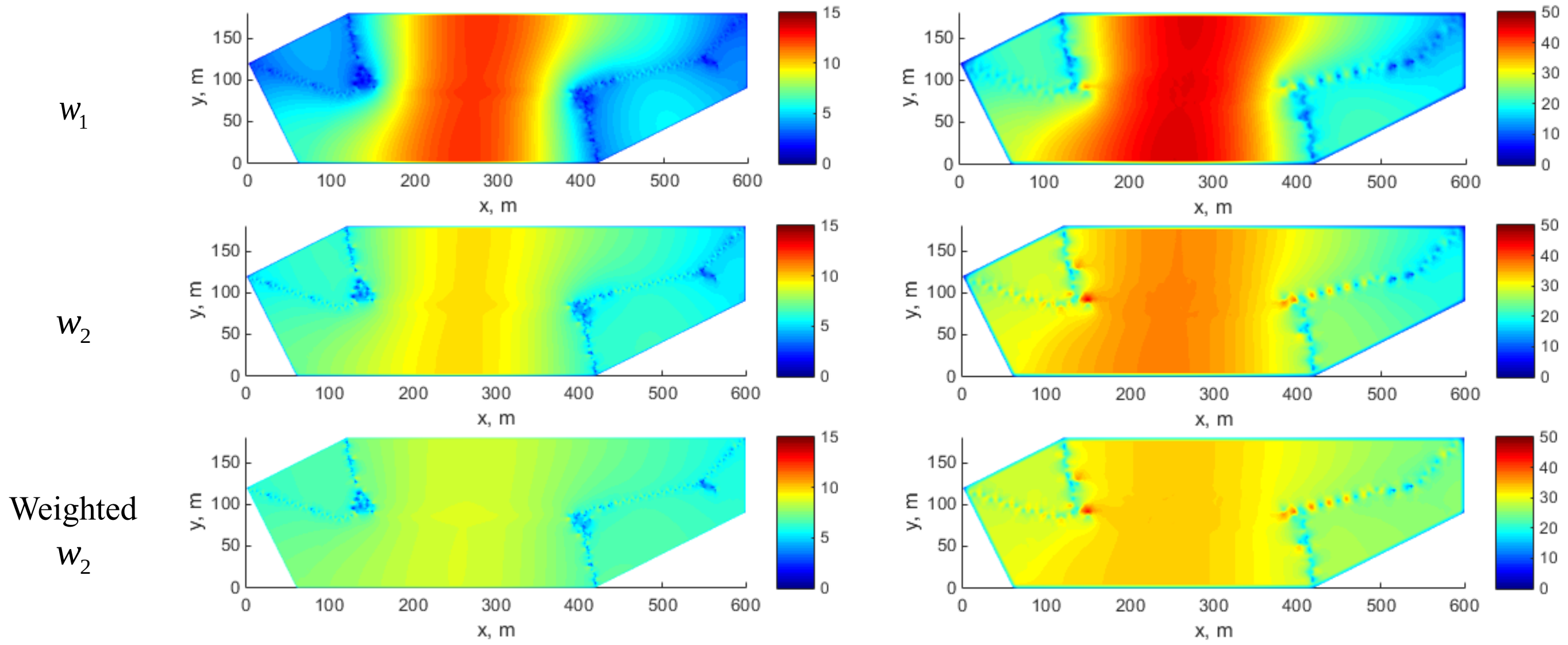

(a) 3m

(b) 6m

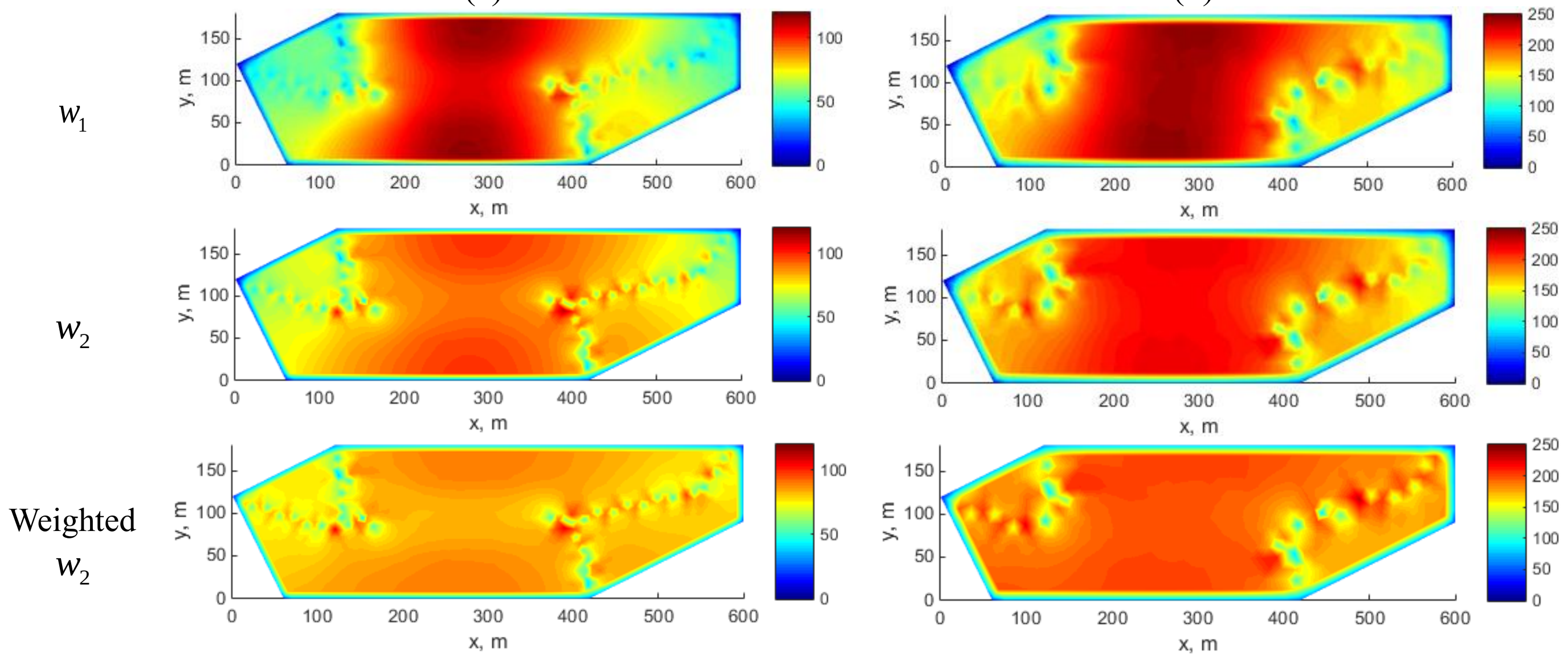

(c) 10m

(d) 15m

Fig. 24 Comparisons of calculated node control volume profiles in different cases

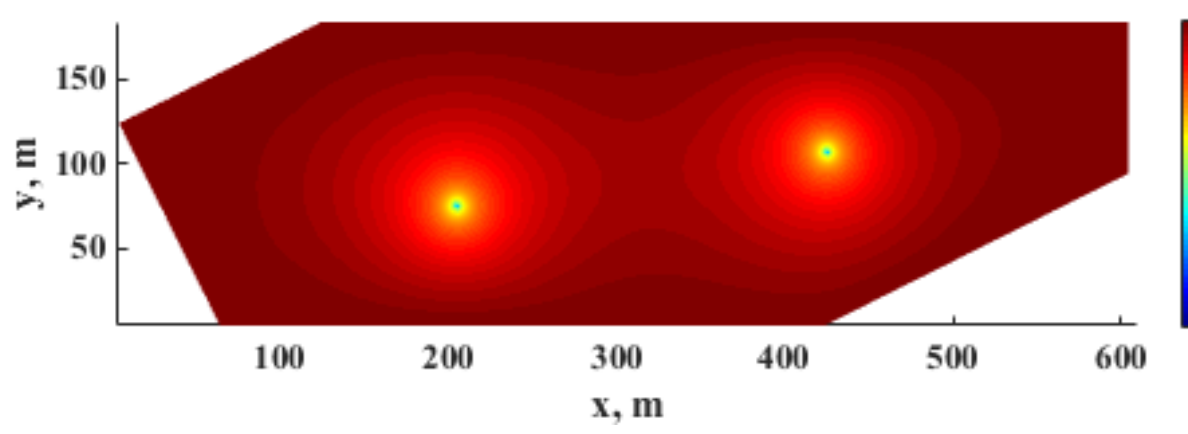

Reference solution, pressure, 300 days

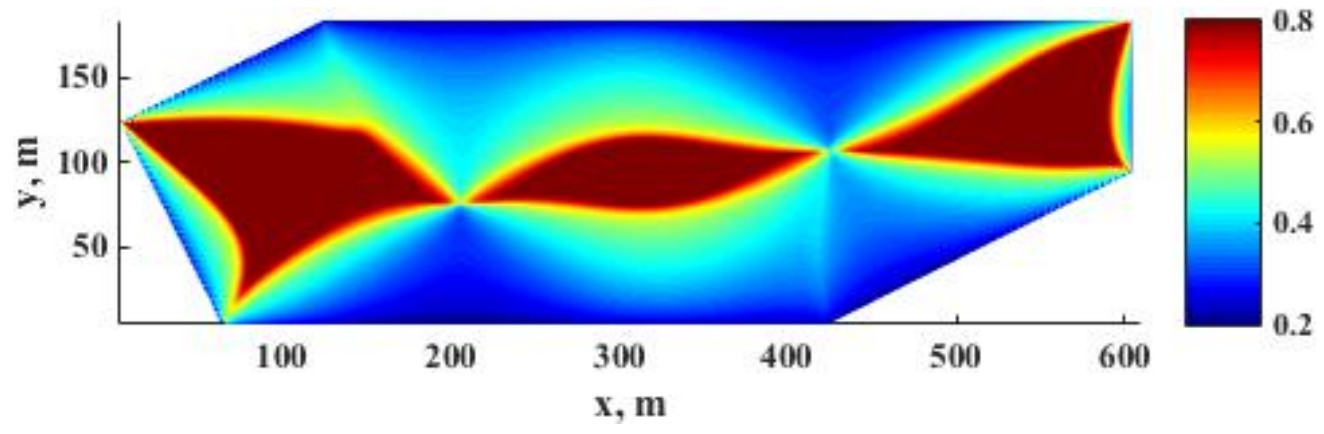

Reference solution, oil saturation, 300 days

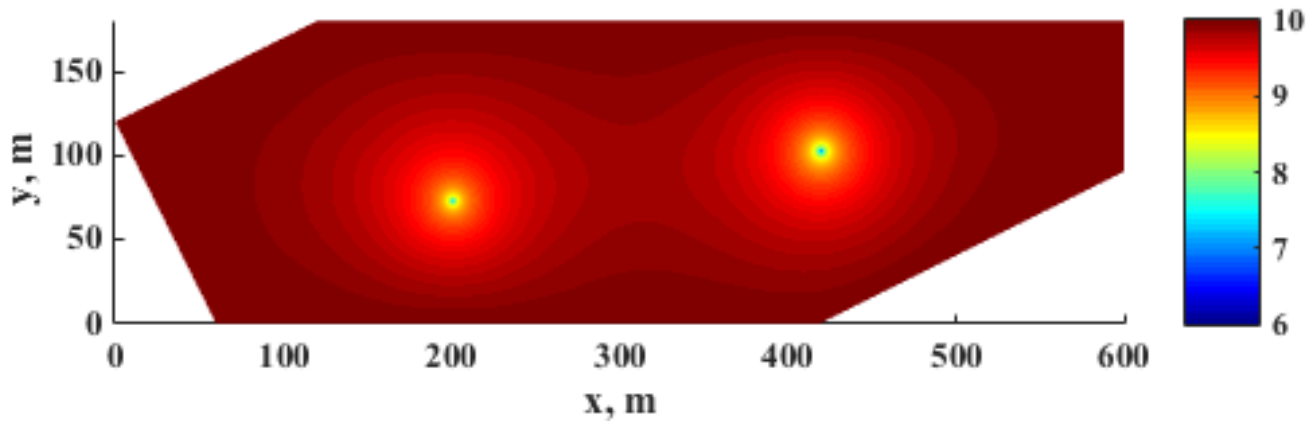

NCDMM solution, pressure, 300 days, 3m

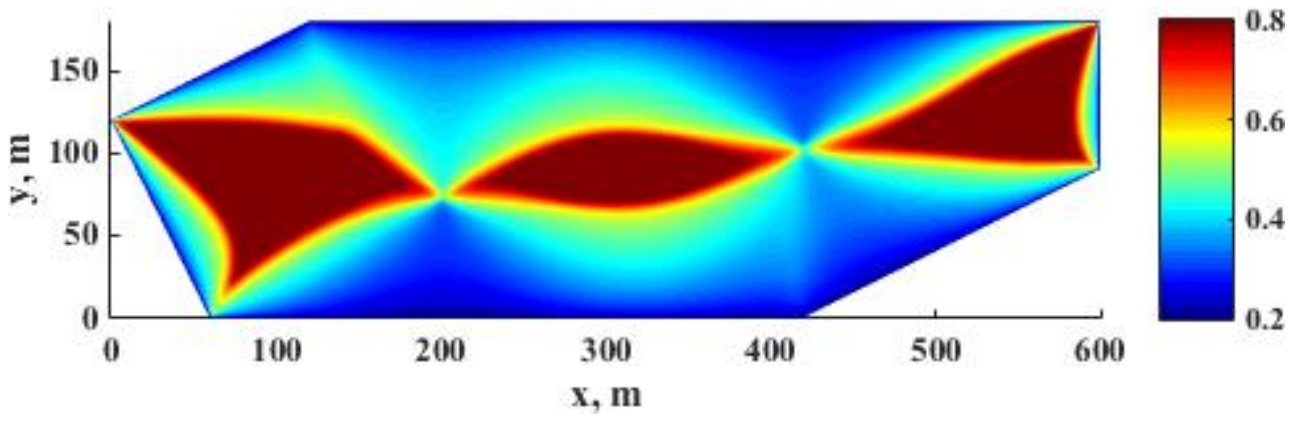

NCDMM solution, oil saturation, 300 days, 3m

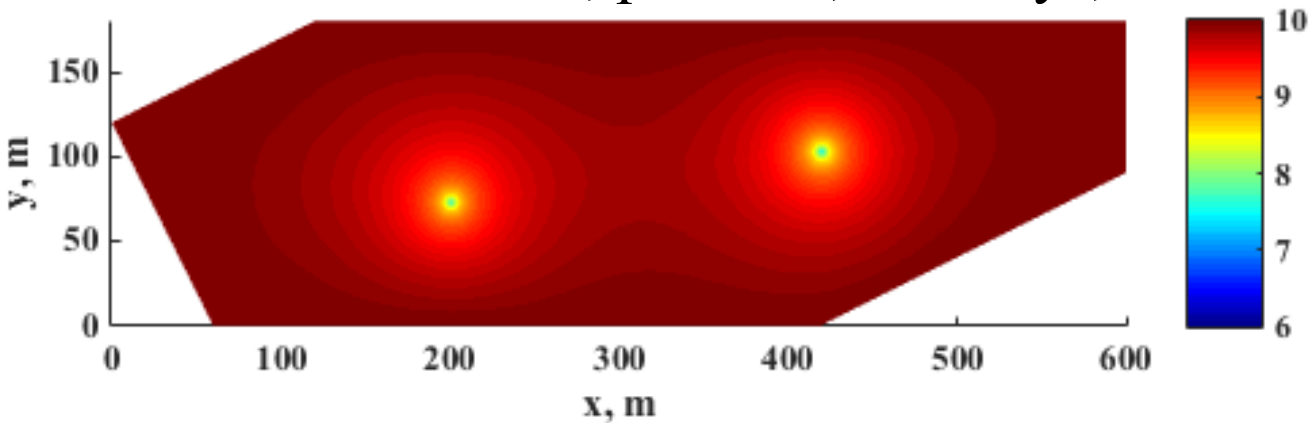

NCDMM solution, pressure, 300 days, 6m

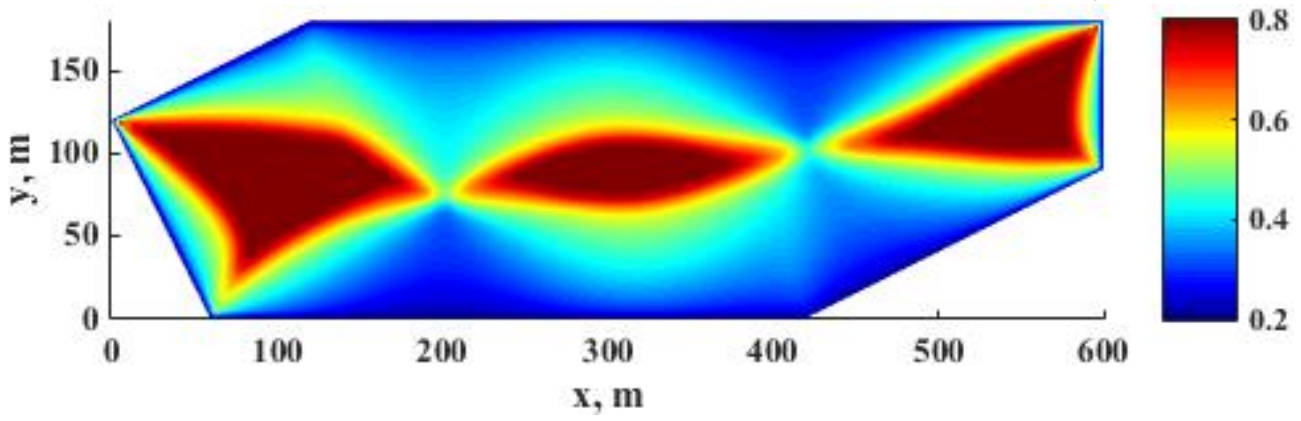

NCDMM solution, oil saturation, 300 days, 6m

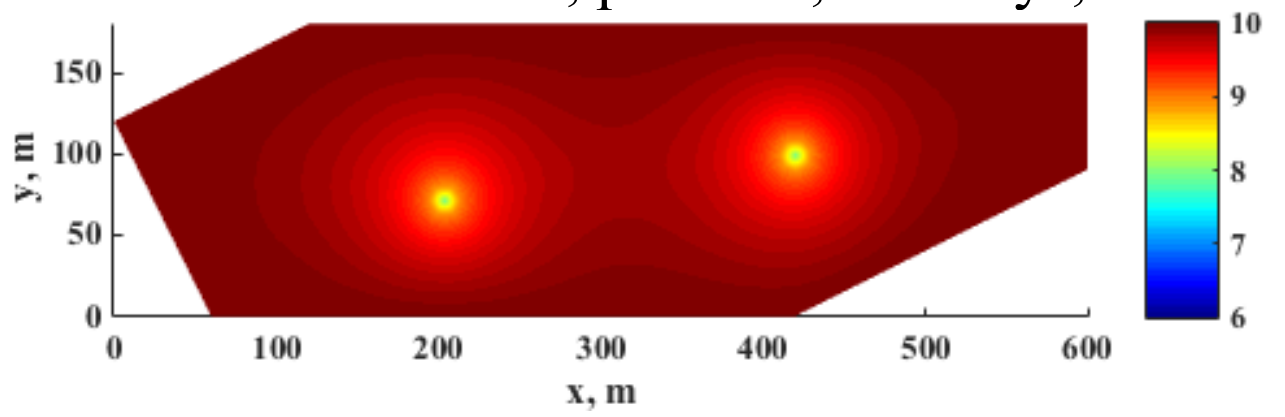

NCDMM solution, pressure, 300 days, 10m

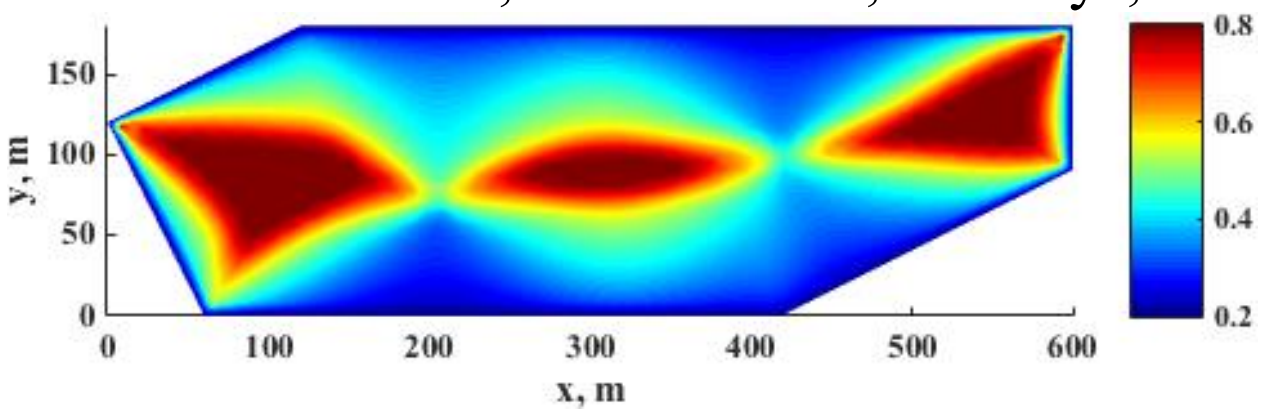

NCDMM solution, oil saturation, 300 days, 10m

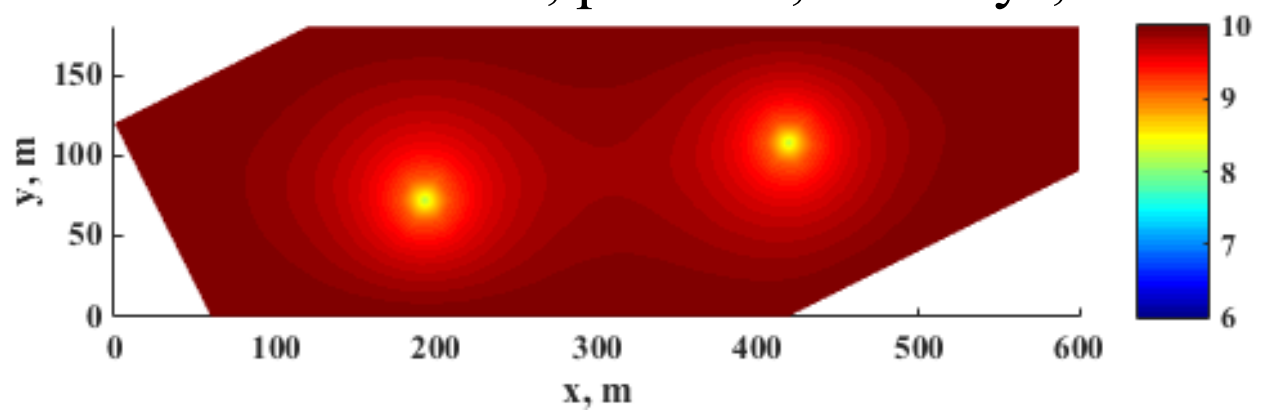

NCDMM solution, pressure, 300 days, 15m

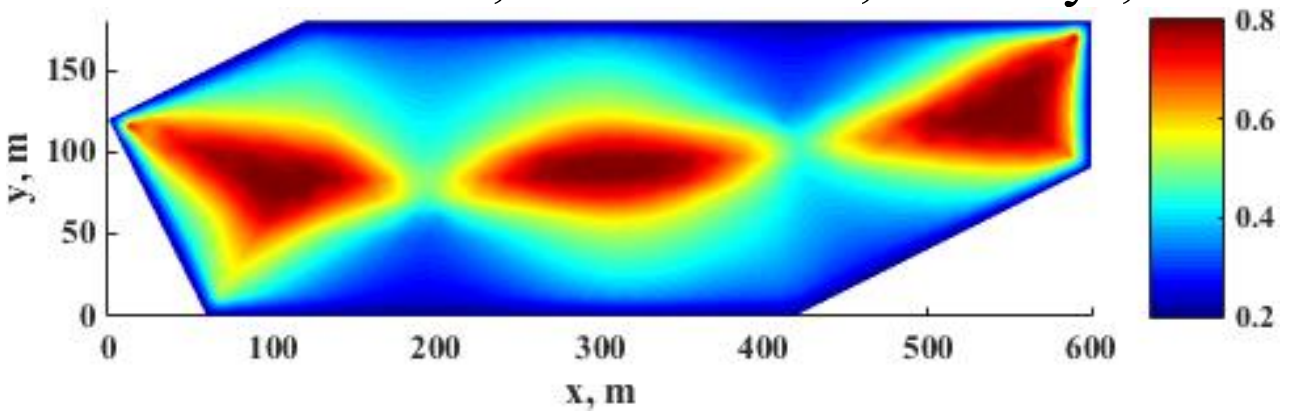

NCDMM solution, oil saturation, 300 days, 15m

Fig. 25 Calculation results of FVM-based reference solution and proposed NCDMM at 300 days

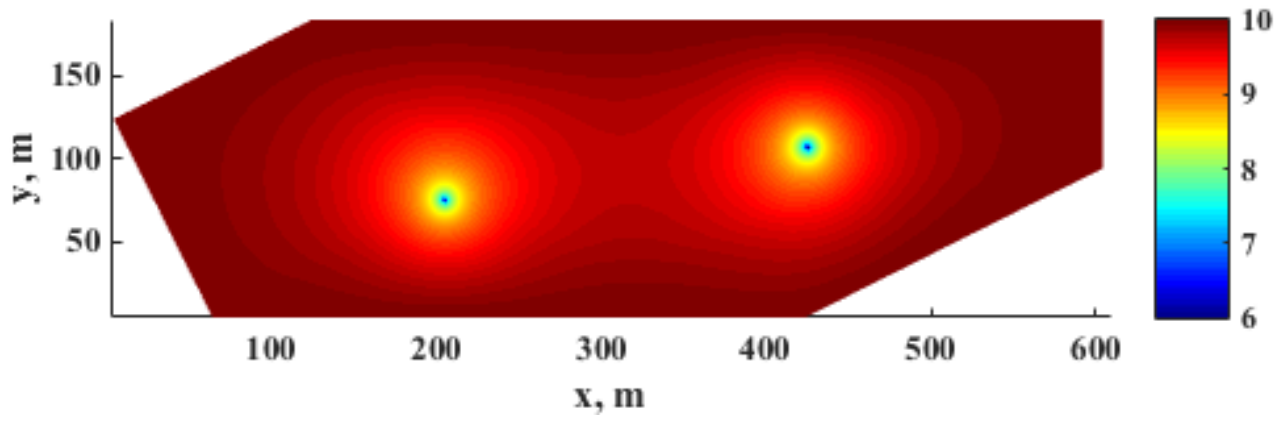

Reference solution, pressure, 100 days

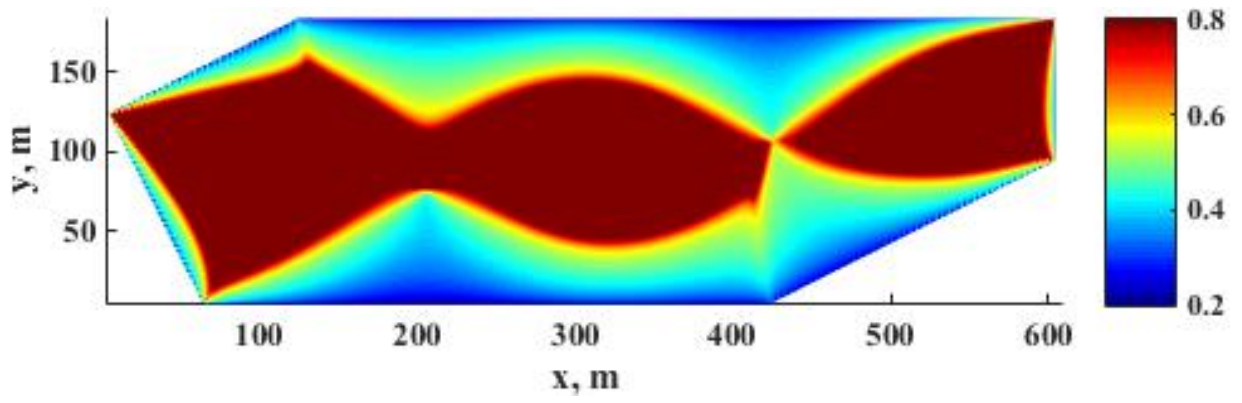

Reference solution, oil saturation, 100 days

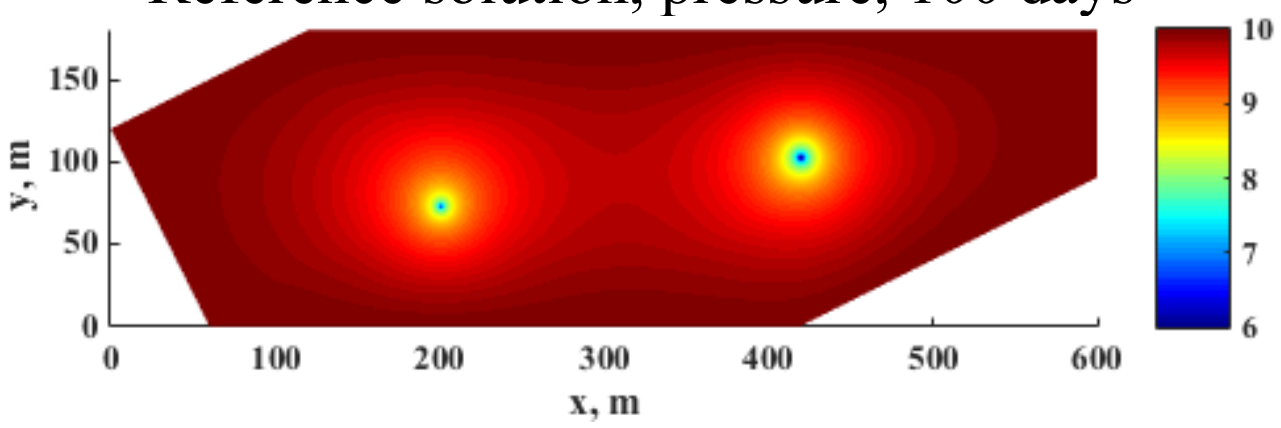

NCDMM solution, pressure, 100 days, 3m

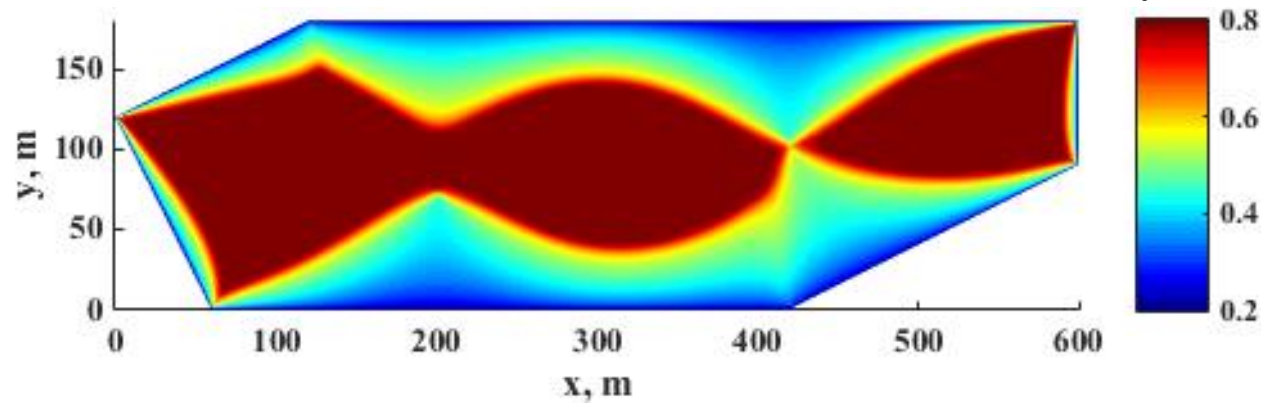

NCDMM solution, oil saturation, 100 days, 3m

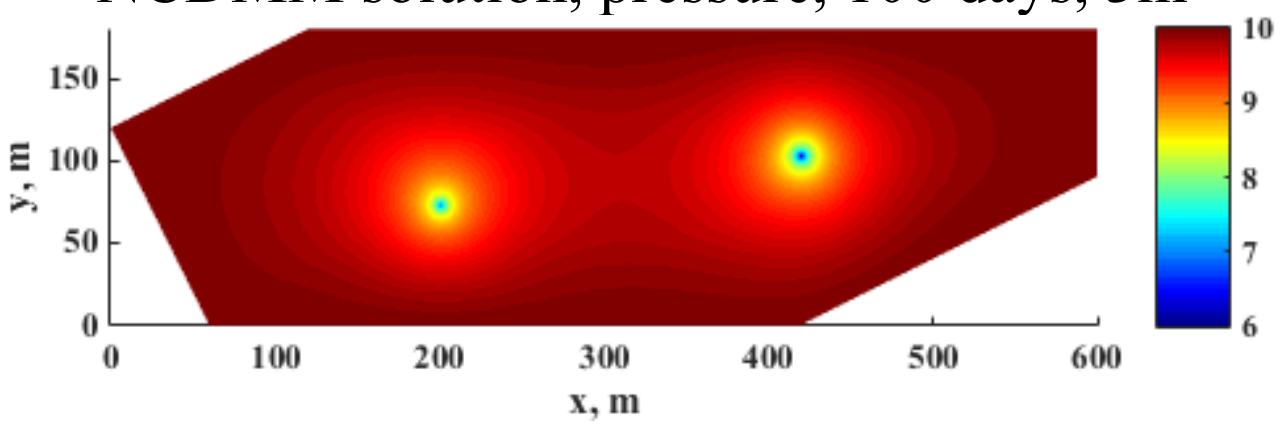

NCDMM solution, pressure, 100 days, 6m

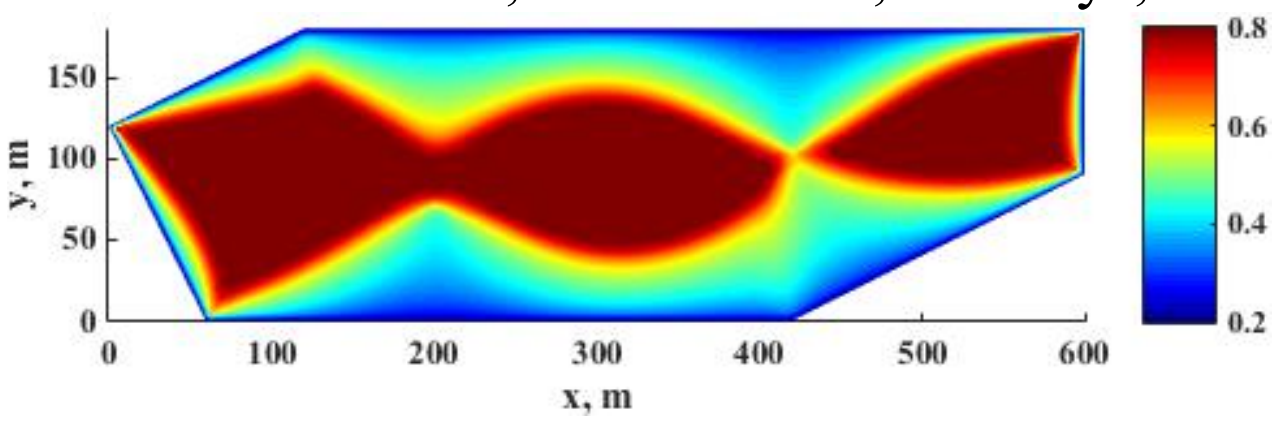

NCDMM solution, oil saturation, 100 days, 6m

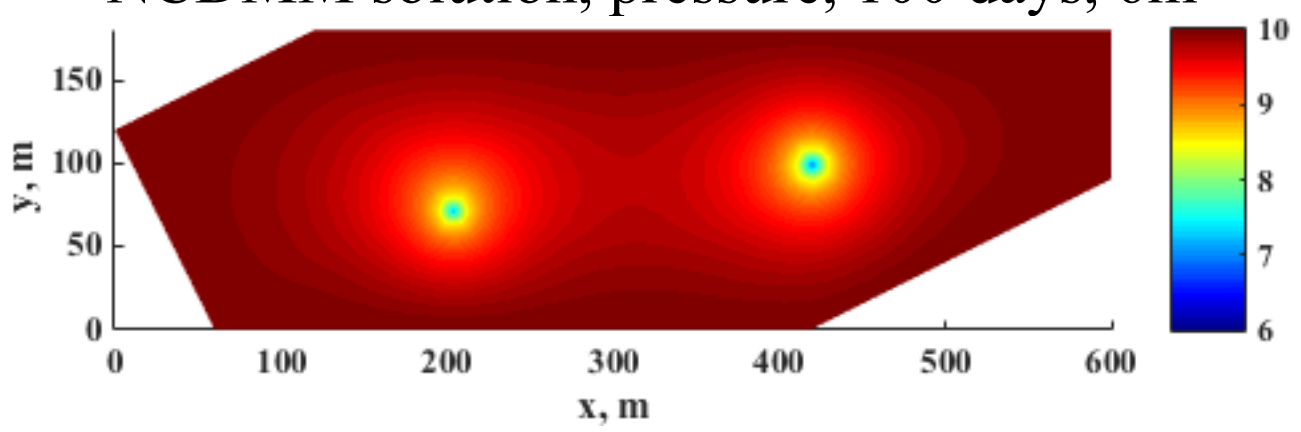

NCDMM solution, pressure, 100 days, 10m

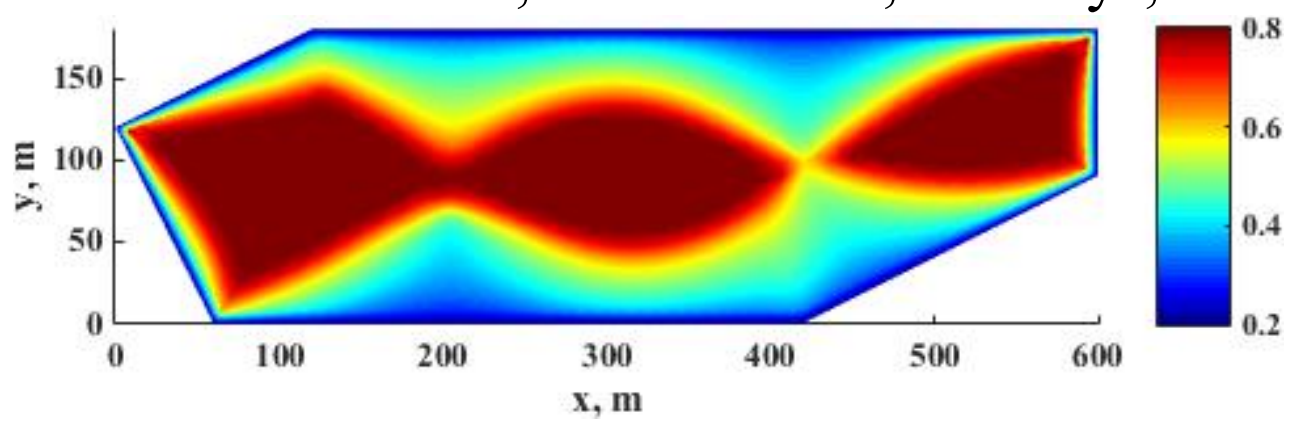

NCDMM solution, oil saturation, 100 days, 10m

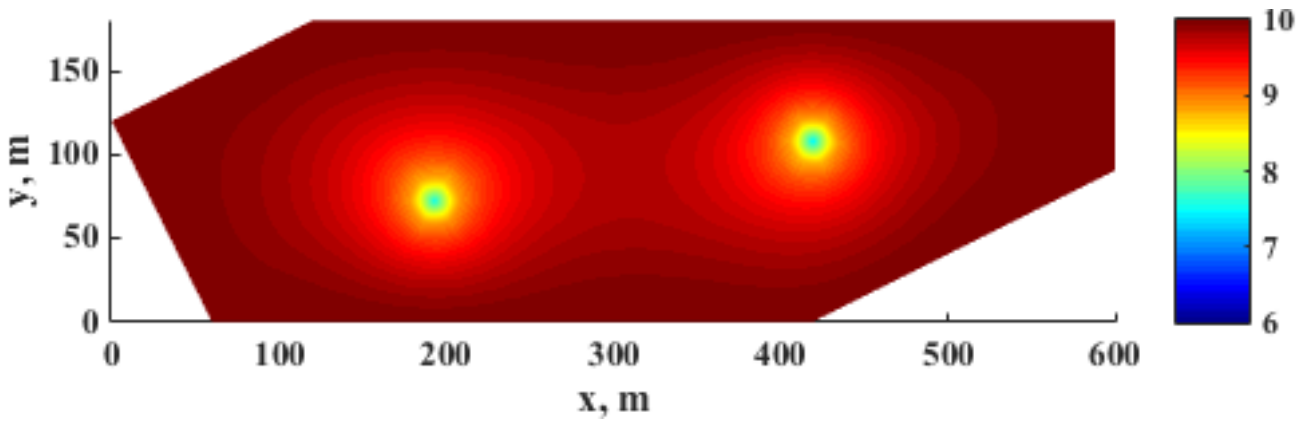

NCDMM solution, pressure, 100 days, 15m

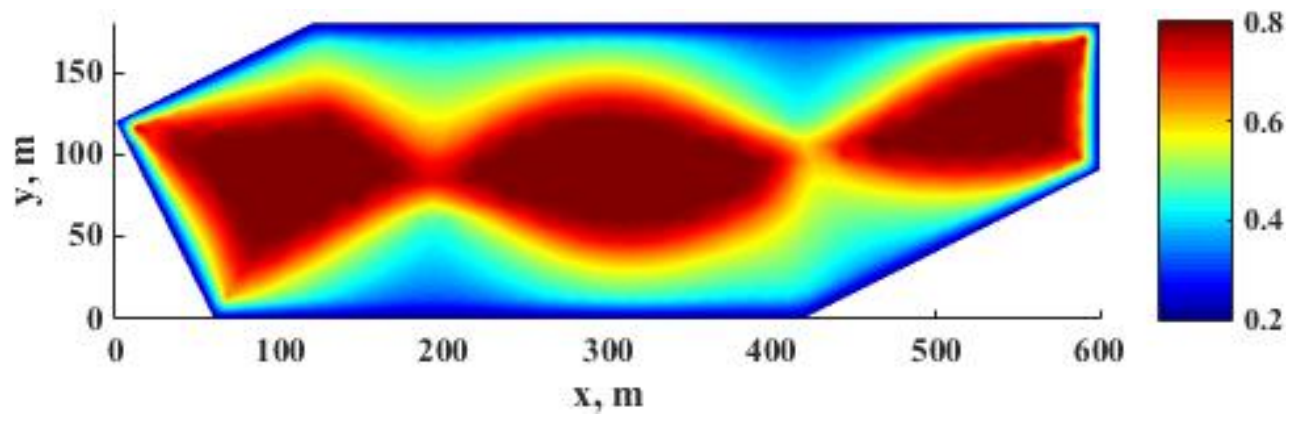

NCDMM solution, oil saturation, 100 days, 15m

Fig. 26 Calculation results of FVM-based reference solution and proposed NCDMM at 100 days

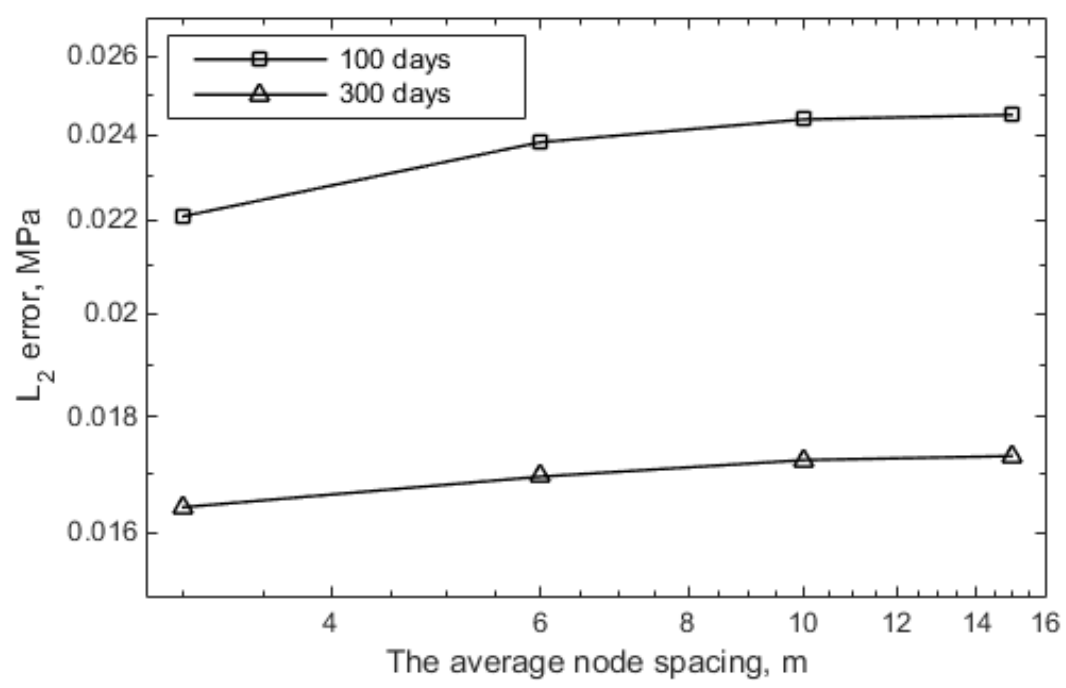

(a) pressure

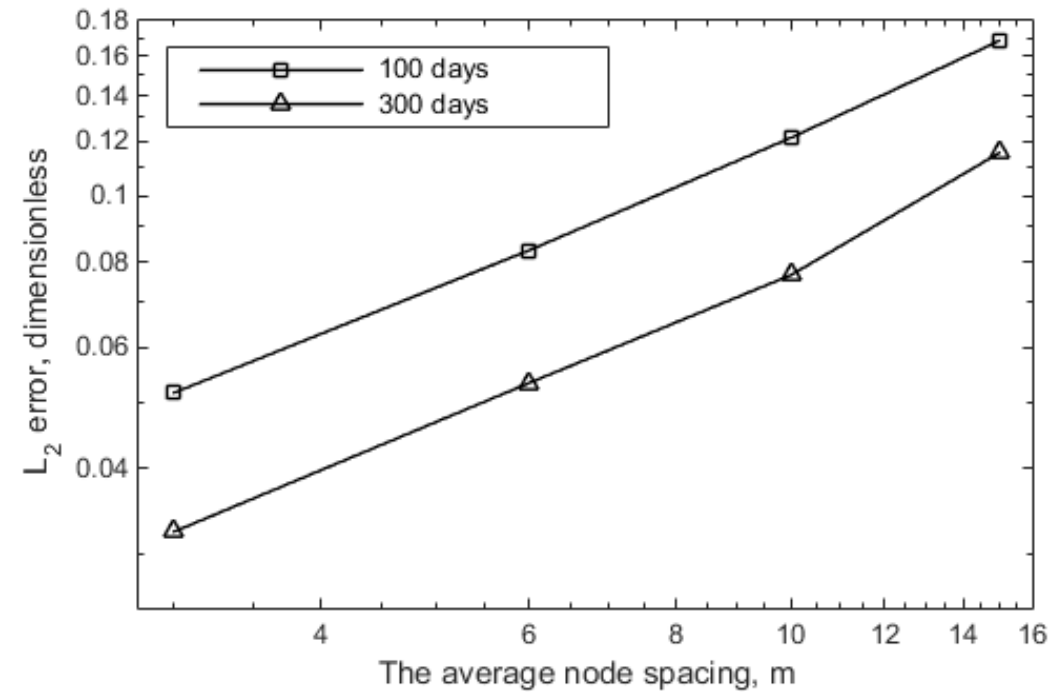

(b) water saturation

Fig. 27 $L_2$ errors of the calculated results

3.4 Two-phase porous flow in a reservoir domain with complex geometry

As shown in Fig. 28 (a), the boundary shape of the reservoir domain is a relatively tortuous curve, and the boundary is closed. There is a water injection well at coordinates (142.5, 52.5), with a constant injection rate of 50m$^3$/d, and a production well at coordinates (277.5,57.5), with a constant liquid production rate of 50m$^3$/d. To discretize the computational domain, the point cloud with added virtual nodes (marked in red) shown in Fig. 28 (b) is given in this example, and the point cloud is composed of the mesh vertices of the triangulation of the reservoir domain. Based on the point cloud, two methods are used to construct the *i-j* pairs (i.e., connectable point cloud), one is to make the radius of the influence domain of each node 7.5m, and the other is to use the triangulation-based method given in Section 2.4. For the convenience of narration, these two connectable point clouds are respectively denoted as $\mathbf{NP}_1$ and $\mathbf{NP}_2$. Similar to Section 3.2, as shown in Fig. 28 (c), a larger rectangular fine-Cartesian-mesh reservoir model that can cover the real reservoir domain is adopted, and the closed boundary condition is constructed by supposing the mesh permeability outside the real reservoir domain to 0 mD. Based on this model, the traditional FVM is used to provide the reference solution for this numerical example.

For these two point clouds, $w_1$, $w_2$, and weighted $w_2$ are respectively used to calculate the corresponding node control volume distributions shown in Fig. 29 (b-g). The mesh volume distribution of the dual PEBI mesh of the triangular mesh is shown in Fig. 29 (a). From the geometric view of the node control volume, the node control volume distribution in NCDMM should be close to the mesh volume distribution of the corresponding PEBI mesh. Therefore, it can be seen from Fig. 29, whether $\mathbf{NP}_1$ or $\mathbf{NP}_2$, with $w_1$ to $w_2$, and then to weighted $w_2$, the calculated node control volume distribution is closer to the mesh volume distribution of the PEBI mesh, which shows that the empirical method of calculating node control volume given in Section 2.3 can indeed achieve higher-accuracy node control volume, and these numerical examples given in Section 3 are also the source of implementation experience pointed out in Section 2.3. Moreover, when the connectable point cloud is $\mathbf{NP}_2$, under the same method, the calculated node control volume is closer to the grid volume distribution of the PEBI grid than when the connectable point cloud is $\mathbf{NP}_1$. Combined with the subsequent comparison of simulation results, we can know that the construction method of the local node point cloud directly affects the accuracy of calculated node control volume calculation, thus affecting the accuracy of simulation results.

Next, the physical parameters in Table 4 and two-phase relative permeability data in Table 5 are used for calculation. Fig. 30, 31, 32, and 33 compare the calculated pressure and oil saturation distribution at 150 days and 420 days in detail. Fig. 34 compares the $L_2$ error of pressure and oil saturation distributions in different cases and the dynamic bottom hole pressure (BHP) data of the injection well. It can be seen that, for calculated pressure profiles and BHP data, $w_1$ to $w_2$, and then to weighted $w_2$, the L2 error becomes smaller and smaller, and when adopting $w_1$, the L2 error is much larger. As also clearly seen from the colorful distribution drawn in Fig. 30 and Fig. 32 that the pressure distribution calculated by adopting $w_1$ is significantly different from the reference solution. Moreover, the L2 error when the connectable point cloud is $\mathbf{NP}_2$ is smaller than that when the connectable point cloud is $\mathbf{NP}_1$. For such an injection-production system, when the net flow into a node control volume is determined, if the node control volume is wrongly enlarged, the pressure change at the node will be less than the actual pressure change. On the contrary, if the node control volume is wrongly reduced, the pressure change at the node will be greater than the actual pressure change. Therefore, the accuracy of pressure distribution directly reflects the calculation accuracy of node control volumes (of course also the accuracy of simulation results). Therefore, the comparison of the above calculation results is consistent with the comparison of the node control volume profiles. NCDMM with adopting the empirical calculation method of node control volumes in Section 2.3 can achieve good calculation performance, and the performance of NCDMM depends on the method of determining the connectable point cloud. So developing a generic method to determine the connectable point cloud that can achieve good calculation performance in the case of any point cloud is an important future work.

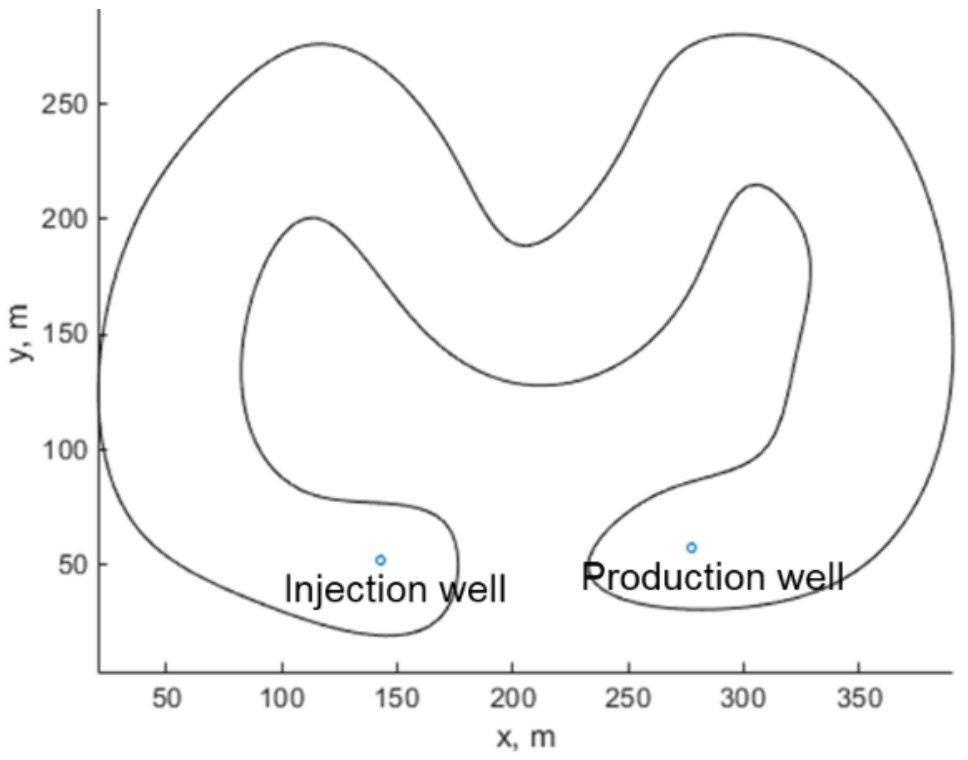

(a) the reservoir domain

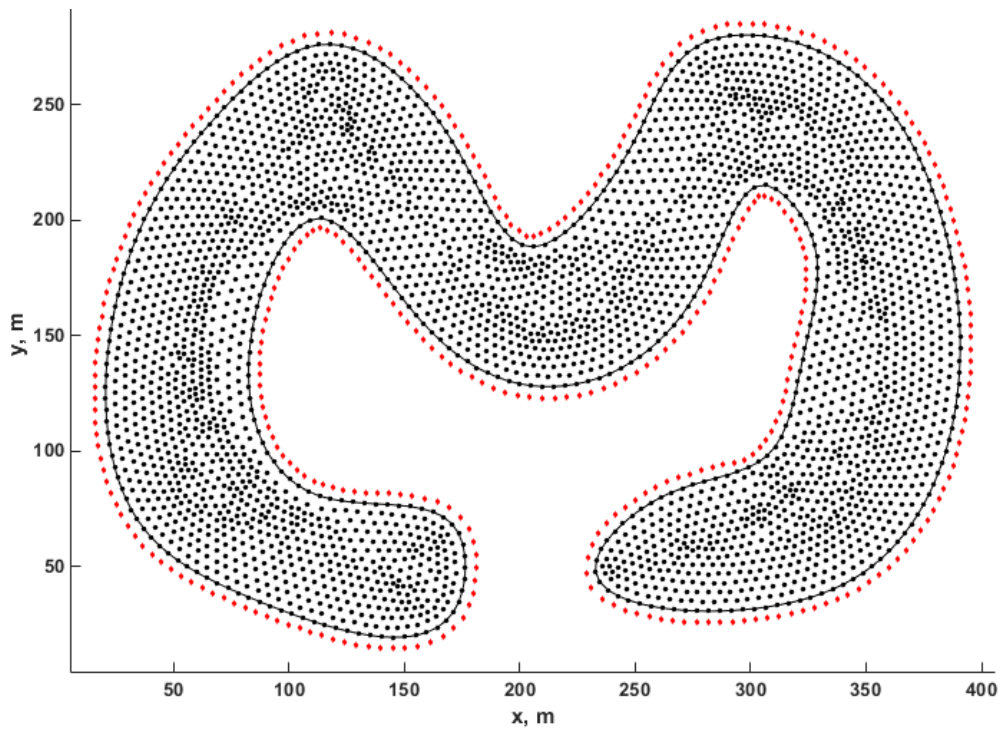

(b) the point cloud

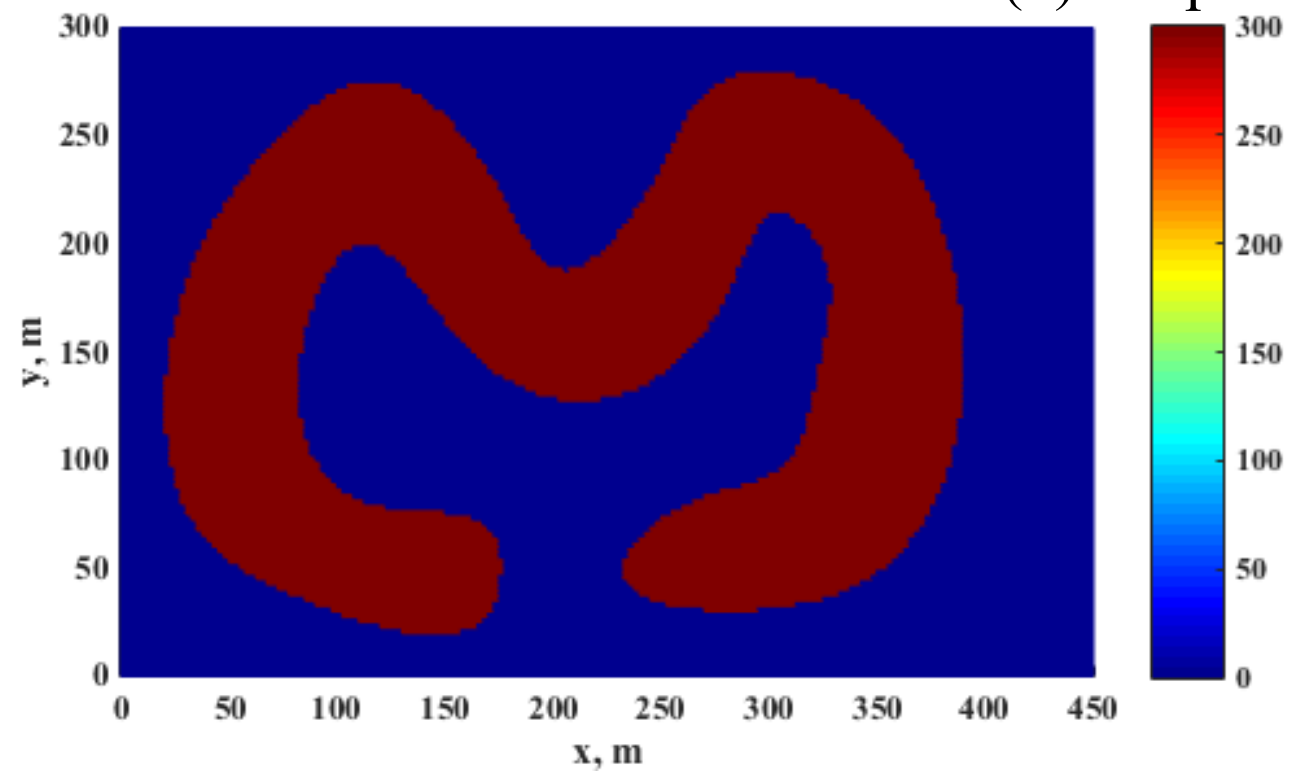

(c) the permeability profile of the rectangular model for reference solution

Fig. 28 The reservoir model and used point cloud in this example

Table 4 Basic physical properties used in this numerical example

| Properties | Values | Properties | Values |
|---|---|---|---|
| Porosity | 0.2 | Permeability | 300 mD |
| Oil compressibility | $8\times10^{-4}$ MPa$^{-1}$ | Water compressibility | $4\times10^{-4}$ MPa$^{-1}$ |
| Rock compressibility | $2\times10^{-4}$ MPa$^{-1}$ | Oil viscosity | 2 mPa·s |
| Water viscosity | 0.6 mPa·s | Initial reservoir pressure | 15 MPa |
| Initial water saturation | 0.20 | Oil volume factor | 1.0 |
| Water volume factor | 1.0 | Reservoir thickness | 5 m |
| Well radius | 0.1 m | Skin factor | 0 |

Table 5 Basic physical properties used in this numerical example

| *Sw* | *Krw* | *Kro* | *Sw* | *Krw* | *Kro* |
|---|---|---|---|---|---|
| 0.15 | 0 | 0.8 | 0.6276 | 0.207 | 0.063 |
| 0.45 | 0.04 | 0.4 | 0.6489 | 0.247 | 0.046 |
| 0.5212 | 0.08 | 0.227 | 0.6701 | 0.292 | 0.031 |
| 0.5425 | 0.103 | 0.184 | 0.6914 | 0.335 | 0.021 |
| 0.5638 | 0.125 | 0.156 | 0.7127 | 0.378 | 0.011 |
| 0.5851 | 0.145 | 0.125 | 0.76 | 0.475 | 0.005 |
| 0.6064 | 0.174 | 0.094 | 0.85 | 0.552 | 0 |

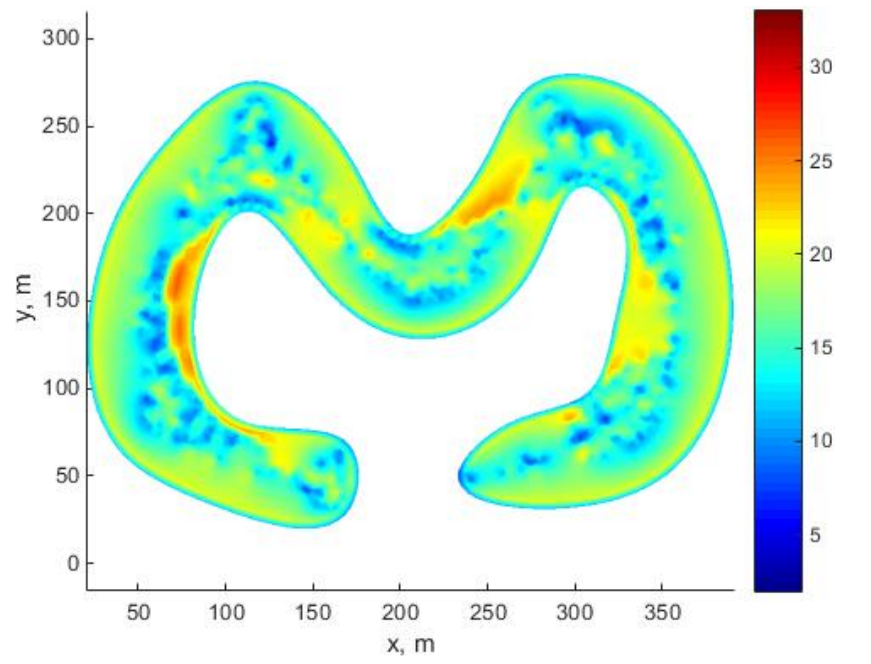

(a) dual PEBI grid volume

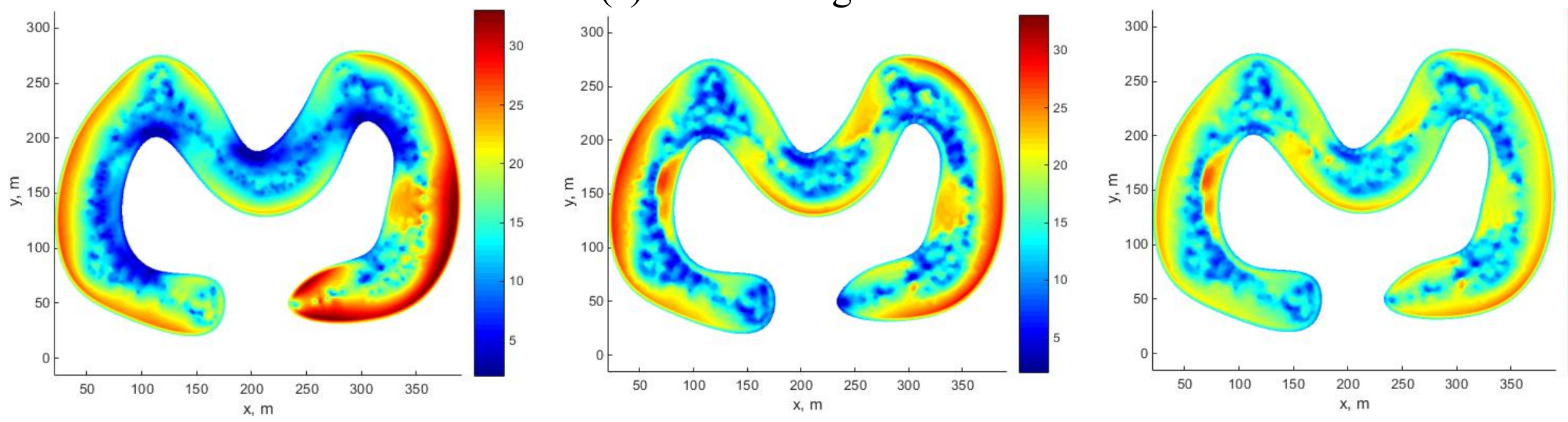

(b) $w_1$, $\mathbf{NP}_1$ (c) $w_2$, $\mathbf{NP}_1$ (d) weighted $w_2$, $\mathbf{NP}_1$

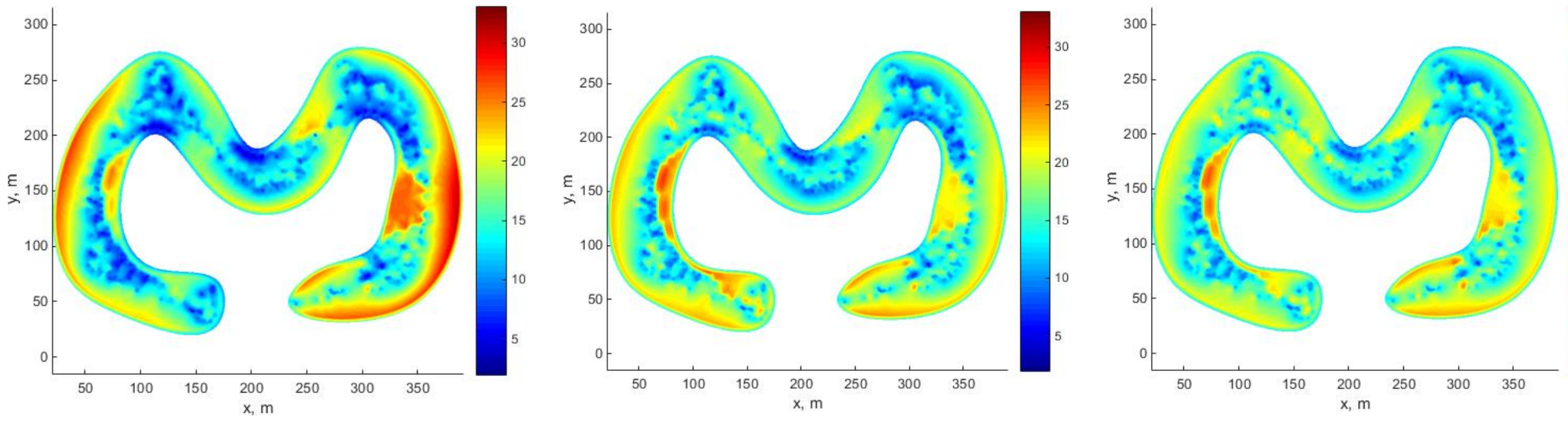

(e) $w_1$, $\mathbf{NP}_2$ (f) $w_2$, $\mathbf{NP}_2$ (g) weighted $w_2$, $\mathbf{NP}_2$

Fig. 29 The calculated node control volumes of the two node clouds by using different weight function

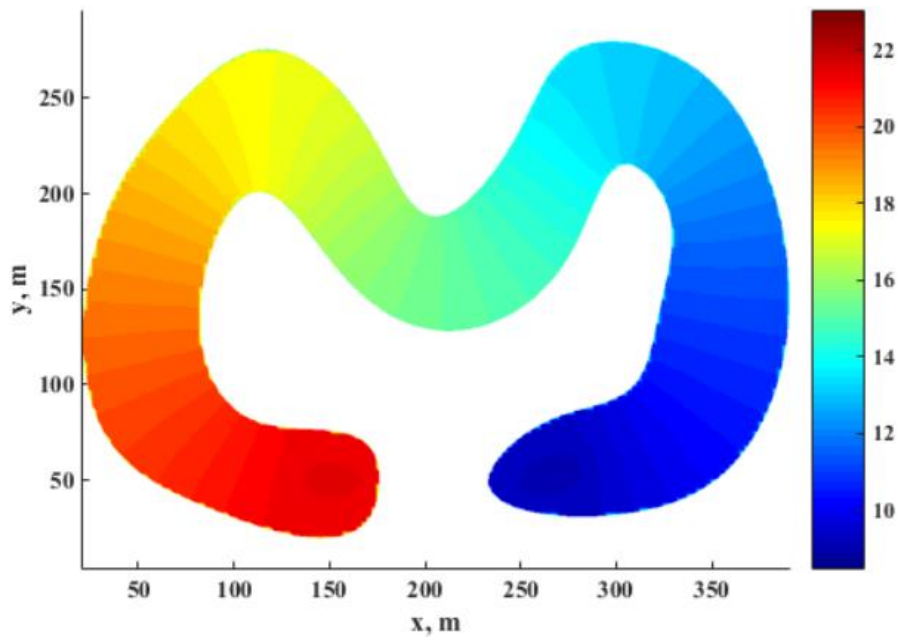

(a) reference solution

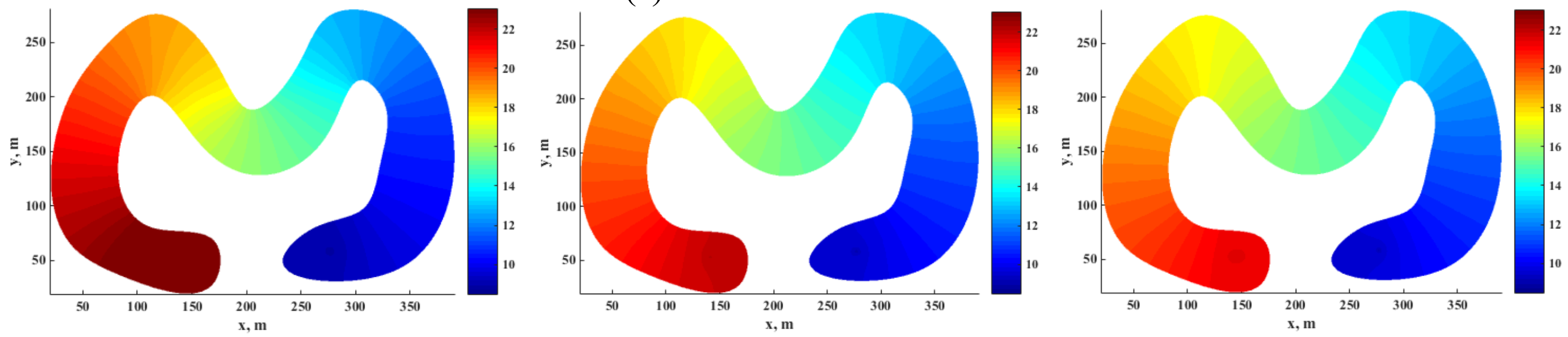

(b) $w_1$, $\mathbf{NP}_1$ (c) $w_2$, $\mathbf{NP}_1$ (d) weighted $w_2$, $\mathbf{NP}_1$

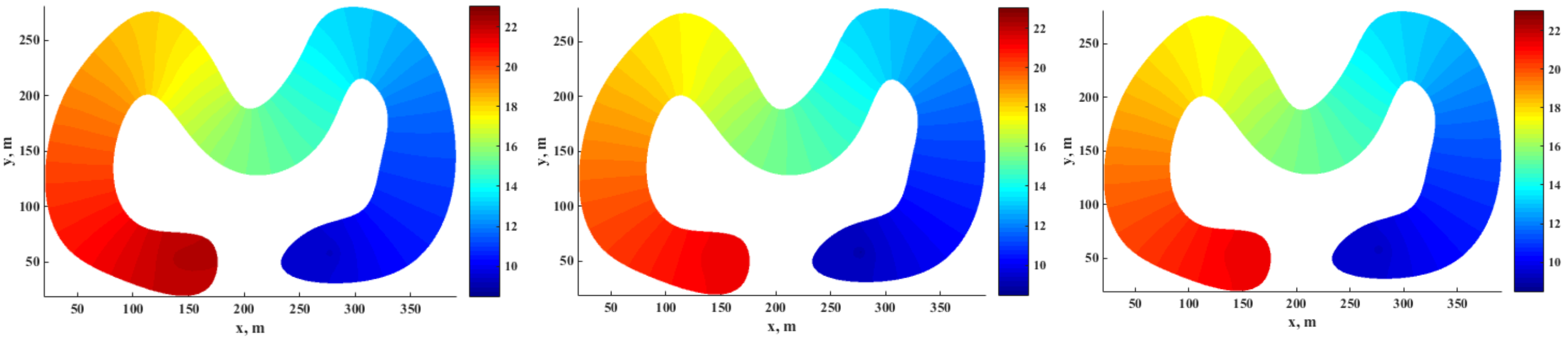

(e) $w_1$, $\mathbf{NP}_2$ (f) $w_2$, $\mathbf{NP}_2$ (g) weighted $w_2$, $\mathbf{NP}_2$

Fig. 30 Comparisons of calculated pressure profiles at 150 days

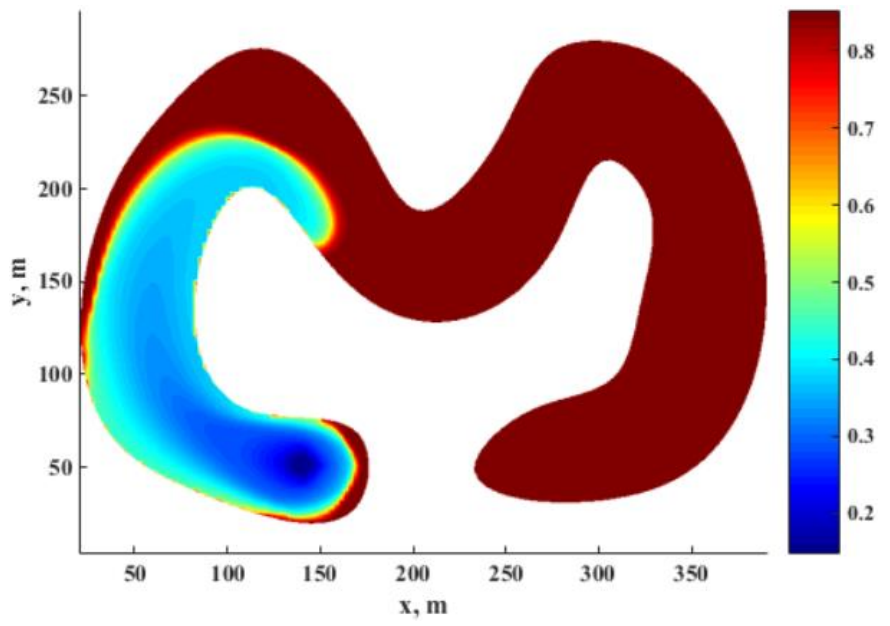

(a) reference solution

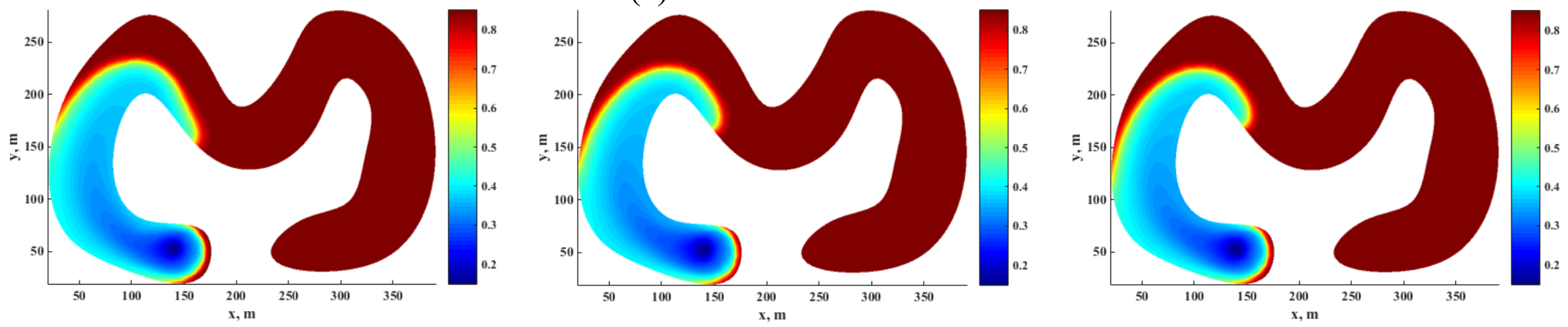

(b) $w_1$, $\mathbf{NP}_1$ (c) $w_2$, $\mathbf{NP}_1$ (d) weighted $w_2$, $\mathbf{NP}_1$

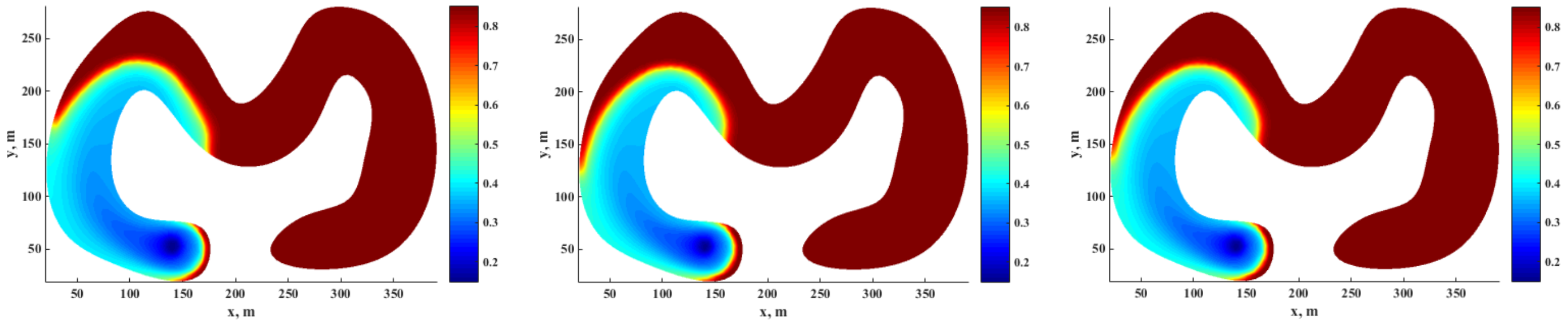

(e) $w_1$, $\mathbf{NP}_2$ (f) $w_2$, $\mathbf{NP}_2$ (g) weighted $w_2$, $\mathbf{NP}_2$

Fig. 31 Comparisons of calculated oil saturation profiles at 150 days

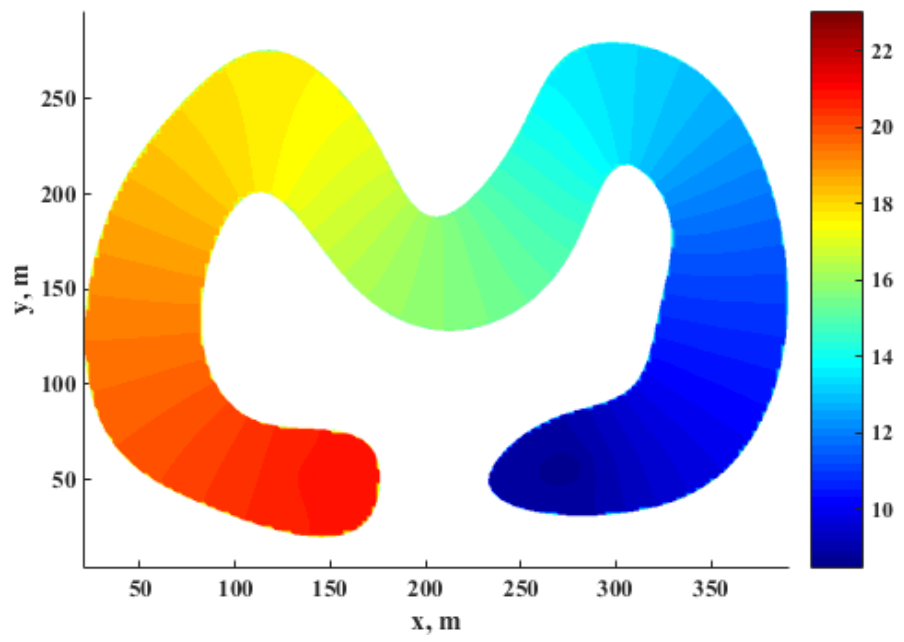

(a) reference solution

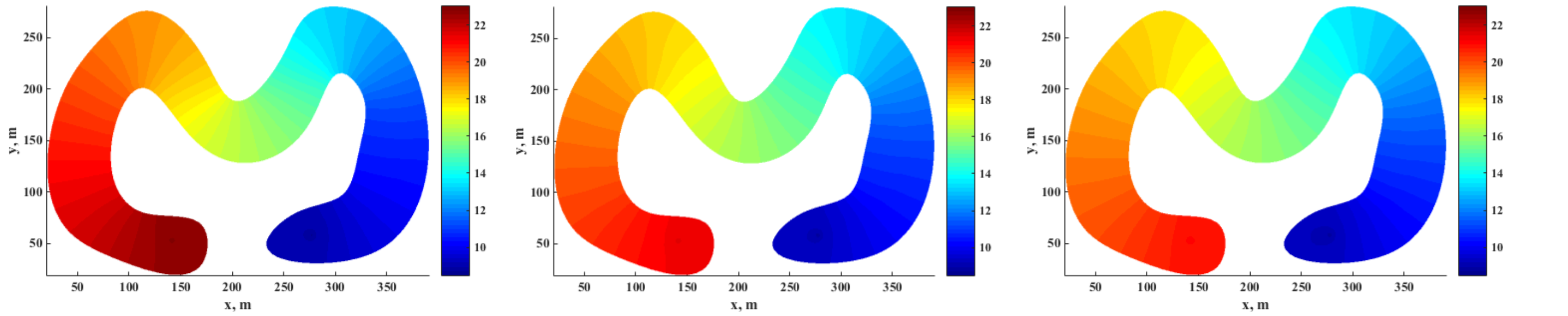

(b) $w_1$, $\mathbf{NP}_1$ (c) $w_2$, $\mathbf{NP}_1$ (d) weighted $w_2$, $\mathbf{NP}_1$

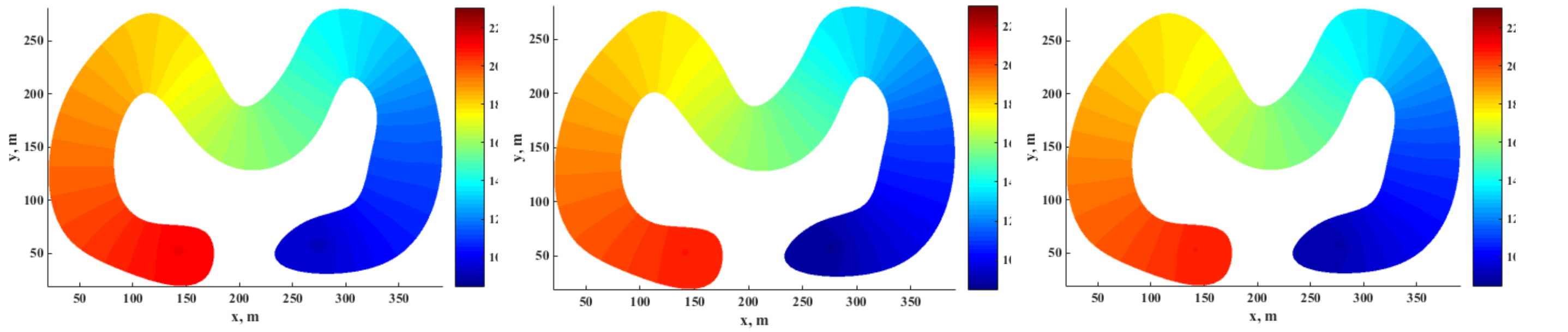

(e) $w_1$, $\mathbf{NP}_2$ (f) $w_2$, $\mathbf{NP}_2$ (g) weighted $w_2$, $\mathbf{NP}_2$

Fig. 32 Comparisons of calculated pressure profiles at 420 days

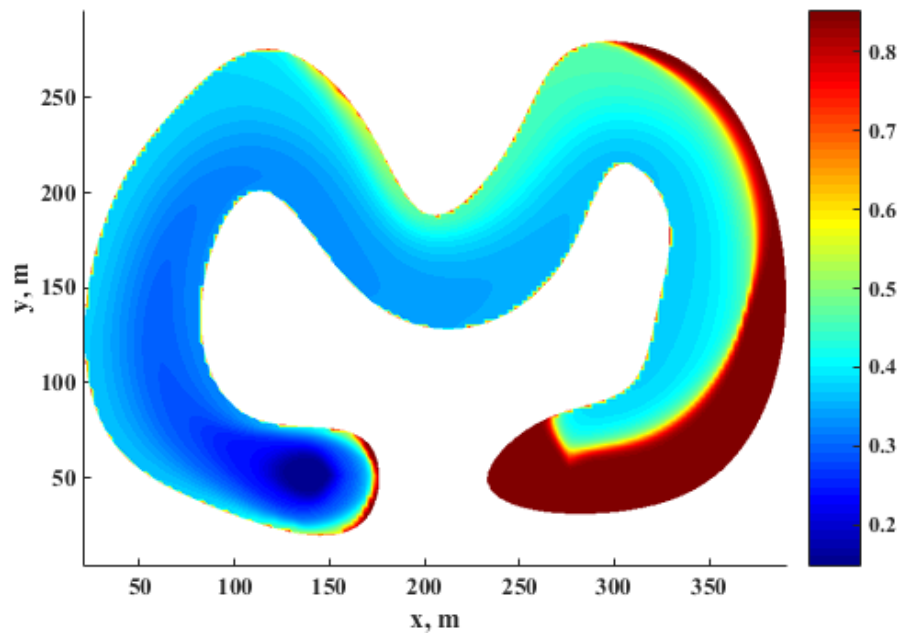

(a) reference solution

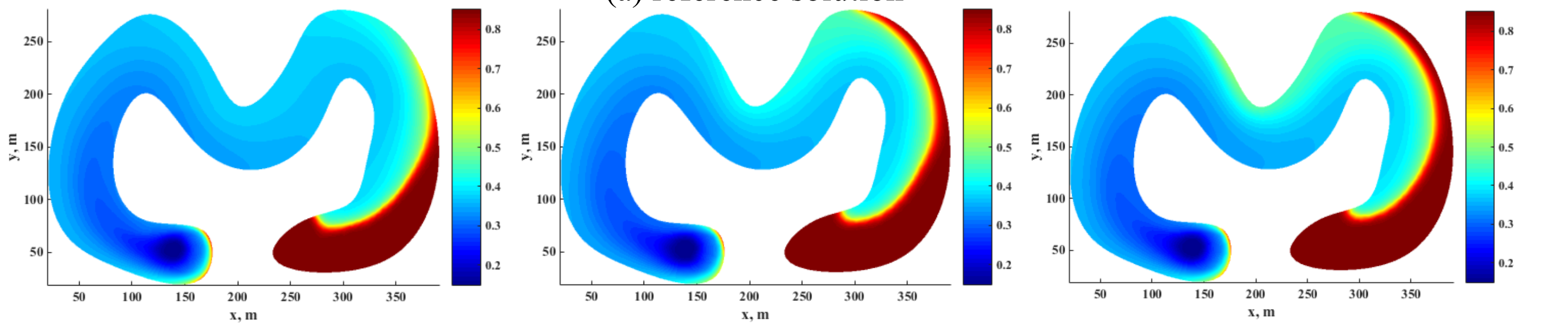

(b) $w_1$, $\mathbf{NP}_1$ (c) $w_2$, $\mathbf{NP}_1$ (d) weighted $w_2$, $\mathbf{NP}_1$

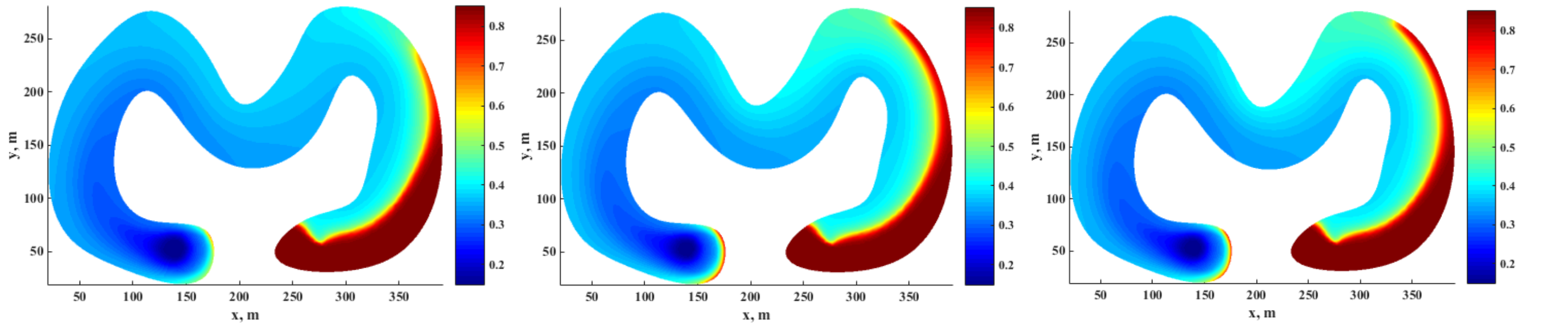

(e) $w_1$, $\mathbf{NP}_2$ (f) $w_2$, $\mathbf{NP}_2$ (g) weighted $w_2$, $\mathbf{NP}_2$

Fig. 33 Comparisons of calculated oil saturation profiles at 420 days

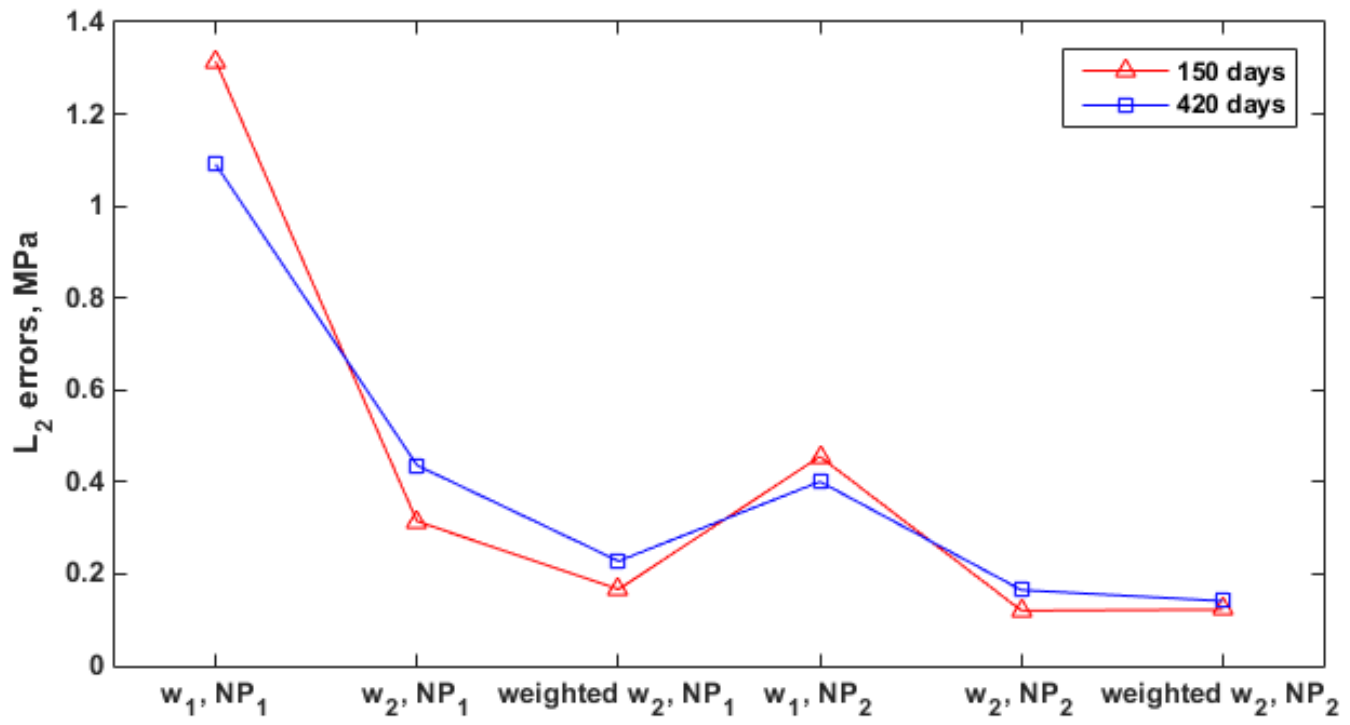

(a) $L_2$ errors of pressure

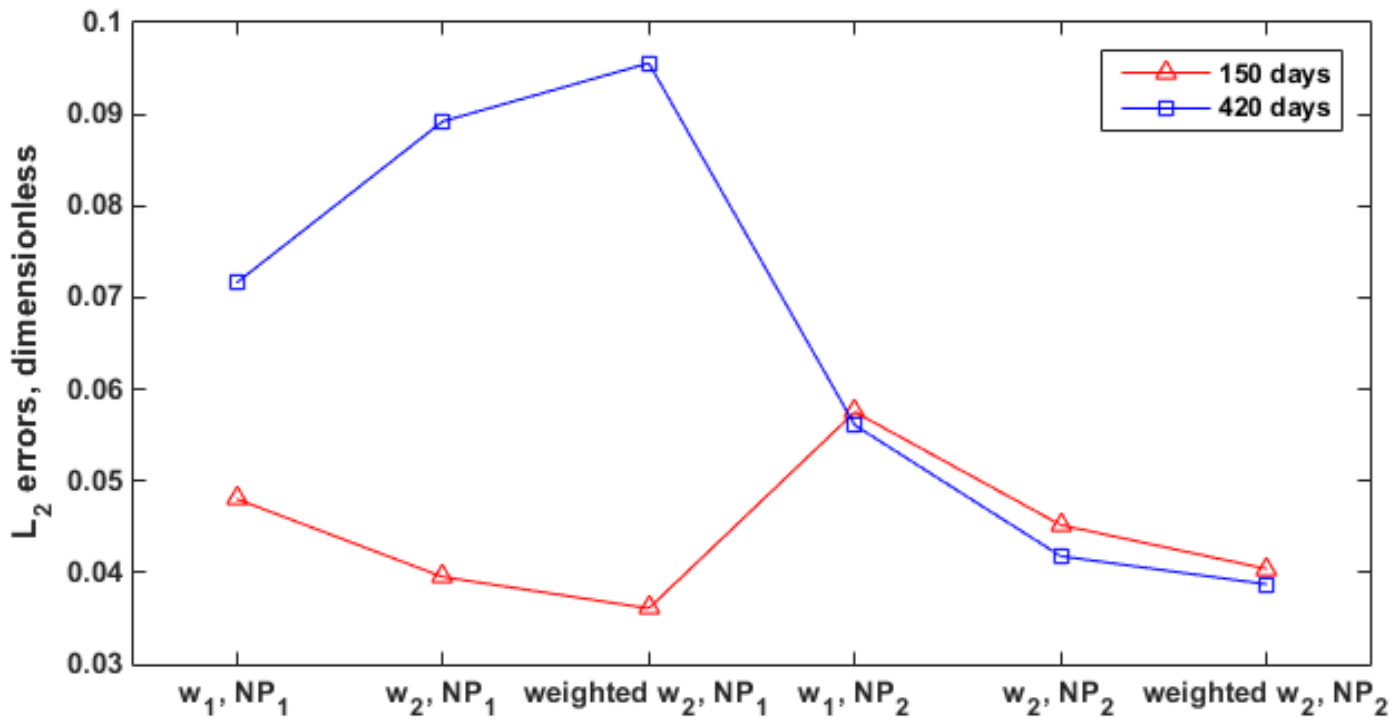

(b) $L_2$ errors of oil saturation

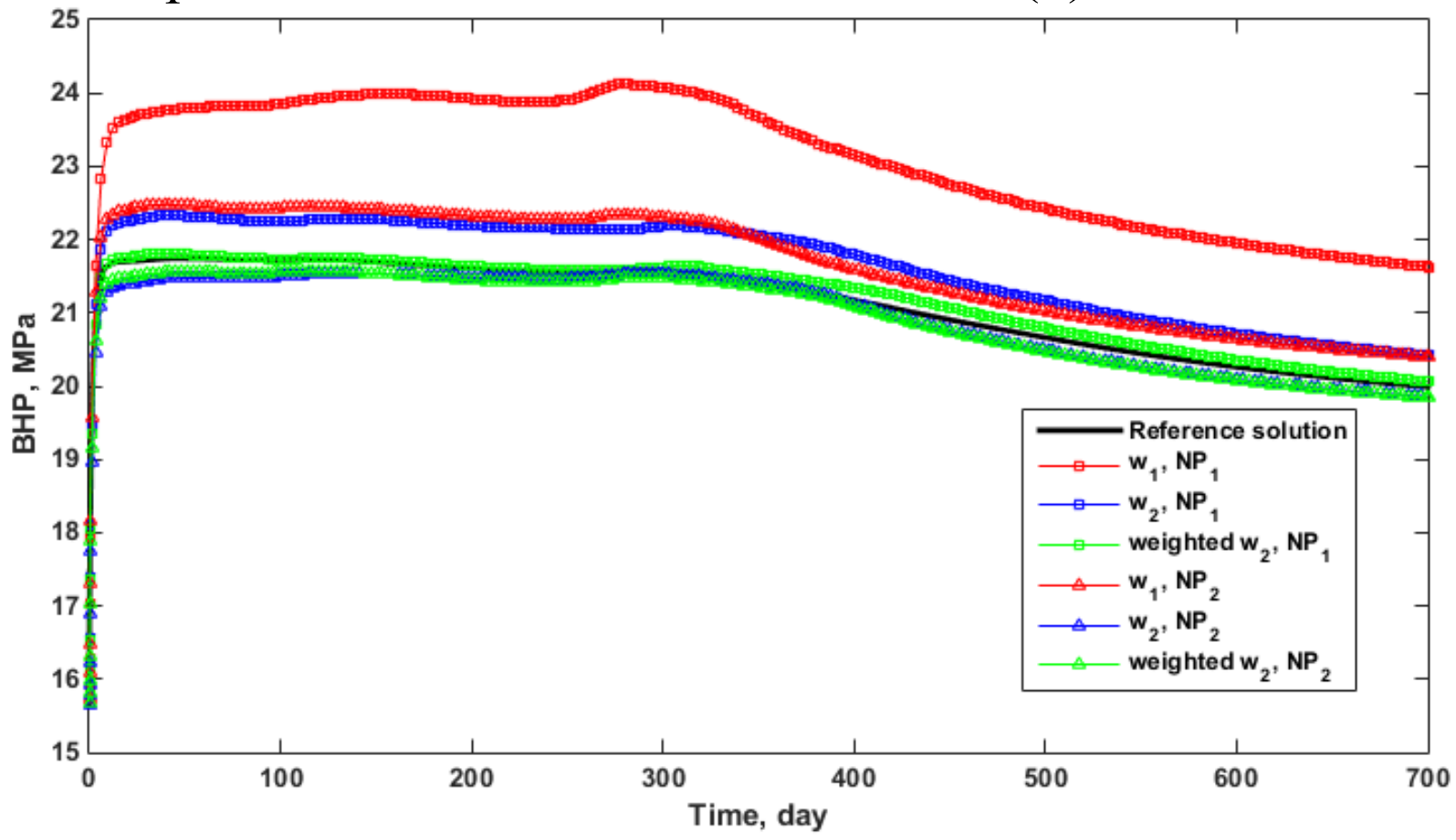

(c) calculated BHPs

Fig. 34 Comparisons of L2 errors of simulation results and BHP data in different cases

## 4. Conclusions

A novel meshless method based on the virtual construction of node control domains is developed, named NCDMM. Throughout the whole paper, seven main conclusions can be obtained as follows:

(1) Concepts of node control volume and connectable point cloud are defined in meshless GFDM, and the defined node control volumes are calculated by derived overdetermined linear equations for the connectable point cloud of the computational domain. The virtual construction of node control domains means that NCDMM only cares about the volume of the node control domain, instead of the geometric shape of the node control domain.

(2) The NCDMM discretization scheme of the governing equations is obtained by integrating the GFDM discretization scheme is integrated on the node control volume, which satisfies the local mass conservation.

(3) The addition of virtual nodes is required for the treatment of derivative boundary conditions and the calculation of node control volumes.

(4) An empirical method of calculating node control volumes is proposed. Compared with using classical quartic spline function, it can yield much higher-accuracy nodes-control-volume distribution, to obtain higher-accuracy simulation results.

(5) Different connectable point clouds constructed based on the same point cloud of the calculation domain have significant differences in the calculation accuracy of node control volumes and simulation results, which shows the importance of an excellent generation method of connectable point cloud for NCDMM to achieve good computational performance. This work presents a connectable point cloud construction method that can achieve good performance in the case of point clouds composed of mesh vertices of the triangulation of the calculation domain.

(6) The solution of the NCDMM-based discrete equations can directly employ the existing nonlinear solver in the FVM-based reservoir simulator, which significantly reduces the cost to construct an NCDMM-based general-purpose reservoir numerical simulator.

(7) Four numerical examples are implemented to show the computational performances of the proposed NCDMM, and the results verify the above theoretical advantages of NCDMM in directly dealing with complex geometries and various boundary conditions.

## 5. Future work

From our viewpoints, there may be three valuable future works:

(1) Although we believe that there is not much technical difficulty in extending the NCDMM to 3D, it is still an interesting topic to extend the NCDMM to 3D and study its computational performances.

(2) Considering the advantages of the NCDMM, the porous flow in the fractured reservoirs and fracture-vug-type reservoirs which have complex geometry and reservoirs with complex boundary conditions may be the ideal potential application scenarios of NCDMM.

(3) NCDMM is expected to be applied to more scientific and engineering fields.

## 6. Acknowledgements

Dr. Rao thanks the support from the National Natural Science Foundation of China (No. 52104017), the Open Fund of Cooperative Innovation Center of Unconventional Oil and Gas (Ministry of Education & Hubei Province) (No. UOG2022-14), and the Open Fund of Hubei Key Laboratory of Drilling and Production Engineering for Oil and Gas (Yangtze University) (Grant No. YQZC202201)